\documentclass[a4paper,11pt]{article}
\usepackage{color,xcolor,ucs}
\usepackage[top=0.3in, bottom=0.5in, left = 0.65in, right = 0.65in]{geometry}
\usepackage{mathtools}

\usepackage{color,xcolor,ucs}
\usepackage{mathtools}   \usepackage{tikz} 
\usepackage{ amssymb }
\usepackage{extarrows} 
\usepackage{pgf,tikz}
\usepackage{float}
\usetikzlibrary{positioning}
\usetikzlibrary{shapes.geometric}
\usetikzlibrary{shapes.misc}
\usetikzlibrary{arrows}
\usepackage{caption}
\usepackage{mathrsfs}
\usetikzlibrary{arrows,shapes,automata,backgrounds,petri,positioning}
\usetikzlibrary{decorations.pathmorphing}
\usetikzlibrary{decorations.shapes}
\usetikzlibrary{decorations.text}
\usetikzlibrary{decorations.fractals}
\usetikzlibrary{decorations.footprints}
\usetikzlibrary{shadows}
\usetikzlibrary{calc}
\usetikzlibrary{spy}
\usepackage{amsmath}
\usepackage{array}
\usepackage{ amssymb }
\usepackage{braket}
\usepackage{qcircuit}
\usepackage{soul}
\usepackage{braket} 
\usepackage{relsize}

\usepackage{amsmath}
\usepackage{ amssymb }
\usepackage{braket}
\usepackage{qcircuit}
\usepackage{soul}
\usepackage{braket}

\title{{\LARGE Operator formalism for discretely holomorphic parafermions of the two-color Ashkin-Teller, loop $\mathrm{O} ( 1 )$, staggered eight-vertex, odd eight-vertex, and abelian sandpile models}}
\author{P. Rigas}
\date{}

\begin{document}

\maketitle

\begin{abstract}
  We extend an operator formalism, developed by  Hongler, Kytölä, and Zahabi in 2012 for the Ising model, to the Ashkin-Teller and loop $\mathrm{O} ( n ) $ models. The formalism is primarily dependent upon notions of massive, and massless, s-holomorphicity of the Ising model, which are respectively satisfied at, and below, the critical temperature, hence allowing for rigorous analyses of the transfer matrix, correlation functions, and lattice fermion operators. From results of another paper, by Tanhayi-Ahari and Rouhani in 2012, which demonstrates that parafermionic observables for the staggered eight-vertex model coincide with those of the Ashkin-Teller model, as well as a 2003 paper, by Wu and Kunz, which establishes a correspondence between the staggered eight-vertex and odd eight-vertex models, we establish further associations. \footnote{\textit{{Keywords}}: Ising model, random-cluster model, parafermionic observable, Cauchy-Rieman equations, massless/massive s-holomorphicity} \footnote{\textbf{MSC Class}: 82B20; 82B44; 82B27; 82B30; 34L25; 60K35}
\end{abstract}

\section{Introduction}

\subsection{Overview}

\noindent Disorder parameters, and variables, which share connections with fermionic and parafermionic observables [2,6], have attracted great interest, with studies examining topics ranging from divergence of the correlation length for FK percolation [6], behavior away from criticality [2], connective constant of the honeycomb lattice [12], operator formalism for the Ising model [13], near-criticality of the FK-Ising model [8], critical temperature [4], connection probabilities [10], as well as several other potential applications to other models, ranging from sharp thresholds [7], self-avoiding walks [9], and conformal invariance [11]. In [6], Duminil-Copin and colleagues raised several open questions, one of which pertains to whether information from parafermionic observables, and related objects, can be extracted from other models, from arguments demonstrating that the correlation length for planar FK percolation diverges. Within the context of the aforementioned works in Statistical Physics, there continues to exist strong motivation for investigating disordered systems, whether ones with continuous, or discontinuous, phase transitions. As such, probabilistic tools that have garnered significant traction, which share connections with Complex Analysis and discrete holomorphicity, can be leveraged for providing an explicit constant for the connective constant of the honeycomb lattice, [12], conformal invariance, [11], long-range order [14], and parafermions for the 8-vertex model by establishing a correspondence with the 6-vertex model [15,16]. Despite the fact that there have already been extensive characterizations of disorder within systems of Statistical Physics, it not only continues to remain of interest for determining additional models for which discrete holomorphicity, and related notions, hold, but also for studying near critical behavior of several models simultaneously. To this end, from a previous work on discrete holomorphicity for the Ising model, [15], we further explore notions that the authors develop, which are useful for capturing critical, and near critical, behaviors, including computations with the transfer matrix, the Cauchy-Riemann equations, and several related objects. 

In other models of interest from Statistical Physics, the Cauchy-Riemann equations, which can be thought of as localized relations on the finite lattice over which correlations, and other interactions, occur, take on different forms than those presented for the Ising model in the aforementioned work. Nevertheless, by leveraging other structures of the equations, [14], which in the case of the Ashkin-Teller, and closely related models that can be obtained by coupling Potts, or other types, of models together, one can leverage operator notions for characterizing low temperature expansions of parafermionic observables. Moreover, the fact that the Cauchy-Riemann equations, and corresponding notions of discrete holomorphicity, for other models in Statistical Physics such as the Ashkin-Teller model, and various loop models,  provides additional avenues for future research, primarily in determining which universal behaviors are common to classes of models. Besides seminal applications of parafermionic observables from Duminil-Copin and Smirnov mentioned previously, further developing arguments for which the parafermionic observable plays a key role can potentially shed light upon the transfer matrix, and related objects, of the quantum inverse scattering framework. Within this inverse scattering framework, the transfer and quantum monodromy matrices, much like discrete holmorphicity and the parafermionic observable, expose the intricate nature of integrability and near critical phenomena.

To explore all such possible research directions, we proceed with the operator formalism, which as an extension of the results provided previously for the Ising model, [13], relies upon a characterization of massive, and massless, s-holomorphicity which is related to the Cauchy-Riemann equations which have been known to vanish for parafermionic observables, from connections with Morera's Theorem. Albeit being able to establish that massive, and massless, s-holomorphicity hold for the Cauchy-Riemann equations associated with parafermionic observables of the Ising model, establishing that similar properties hold for parafermionic, or fermionic, observables of other models depends upon: identifying critical points of the model, which in the case of the two-color Ashkin-Teller model, amount to a line of possible critical points, the smallest of which is $\frac{1}{4} \mathrm{log} \big( 3 \big)$; analyzing the propagation mechanism under the Cauchy-Riemann equations; and finally, making use of the transfer matrix for the two-color Ashkin-Teller model to perform computations with the Pfaffian. Besides being able to perform these computations, we make use of a completely different duality relation for the Ashkin-Teller than that of the Ising model, which informs several quantities that are obtained for the propagation mechanism. These steps not only hold promise for other models, in the sense that accompanying notions of massive, and massless, strong holomorphicity depend upon the sites of the lattice at which a suitable observable is evaluated, but also illustrate connections between parafermionic observables that are bosons and fermions, which was not previously discussed for the operator formalism of the Ising model. As a starting point, developing the operator formalism for the Ising model provides information on how notions of discrete holomorphicity, along with the accompanying observable, illustrate how correlations, the correlation length, and related quantities behave when approaching the critical point from above and below. For other models of Statistical Physics that are less classical than the Ising model, the operator formalism provided in [13] is capable of describing critical, and near critical, phenomena for other models given their respective phase diagrams and critical points. From the wide variety of models in this work which have parafermionic observables that are bosons and fermions, the operator formalism for the loop model, given $n \equiv 1$, is particularly interesting. The loop model, in comparison to the two-color Ashkin-Teller, and abelian sandpile, models, has connections to vertex models in the full packing limit which makes it different, from the other models examined in this work.

Extending the operator formalism to the two-color Ashkin-Teller model first allows for connections with loop models, in particular the loop $\mathrm{O} \big( n \big)$ model for $n \equiv 1$, in addition to vertex models, in which there exists a correspondence between the 8-vertex, under staggering, to the odd 8-vertex model [16]. Within the model operator framework, parafermionic observables that are bosons can be further examined from previous computations, [15], which demonstrated that holomorphic parafermions exist in the 8-vertex model. Demonstrating the existence of such observables for the 8-vertex model is primarily reliant on computations surrounding a parameter redefinition, hence allowing the authors to deduce the forms of observables in the staggered 6-vertex, Ashkin-Teller, and loop $\mathrm{O} \big( n \big)$ models. Establishing connections between vertex, and lattice, models has emerged as a research area of intense interest within Probability theory and Mathematical Physics, as exhibited in previous work of the author which established connections between Russo-Seymour-Welsh theory, and crossing probabilities, for the 6-vertex model under sloped boundary conditions that can be mapped onto the Ashkin-Teller model on the self dual line of the phase diagram. On more of the discrete holomorphicity and disorder variable side, the operator formalism proposed by Hongler, Kytola and Zahabi states that the parafermionic observable for the Ising model is a kernel of the convolution operator, which can also be studied in the case of the two-color Ashkin-Teller, and loop models. For other models, we make use of this framework to obtain the kernel of other convolution operators, within the context of the staggered eight-vertex and odd eight-vertex models [16].

\subsection{Random cluster model objects}

\noindent We provide a brief overview of properties of the random-cluster and Ising models from [2]. Fix $G \equiv \big( V , E \big)$. Over the support of this graph, denote the random cluster measure as,

\begin{align*}
  \phi^{\xi}_{p,q,G} \big( \omega \big) = \frac{p^{\mathrm{o}(\omega)} \big( 1 - p \big)^{\mathrm{c}(\omega)} q^{k ( \omega)}}{Z^{\xi} \big( p , q , G \big)}  \text{ } \text{ , } 
\end{align*}

\noindent in which the probability of sampling a \textit{random cluster configuration}, $\omega$, under $\phi$ with boundary conditions $\xi$ is dependent upon the product of the $p$, $1-p$, and the number of clusters $q$. In the denominator of the probability measure above, the partition function $Z$ ensures that $\phi^{\xi}_{p,q,G} \big( \cdot \big)$ is a probability measure. From the random cluster measure, there exista a bijection to Eulerian loop configurations, in which the probability measure, under such a transformation, instead takes the form,

\begin{align*}
    \phi^{a,b}_{p,\sqrt{2},G} \big( \omega \big)    = \frac{\sqrt{2}^{|\mathrm{loops}|} x \big( p \big)^{|\mathrm{open \text{ } bonds }|}}{\widetilde{Z}^{a,b} \big( p , G \big) }   \text{ } \text{ , } 
\end{align*}

\noindent which can be obtained by Euler's formula, for the normalizing constant $\widetilde{Z}$, and for,

\begin{align*}
  x \big( p \big) \equiv \frac{p}{\sqrt{2} \big( 1 - p \big)}  \text{ } \text{ . } 
\end{align*}

\noindent Aside from introducing the random-cluster model and its probability measure, the Ising model is another celebrated model of statistical mechanics, which can be defined first with the Hamiltonian,

\begin{align*}
      \mathcal{H} \big( \sigma , \Lambda \big) \equiv \mathcal{H} \big( \sigma \big) \equiv \underset{i,j\in \Lambda}{\underset{i \sim j}{\sum}}   J_{ij} \sigma_i \sigma_j + \underset{i \in \partial \Lambda}{\sum}   h_i \sigma_i    \text{ } \text{ , } 
\end{align*}

\noindent with the corresponding probability measure,

\begin{align*}
  \textbf{P}^{\chi}_{G} \big( \sigma \big) \equiv    \frac{\mathrm{exp} \big( - \beta \mathcal{H} \big( \sigma \big) \big) }{Z^{\chi} \big( \sigma , G \big) }     \text{ } \text{ , } 
\end{align*}

\noindent for boundary conditions $\chi$ at inverse temperature $\beta$, and partition function $Z$ so that the expression above is a probability measure. With respect to $\textbf{P}^{\chi}_{G} \big( \cdot \big)$, there exists a relation between the probability measures of the Ising and random-cluster models, in which,

\begin{align*}
  \textbf{E}^{\mathrm{free}}_G \big[ \sigma \big( 0 \big) \sigma \big( a \big)  \big]  = \phi^0_{p,2,G} \big( 0 \longleftrightarrow a \big)   \text{ } \text{ , } 
\end{align*}

\noindent under free boundary conditions, in which the spin-spin correlations at the origin and point $a$ under the expectation of the Ising model are equivalent to the probability of a connectivity event, $\big\{ 0 \longleftrightarrow a \big\}$, between the points $0$ and $a$, occurring under free boundary conditions in the random-cluster model.

\bigskip

\noindent The notion of duality has played a significant role in several models of statistical physics. For the random-cluster model, the self-dual point, $p_{\mathrm{sd}} \big( q \big)$, takes the form,

\begin{align*}
     p_{\mathrm{sd}} \big( q \big) \equiv  \frac{\sqrt{q}}{\sqrt{q}+1}  \text{ } \text{ , } 
\end{align*}

\noindent arises from the conditions,

\begin{align*}
    p^{*} \big( p , q \big) \equiv \frac{\big( 1 - p \big) q}{\big( 1 - p \big) q + p }    \text{ } \text{ , } 
\end{align*}

\noindent and,

\begin{align*}
   \frac{ p^{*} p }{\big(1 - p^{*} \big) \big( 1 - p \big) } = q    \text{ } \text{ , } 
\end{align*}

\noindent being satisfied. Under the assumption that $q \equiv 2$ in the random-cluster model, the fermionic observable takes the form,

\begin{align*}
  F \big( e \big) \equiv \phi^{a,b}_{p,G} \big[ \mathrm{exp} \big( \frac{i}{2} W_{\gamma} \big( e , e_b \big) \big) \textbf{1}_{e \in \gamma} \big]   \text{ } \text{ , } 
\end{align*}

\noindent where the quantity $W_{\gamma} \big( \cdot , \cdot \big)$ in the exponential denotes the winding number of the path $\gamma$, namely the number of left and right turns of the exploration path, ie the total rotation of the path measured in radians between $e$ and $e_b$.

\bigskip

\noindent Observables of the form, and similar forms, above have many applications. In [6], it was shown that the correlation length diverges as $p \longrightarrow p_c \big( q \big)$ from above, which is given by the reciprocal of the infimum,

\begin{align*}
  \frac{1}{\xi \big( p \big)} \equiv   - \underset{n > 0 }{\mathrm{inf}}   \text{ } \frac{1}{n} \text{ } \mathrm{log} \big[        \phi^0_{\textbf{Z}^2,p,q} \big( 0 \longleftrightarrow \big( n , 0 \big) \big)    \big]       \longrightarrow 0\text{ } \text{ . } 
\end{align*}

\noindent On the other hand, in the Ising model, an observable can be defined, which shares connections with massive, and massless s-holomorphicity, [13], which is a function of the form,

\begin{align*}
      f \big( a, z \big) \equiv  \frac{1}{Z^{\chi}} \underset{\gamma : a \longrightarrow z}{\sum} \mathrm{exp} \big[    - 2 \beta \big|  \mathrm{edges} \big( \gamma \big)       \big| - \frac{i}{2} W_{\gamma} \big( a , z \big)      \big]   \text{ } \text{ , } 
\end{align*}

\noindent over a discrete domain $\Omega$, the union of faces of $\textbf{Z}^2$, for the midpoints of edges $a$ and $z$. In addition to the observables $F \big( e \big)$ and $f \big( a , z \big)$ introduced above, it continues to remain of interest to determine properties of observables for other models [6].

\subsection{Observables for other models of statistical physics}

\noindent Under the assumption that such observables for other models exist, in [15] it was demonstrated that observables exist for the staggered six-vertex, eight-vertex, and Ashkin-Teller models, which we recount below. The fact that observables exist for the three other models mentioned follows from the fact that contour integrals can be expressed, discretely, as,

\begin{align*}
  \underset{(i,j)\in C}{\sum} F \big( z_{ij} \big) \big( z_j - z_i \big) \equiv 0  \text{ } \text{ , }
\end{align*}

\noindent for a contour $C$, and $F \big( z_{ij} \big)$ a complex-valued function defined over the midpoints $z_{ij}$ of the edges $\big(ij\big)$. From the condition above, discrete parafermions can be defined for the self-dual Ising, Potts, $\mathrm{O}\big( n \big)$, Ashkin-Teller, and eight-vertex models. For the Ashkin-Teller model, the parafermionic observable, $\psi$, takes the form,

\begin{align*}
  \psi = \mathrm{exp} \big( i W_{\gamma} \big) \sigma \mu_{\tau^{\prime}}  \text{ } \text{ , } 
\end{align*}

\noindent where, as in the definitions of the parafermionic observables for random-cluster and FK percolation, the power of the exponential is proportional to the winding number, in addition to the factors $\sigma$ and $\mu_{\tau^{\prime}}$, which respectively denote the cluster between the points $z_1$ and $z_2$ of the square lattice, and a $\tau^{\prime}$ domain wall, while for the staggered eight-vertex model, the parafermionic observable takes the same form.

\bigskip

\noindent For the $\mathrm{O}\big( n \big)$ model, the parafermionic observable was shown to take the form, for $\Omega \subsetneq \textbf{H}, 
 $, [12],

\begin{align*}
      F^{\mathrm{loop}} \big( a , z , x ,\sigma \big) \equiv F^{\mathrm{loop}} \big( z \big)  \equiv \underset{\gamma \subset \Omega}{\underset{\gamma : a \longrightarrow z}{\sum}} \mathrm{exp} \big( - i \sigma W_{\gamma} \big( a , z \big) \big) x^{l ( \gamma ) }  \text{ } \text{ , } 
\end{align*}

\noindent where $\sigma$ denotes a real parameter appearing in the winding term of the path in the power of the exponential. Lastly, for the $Z_N$ model, the observable takes a similar form as that of the parafermionic observable for the Ashkin-Teller model, in which, [1],

\begin{align*}
  \psi \big( r , \vec{r} \big) \equiv \psi \big( r \big)  \equiv \mathrm{exp} \big( - i \sigma_m \theta \big( r , \vec{r} \big) \big) s_m \big( r \big) \mu_m \big( r \big)   \text{ } \text{ , } 
\end{align*}

\noindent where in the definition of the observable above, $\mu_m \big( r \big)$ denotes a disorder variable, and $m \in \big\{ 1 , \cdots , N-1 \big\}$, and $s_m \big( r \big)$ denotes a function of the spin of the observable. To apply the operator formalism developed for the Ising model at the critical temperature $\beta_c$, it is also of importance to identify the critical parameters for the $\mathrm{O} \big( n \big)$ and Ashkin-Teller models, which are respectively given by, [12], 

\begin{align*}
  x_c \big( n \big) \equiv \frac{1}{\sqrt{2+\sqrt{2-n}}}  \text{ } \text{ , } 
\end{align*}

\noindent for $0 \leq n < 2$, which has a probability measure of the form,

\begin{align*}
        \textbf{P}^{\mathrm{loop} , \xi}_{\Lambda_{\textbf{H}}} \big[ \sigma \big] = \frac{x^{n ( \sigma)} n^{l(\sigma)}}{Z^{\mathrm{loop},\xi}_{\Lambda_{\textbf{H}}} \big( \sigma \big) }   \text{ } \text{ , } 
\end{align*}

\noindent under boundary conditions $\xi$. Later, in \textit{1}, we will apply the operator framework to the high-temperature expansion of the loop probability measure above, which takes the form, 

\begin{align*}
\textbf{P}^{\mathrm{loop}, \xi}_{\Lambda_{\textbf{H}}} \big[ \sigma \big] \overset{\mathrm{HTE}}{\sim }         \frac{n^{k ( \sigma )} x^{e(\sigma)} \mathrm{exp} \big(    h r \big( \sigma \big)  + h^{\prime}    r^{\prime} \big( \sigma \big)     \big) }{Z^{\mathrm{loop},\xi}_{\Lambda_{\textbf{H}}}\big( \sigma \big)  }          
\underset{n \equiv 1 , h^{\prime} \equiv 0}{\overset{\beta \equiv \frac{1}{2} | \mathrm{log}  x |}{\longleftrightarrow}} \textbf{P}^{\mathrm{Ising} , \chi}_{\Lambda_{\textbf{Z}^2}} \big[ \sigma^{\mathrm{Ising}} \big] \equiv \mathrm{O} \big( 1 \big) \text{ } \mathrm{measure}  \text{ } \text{ , } 
\end{align*}

\noindent The second critical point, along a line of critical points, for the Ashkin-Teller model takes the form, [5],

\begin{align*}
  J \equiv U \equiv  \frac{1}{4} \mathrm{log} \big( 3 \big) \text{ } \text{ , } 
\end{align*}

\noindent from the fact that the Ashkin-Teller and Ising models belong to the same universality class, where $J$ denotes one of the coupling constants appearing in the Ashkin-Teller Hamiltonian,

\begin{align*}
 \mathcal{H} \equiv \mathcal{H}^{\mathrm{AT}} \equiv - \underset{i \sim j}{\sum} \bigg[     J \big( \tau \big( i \big) \tau \big( j \big) + \tau^{\prime} \big( i \big) \tau^{\prime} \big( j \big) \big) + U \big( \tau \big( i \big) \tau \big( j \big) \tau^{\prime} \big( i \big) \tau^{\prime} \big( j \big) \big)   \bigg]  \text{ } \text{ , } 
\end{align*}

\noindent for two sites $i,j$ on the square lattice, with $\big( \tau , \tau^{\prime} \big) \in \big\{ \pm 1 \big\}^V \times \big\{ \pm 1 \big\}^V$, with the corresponding probability measure,

\begin{align*}
  \textbf{P}^{\mathrm{AT},\xi}_{\Lambda} \big[ \omega \big] \equiv \frac{\mathrm{exp} \big( \mathcal{H} \big( \omega \big) \big)}{Z^{\xi}_{\Lambda} \big( \omega \big) }  \text{ } \text{ , } 
\end{align*}

\noindent under boundary conditions $\xi$, $\Lambda \subsetneq \textbf{Z}^2$, an Ashkin-Teller configuration $\omega$, and partition function,

\begin{align*}
 Z^{\xi}_{\Lambda} \big( \omega \big)  \equiv \underset{\omega \in \Omega}{\sum}  \mathrm{exp} \big( \mathcal{H} \big( \omega \big) \big) \text{ } \text{ . } 
\end{align*}

\noindent To make the statement of several definitions the easiest in future sections (including, the function $\tau_x$ in \textit{2.1}, and the transfer matrix construction in \textit{2.2.2}, which is dependent upon the construction of several bases spanning an analogue of the Clifford algebra identified in [13]), we restrict ourselves to coupling two Potts models together with two colors in $\textbf{P}^{\mathrm{AT},\xi}_{\Lambda} \big[ \cdot \big]$. 

\subsection{Paper overview}

\noindent From the background of the operator formalism for the Ising model, in the next section we introduce the notions of massive, and massless, s-holomorphicity. In previous work, [13], such notions were applied to study the Ising model, in which the parafermionic observable satisfied the massless s-holmorphic equations at criticality, and the massive s-holomorphic equations away from criticality. Similar notions of holomorphicity are available not only for the Askin-Teller model, but also for the loop model. More specifically, in the next section, before discussing the staggered eight-vertex and odd eight-vertex models in \textit{3}, has the following breakdown:

\begin{itemize}
\item[$\bullet$] \textit{2.1}: Holomorphicity at, and below, the critical temperature $\beta_c$, 

\item[$\bullet$] \textit{2.2}: Transfer matrices,

\item[$\bullet$] \textit{2.3}: Fock representation, 

\item[$\bullet$] \textit{2.4}: Operator correlations and observables,

\item[$\bullet$] \textit{2.5}: Low temperature expansions of parafermionic observables,

\item[$\bullet$] \textit{2.6}: Riemann-Poincare-Steklov operators.

\end{itemize}

\noindent The main result of this effort culminates in two convolution operators, denoted $\big( U^{\textbf{b}}_{\Omega} \big)^{\mathrm{AT}}$, and $\big( U^{\textbf{b}}_{\Omega} \big)^{\mathrm{loop}}$, introduced in \textit{2.6}, whose kernels are the parafermionic observables for the Ashkin-Teller and loop $\mathrm{O} \big( 1 \big)$ models.

\section{Massive, and massless, s-holomorphicity}

\noindent As a beginning case, we restate results of s-holomorphicity for the Ising model.

\subsection{Holomorphicity at, and below, the critical temperature $\beta_c$}

\noindent For the Ising model, it is know that the critical temperature $\beta_c \equiv \frac{1}{2} \mathrm{log} \big( \sqrt{2} + 1 \big)$. The statement of the two results below describes how massless s-holomorphicity occurs at the critical temperature, while below the critical temperature, in the subcritical regime, massive s-holomorphicity occurs.

\bigskip

\noindent \textbf{Definition} \textit{1} (\textit{massive s-holomorphicity below the critical temperature}, [13]). Fix $\beta >0$. For $\beta < \beta_c$, for a complex number $\nu \equiv \bar{\lambda}^3 \frac{\alpha+i}{\alpha-i}$, with $\alpha \equiv \mathrm{exp} \big( - 2 \beta \big)$ and $\lambda \equiv \mathrm{exp} \big( \frac{i \pi}{4} \big)$, the function $F : \Omega \longrightarrow \textbf{C}$ is said to be massive s-holomorphic if,

\begin{align*}
   F \big( N \big) + \nu^{-1} \lambda \bar{F \big( N \big) } = \nu^{-1} F \big( E \big) + \lambda \bar{F \big( E \big) } \text{ } \text{ , } \\ F \big( N \big) + \nu \lambda^{-1} \bar{F \big( N \big)} = \nu F \big( W \big) + \lambda^{-1} \bar{F \big( W \big)} \text{ } \text{ , } \\ F \big( S \big) + \nu \lambda^3 \bar{F \big( S \big) } = \nu F \big( E \big) + \lambda^3 \bar{F \big( E \big)} \text{ } \text{ , } \\ F \big( S \big) + \nu^{-1} \lambda^{-3} \bar{F \big( S \big)} = \nu^{-1} F \big( W \big) + \lambda^{-3} \bar{F \big( W \big) } \text{ } \text{ , } 
\end{align*}

\noindent for any face of $\Omega$ with edges E,N,W,S.

\bigskip

\noindent A similar definition holds for massless s-holomorphicity at the critical temperature, which is stated through a second definition below.

\bigskip

\noindent \textbf{Definition} \textit{2} (\textit{massless s-holomorphicity at the critical temperature}, [13]). Fix $\beta >0$. For $\beta \equiv \beta_c$, given the same quantities $\alpha$ and $\lambda$, the function $F : \Omega \longrightarrow \textbf{C}$ is said to be massless s-holomorphic if,

\begin{align*}
   F \big( N \big) +  \lambda \bar{F \big( N \big) } =  F \big( E \big) + \lambda \bar{F \big( E \big) } \text{ } \text{ , } \\ F \big( N \big) + \lambda^{-1} \bar{F \big( N \big)} =  F \big( W \big) + \lambda^{-1} \bar{F \big( W \big)} \text{ } \text{ , } \\ F \big( S \big) + \lambda^3 \bar{F \big( S \big) } =  F \big( E \big) + \lambda^3 \bar{F \big( E \big)} \text{ } \text{ , } \\ F \big( S \big) + \lambda^{-3} \bar{F \big( S \big)} =  F \big( W \big) + \lambda^{-3} \bar{F \big( W \big) } \text{ } \text{ , } 
\end{align*}

\noindent for any face of $\Omega$ with edges E,N,W,S.

\bigskip

\noindent Besides the relations above which hold about a single face of $\Omega$, one must also define the propagation operators for $\beta < \beta_c$, and also for $\beta \equiv \beta_c$.

\bigskip

\noindent \textbf{Lemma} \textit{1} (\textit{propagation mechanism for the Ising model away from the critical temperature}, \textbf{Lemma} \textit{6}, [13]). For $\beta < \beta_c$, and the interval,

\begin{align*}
  \textbf{I}^{*} \equiv \big\{  a + \frac{1}{2} , a + \frac{3}{2} , \cdots , b - \frac{1}{2}        \big\}   \text{ } \text{ , } 
\end{align*}

\noindent with left, and right, endpoints respectively given by $k_L \equiv a + \frac{1}{2}$ and $k_R \equiv b - \frac{1}{2}$, dual to an interval of $\textbf{Z}$,

\begin{align*}
  \textbf{I} \equiv \big\{ a , a+1 , \cdots , b \big\}  \text{ } \text{ , } 
\end{align*}

\noindent the s-holomorphic propagator $P : \big( \textbf{R}^2 \big)^{\textbf{I}^{*}} \longrightarrow \big( \textbf{R}^2 \big)^{\textbf{I}^{*}}$ satisfies, for $\lambda \equiv \mathrm{exp} \big( \frac{i \pi}{4}\big)$, and $k$ such that $k \in \textbf{I}^{*} \backslash \big\{ k_L , k_R \big\}$,

\begin{align*}
     \big( P f \big) \big( k \big)   \equiv        \frac{1}{\sqrt{2} \lambda^3}      f \big( k - 1 ) + 2 f \big( k \big) + \frac{\lambda^3}{\sqrt{2}} f \big( k +1 \big)  + \frac{1}{\sqrt{2}} \bar{f \big( k - 1 \big)} - \sqrt{2} \bar{f \big( k \big)} + \frac{1}{\sqrt{2}} \bar{f \big( k + 1\big)}        \text{ } \text{ , } \\     \big( P f \big) \big( k_L \big) \equiv     \big( 1 + \frac{1}{\sqrt{2}} \big)  f \big( k_L \big) + \frac{\lambda^3}{\sqrt{2}} f \big( k_L + 1 \big) +      \big( \lambda^3 + \frac{1}{\sqrt{2} \lambda^3} \big)             \bar{f \big( k_L \big)} + \frac{1}{\sqrt{2} } \bar{f \big( k_L + 1 \big) }                    \text{ } \text{ , } \\         \big( P f \big) \big( k_R \big) \equiv          \frac{1}{\sqrt{2} \lambda^3} f \big( k_R - 1 \big) + \big( 1 + \frac{1}{\sqrt{2}} \big) f \big( k_R\big) + \frac{1}{\sqrt{2}} \bar{f \big( k_R - 1 \big)} + \big( \frac{1}{\lambda^3} + \frac{\lambda^3}{\sqrt{2}} \big) \bar{f \big( k_R \big) }  \text{ } \text{ . } 
\end{align*}

\bigskip

\noindent \textbf{Lemma} \textit{2} (\textit{propagation mechanism for the Ising model at the critical temperature}, \textbf{Lemma} \textit{6}, [13]). For $\beta \equiv \beta_c$, and the interval $\textbf{I}^{*}$ defined in \textbf{Lemma} \textit{1}, the s-holomorphic propagator $P_{\beta} : \big( \textbf{R}^2 \big)^{\textbf{I}^{*}} \longrightarrow \big( \textbf{R}^2 \big)^{\textbf{I}^{*}}$ satisfies, for $k$ such that $k \in \textbf{I}^{*} \backslash \big\{ k_L , k_R \big\}$,

\begin{align*}
     \big( P_{\beta} f \big) \big( k \big)   \equiv     
        \big( \frac{-S-i}{2S} \big)      f \big( k -1 \big) + \frac{C^2}{S}  f \big( k \big) + 
\frac{i - S}{2S} f \big( k +1 \big)   +  \frac{C}{2S} \bar{f \big( k -1 \big)} - C \bar{f \big( k \big)} + \frac{C}{2S} \bar{f \big( k +1 \big) }                  \text{ } \text{ , } \\     \big( P_{\beta} f \big) \big( k_L \big) \equiv  \frac{\big( S+C \big) C}{2S}   f \big( k_L \big) + \frac{i-S}{2S} f \big( k_L + 1 \big) +                     \bar{f \big( k_L \big) } +    \frac{C}{2S}       \bar{f \big( k_L + 1 \big) }                     \text{ } \text{ , } \\         \big( P_{\beta} f \big) \big( k_R \big) \equiv - \frac{S_i}{2S}  f \big( k_R - 1 \big) + 
\frac{\big( S+ C \big) C}{2S} f \big( k_R \big) + \frac{C}{2S} \bar{f \big( k_R - 1 \big)} + \frac{- \big( S+ C\big)S + i \big( C - S \big)}{2S} \bar{f \big( k_R \big)}       \text{ } \text{ , } 
\end{align*}

\noindent with $S \equiv \mathrm{sinh} \big( 2 \beta \big)$ and $C \equiv \mathrm{cosh} \big( 2 \beta \big)$.

\bigskip

\noindent Over the boundary $\partial \Omega$ of the finite volume, boundary conditions satisfy the following property.

\bigskip

\noindent \textbf{Definition} \textit{3} (\textit{Riemann boundary conditions for the Ising model}, \textbf{Definition} \textit{3}, [13]). A function $f : \Omega \longrightarrow \textbf{C}$ is said to satisfy Riemann boundary conditions, if, for an edge $z \in \partial \Omega$,

\begin{align*}
 f \big( z \big) \big|\big|  \frac{1}{\sqrt{\tau_{\mathrm{cw}} \big( z \big)}}  \text{ } \text{ , } 
\end{align*}

\noindent namely that $f \big( z \big)$ is a multiple of $\tau_{\mathrm{cw}} \big( z \big)$, in which,

\begin{align*}
e \text{ } \mathrm{vertical} \text{ } \Longleftrightarrow   \tau_{\mathrm{cw}} \big( e \big) \in \big\{ \pm 1  \big\} \text{ } \text{ , } \\   e \text{ } \mathrm{horizontal} \text{ }  \Longleftrightarrow   \tau_{\mathrm{cw}} \big( e \big) \in \big\{ \pm i  \big\}  
 \text{ } \text{ . } \tag{**}
\end{align*}

\bigskip

\noindent With $P$ and $P_{\beta}$, another component of the operator formalism for the Ising model is the transfer matrix. As a mapping from the state space of the Ising model into itself, the transfer matrix takes the following form.

\bigskip

\noindent \textbf{Definition} \textit{4} (\textit{transfer matrix of the Ising model}, [13]). For the Ising model, the transfer matrix takes the form,

\begin{align*}
               V \equiv \big( V^h \big)^{\frac{1}{2}} V^V \big( V^h \big)^{\frac{1}{2}}  \text{ } \text{ , } 
\end{align*}

\noindent for, 

\begin{align*}
\big( V^h \big)^{\frac{1}{2}} \equiv   \mathrm{exp} \big( \beta \underset{a \leq i \leq b}{\sum} \sigma_i \rho_i \big)   \text{ } \text{ , } \\  V^V \equiv    \mathrm{exp} \big( \frac{\beta}{2} \underset{a \leq i \leq b-1}{\sum}   \sigma_i \sigma_{i+1} \big)      \text{ } \text{ , } 
\end{align*}

\noindent under the basis $\big\{ e_{\sigma}\big\}$, in which $\big( V^h \big)^{\frac{1}{2}} \equiv 0$ if $\sigma \neq \rho$, and similarly, $V^V \equiv 0$ if $\sigma_a \neq \rho_a$ and $\sigma_b \neq \rho_b$.

\bigskip

\noindent The fact that the Askin-Teller model, and Ising model, belong to the same universality class allows for us to make use of the notions of massive, and massless, s-holomorphicity above, with the following.

\bigskip

\noindent \textbf{Definition} \textit{5} (\textit{massive s-holomorphicity below the critical temperature of the Ashkin-Teller model}). Fix $J,U \in \textbf{R}$. Over the square grid $\Omega$, the Ashkin-Teller observable $\psi \big( r \big)$, for couplings $J,U < \frac{1}{4} \mathrm{log} \big( 3 \big)$, satisfies,

\begin{align*}
   \psi \big( N \big) + \nu^{-1} \lambda \bar{\psi \big( N \big) } = \nu^{-1} \psi \big( E \big) + \lambda \bar{\psi \big( E \big) } \text{ } \text{ , } \\ \psi \big( N \big) + \nu \lambda^{-1} \bar{\psi \big( N \big)} = \nu \psi \big( W \big) + \lambda^{-1} \bar{\psi \big( W \big)} \text{ } \text{ , } \\ \psi \big( S \big) + \nu \lambda^3 \bar{\psi \big( S \big) } = \nu \psi \big( E \big) + \lambda^3 \bar{\psi \big( E \big)} \text{ } \text{ , } \\ \psi \big( S \big) + \nu^{-1} \lambda^{-3} \bar{\psi \big( S \big)} = \nu^{-1} \psi \big( W \big) + \lambda^{-3} \bar{\psi \big( W \big) } \text{ } \text{ , } 
\end{align*}

\noindent for parameters $\lambda$ and $\nu$ provided in \textbf{Definition} \textit{1}, and any face of $\Omega$ with edges E,N,W,S.

\bigskip

\noindent \textbf{Definition} \textit{6} (\textit{massless s-holomorphicity at the critical temperature of the Ashkin-Teller model}). Fix $J, U \in \textbf{R}$. Over the square grid $\Omega$, the Ashkin-Teller observable $\psi \big( r \big)$, for couplings $J \equiv U \equiv \frac{1}{4} \mathrm{log} \big( 3 \big)$, satisfies, 

\begin{align*}
   \psi \big( N \big) +  \lambda \bar{\psi \big( N \big) } = \psi \big( E \big) + \lambda \bar{\psi \big( E \big) } \text{ } \text{ , } \\ \psi \big( N \big) + \lambda^{-1} \bar{\psi \big( N \big)} = \psi \big( W \big) + \lambda^{-1} \bar{\psi \big( W \big)} \text{ } \text{ , } \\ \psi \big( S \big) + \lambda^3 \bar{\psi \big( S \big) } =  \psi \big( E \big) + \lambda^3 \bar{\psi \big( E \big)} \text{ } \text{ , } \\ \psi \big( S \big) + \lambda^{-3} \bar{\psi \big( S \big)} = \psi \big( W \big) + \lambda^{-3} \bar{\psi \big( W \big) } \text{ } \text{ , } 
\end{align*}

\noindent for the parameter $\lambda$ provided in \textbf{Definition} \textit{1}, and any face of $\Omega$ with edges E,N,W,S.

\bigskip

\noindent For the loop $\mathrm{O} \big( n \big)$ model, notions of discrete holomorphicity have already been employed by Duminil-Copin and Smirnov to demonstrate that the connective constant of the honeycomb lattice is $\sqrt{2+\sqrt{2}}$, [12], at the critical point, which is a weaker form of s-holomorphicity.

\bigskip

\noindent \textbf{Lemma} (\textit{discrete holomorphicity of the loop parafermionic observable for the connective constant of the honeycomb lattice}, \textbf{Lemma} \textit{1}, [12]). Fix $\sigma \equiv \frac{5}{8}$. For the loop $\mathrm{O}\big( n \big)$ observable, at the critical point $x_c \big( n \big)$, given any vertex within a finite volume of the hexagonal lattice,

\begin{align*}
    \big( p - v \big) F^{\mathrm{loop}} \big( p \big) + \big( q - v \big) F^{\mathrm{loop}} \big( q \big) + \big( r - v \big) F^{\mathrm{loop}} \big( r \big) = 0    \text{ } \text{ , } 
\end{align*}

\noindent for the midpoints of edges $p,q,r$ adjacent to $v$.

\bigskip

\noindent From the notion of discrete holomorphicity above at the critical point of the loop $\mathrm{O} \big( n \big)$ model, to introduce notions of massive, and massless, s-holomorphicity for the loop $\mathrm{O} \big( n \big)$ model, one needs to not only make use of different parameters $\lambda$ and $\alpha$ than those which are introduced in \textbf{Definition} \textit{1} for the Ising model, but also to account for the total number of terms appearing in the Cauchy-Riemann equations. Immediately, one can observe that the Cauchy-Riemann equations for the Ising and Ashkin-Teller models differ from that of the loop model from the fact that in the former there are four terms appearing for edges oriented along the E,N,W,S portions of each face of the square lattice, in comparison to there only being three edges oriented along each face of the triangular lattice in the latter.

\bigskip

\noindent For the following definition of s-holomorphicity, along the lines of previous comments above, we make use of the relation originally provided for midpoints of the square lattice, [13],

\begin{align*}
        f \big( e_v \big) + \frac{i}{\theta} \bar{f \big( e_v \big) } = f \big( e_w \big) + \frac{i}{\theta} \bar{f \big( e_w \big)  }  \text{ } \text{ , } 
\end{align*}

\bigskip

\noindent for some function $f$, angle,

\begin{align*}
   \theta \equiv \frac{2u-v-w}{\big|2u - v - w  \big|} \text{ } \text{ , } 
\end{align*}

\noindent and edges,

\begin{align*}
 e_v \equiv  vu   \text{ } \text{ , } \\ e_w \equiv  wu \text{ } \text{ . } 
\end{align*}

\bigskip

\noindent For the following definitions of s-holomorphicity below for the loop $\mathrm{O} \big( n \big)$ model, introduce,

\begin{align*}
  e_1 \equiv \bar{e_1} \equiv -1   \text{ } \text{ , } \\ e_2 \equiv \mathrm{exp} \big( i \frac{4 \pi}{3} \big)   \text{ } \text{ , } \\  \bar{e_2} \equiv \mathrm{exp} \big(  - i \frac{4 \pi}{3} \big)   \text{ } \text{ , } \\ e_3 \equiv \bar{e_3} \equiv -1 \text{ } \text{ , } \\   e_4 \equiv \mathrm{exp} \big( i \frac{\pi}{3} \big)  \text{ } \text{ , }  \\ \bar{e_4} \equiv \mathrm{exp} \big( -  i \frac{\pi}{3} \big)  \text{ } \text{ . } 
\end{align*}

\noindent We postpone introducing the construction of the transfer matrix for the loop model to \textit{2.2.2}. For the definition of the observable below that is defined in terms of the loop transfer matrix, this quantity is equivalent to the expected value of the observable,

\begin{align*}
  \mathscr{O} \equiv K^{|\gamma|} \lambda^{t_r} \bar{\lambda}^{t_l}  \text{ } \text{ , } 
\end{align*}

\noindent for the number of paths in a loop configuration, $\big| \gamma \big|$, and the number of left, and right, turns of the path, $t_r$ and $t_l$, respectively.

\bigskip

\noindent \textbf{Definition} \textit{7} (\textit{massive s-holmorphicity of the loop parafermionic observable below Nienhuis' critical point}, [14]). For $x < x_c \big( n \big)$, the observable,

\begin{align*}
   F \big( k , m \big) \equiv \frac{1}{Z} \bra{\alpha} T_3^{M-k} T_2 \big( m \big) T^k_1 \ket{\beta} \text{ } \text{ , } 
\end{align*}

\noindent defined in terms of the transfer matrix $T$, $m \in \big\{ 1 , 2 , \cdots , M \big\}$, and $\alpha, \beta>0$, satisfies,

\begin{align*}
       F \big( z_1 \big) + \bar{e_1}^{2s} \bar{F \big( z_1 \big) } = F \big( z_2 \big) + \bar{e_1}^{2s} \bar{F \big( z_2 \big) }      \text{ } \text{ , }   \end{align*}

       \begin{align*}
     F \big( z_2 \big) + \bar{e_2}^{2s} \bar{F \big( z_2 \big)} = F \big( z_3 \big) + \bar{e_2}^{2s} \bar{F \big( z_3 \big)} \text{ } \text{ , }  \\    F \big( z_3 \big) + \bar{e_3}^{2s} \bar{F \big( z_3 \big)} = F \big( z_4 \big) + \bar{e_3}^{2s} \bar{F \big( z_4 \big) }     \text{ } \text{ , }  \\ F \big( z_4 \big) + \bar{e_4}^{2s} \bar{F \big( z_4 \big) } = F \big( z_1 \big) + \bar{e_4}^{2s} \bar{F \big( z_1 \big) } \text{ } \text{ , }
\end{align*}

\noindent for any face of $\textbf{H}$ with edges $e_1 , e_2, e_3, e_4$.

\bigskip

\noindent \textbf{Definition} \textit{8} (\textit{massless s-holmorphicity of the loop parafermionic observable at Nienhuis' critical point}, [14]). For $x \equiv x_c \big( n \big)$, the observable,

\begin{align*}
   F \big( k , m \big) \equiv \frac{1}{Z} \bra{\alpha} T_3^{M-k} T_2 \big( m \big) T^k_1 \ket{\beta} \text{ } \text{ , } 
\end{align*}

\noindent defined in terms of the transfer matrix $T$, $m \in \big\{ 1 , 2 , \cdots , M \big\}$, and $\alpha, \beta>0$, given a parameter,

\begin{align*}
 \nu^{\mathrm{loop}} \equiv  \bar{\lambda}^2 \frac{\alpha^{\prime}+i}{\alpha^{\prime}-i } \equiv \bigg[ \mathrm{exp} \big( - i \frac{\pi}{4} \big) \bigg]^2 \frac{n+i}{n-i}  \text{ } \text{ , } 
\end{align*}

\noindent for $0 \leq n <2$, satisfies,

\begin{align*}
       F \big( z_1 \big) + \big( \nu^{\mathrm{loop}}\big)^{-1}\bar{e_1}^{2s} \bar{F \big( z_1 \big) } = \big( \nu^{\mathrm{loop}}\big)^{-1}  F \big( z_2 \big) + \bar{e_1}^{2s} \bar{F \big( z_2 \big) }      \text{ } \text{ , }   \\  
     F \big( z_2 \big) + \big( \nu^{\mathrm{loop}}\big)^{-1}\bar{e_2}^{2s} \bar{F \big( z_2 \big)} =  \big( \nu^{\mathrm{loop}}\big)^{-1}  F \big( z_3 \big) + \bar{e_2}^{2s} \bar{F \big( z_3 \big)} \text{ } \text{ , }  \\    F \big( z_3 \big) + \big( \nu^{\mathrm{loop}}\big)^{-1}\bar{e_3}^{2s} \bar{F \big( z_3 \big)} = \big( \nu^{\mathrm{loop}}\big)^{-1}  F \big( z_4 \big) + \bar{e_3}^{2s} \bar{F \big( z_4 \big) }     \text{ } \text{ , }  \\ F \big( z_4 \big) + \big( \nu^{\mathrm{loop}}\big)^{-1}\bar{e_4}^{2s} \bar{F \big( z_4 \big) } = \big( \nu^{\mathrm{loop}}\big)^{-1} F \big( z_1 \big) + \bar{e_4}^{2s} \bar{F \big( z_1 \big) } \text{ } \text{ , }
\end{align*}

\noindent for any face of $\textbf{H}$ with edges $e_1 , e_2, e_3, e_4$.

\bigskip

\noindent In the setting of the Ising model, the two generators of the Clifford algebra,

\begin{align*}
   p_k  \text{, } \\ q_k \text{, } 
\end{align*}

\noindent act on the basis elements $e_{\sigma}$, through the action,

\begin{align*}
       p_k \big(   e_{\sigma}     \big)   =   \sigma_{k + \frac{1}{2}} e_{\tau}  \text{ } \text{ , } \\     q_k \big(  e_{\sigma}      \big)  = i \sigma_{k - \frac{1}{2}} e_{\tau} \text{ } \text{ , } \tag{*}
\end{align*}

\noindent for a configuration indexed with $k$ over the dual interval $\textbf{I}^{*}$, namely over the sample space over $\textbf{I}$, $\big\{ \pm 1 \big\}^{\textbf{I}}$, where $\tau$ is given by,

\[
\tau_x  \equiv  \text{ } 
\left\{\!\begin{array}{ll@{}>{{}}l}     \sigma_x,   &  \text{if} \text{ }  x > k  \text{, } \\
- \sigma_x,   & \text{if} \text{ }  x < k  \text{, } 
\end{array}\right. 
\]

\noindent in which, from the definition of $\tau_x$ above, for $x<k$ the spin is flipped to any of the $q-1$ remaining colors of the Potts model. From definitions of s-holomorphicity for the Askin-Teller and $\mathrm{O} \big( n \big)$ models, one must also encode boundary conditions for the discrete boundary value problems, which is achieved with the following.

\bigskip

\noindent \textbf{Definition} \textit{9} (\textit{Riemann boundary conditions for the Ashkin-Teller model}). For the Ashkin-Teller model, a function $f : \Omega \longrightarrow \textbf{C}$ is said to satisfy Riemann boundary conditions, if for an edge $z$ incident to the boundary,

\begin{align*}
 f \big( z \big) \big|\big|  \frac{1}{\sqrt{\tau_{\mathrm{cw}} \big( z \big)}}  \text{ } \text{ , } 
\end{align*}

\noindent under the convention provided in (**).

\bigskip

\noindent \textbf{Definition} \textit{10} (\textit{Riemann boundary conditions for the loop model}). For the loop model, a function $f : \Omega \longrightarrow \textbf{C}$ is said to satisfy Riemann boundary conditions, if for an edge $z$ incident to the boundary of some finite volume of the hexagonal lattice,

\begin{align*}
 f \big( z \big) \big|\big|  \frac{1}{\sqrt{\tau_{\mathrm{cw}} \big( z \big)}}  \text{ } \text{ , } 
\end{align*}

\noindent under the convention provided in (**).

\bigskip

\noindent Equipped with massive, and massless, s-holomorphicity, as well as boundary conditions for discrete Riemann boundary value problems, below we also provide properties of the propagation matrices for the Ashkin-Teller and loop models.

\bigskip

\noindent \textbf{Lemma} \textit{3} (\textit{propagation mechanism for the Ashkin-Teller model at the critical} $\frac{1}{4} \mathrm{log} \big( 3 \big)$ \text{threshold}). For $J , U < \frac{1}{4} \mathrm{log} \big( 3 \big)$ and some function $f \equiv f^{\mathrm{AT}}$, the propagation mechanism $P^{\mathrm{AT}}: \big( \textbf{R}^2 \big)^{\textbf{I}^{*}} \longrightarrow \big( \textbf{R}^2 \big)^{\textbf{I}^{*}}$ satisfies, for $k \in \textbf{I}^{*} \backslash \big\{ k_L, k_R \big\}$,

\begin{align*}
  \big( P^{\mathrm{AT}} f \big) \big( k \big) \equiv         \frac{\lambda^{-3}}{\sqrt{2}} f \big( k -1 \big) + 2 f \big( k \big) + \frac{\lambda^3}{\sqrt{2}} f \big( k + 1 \big) + \frac{1}{\sqrt{2}} \bar{f \big( k -1 \big) } - \sqrt{2} \bar{f \big( k \big) } + \frac{1}{\sqrt{2}} \bar{f \big( k +1 \big) }  \text{ } \text{ , } \\   \big( P^{\mathrm{AT}} f \big) \big( k_L  \big) \equiv  \big( 1 + \frac{1}{\sqrt{2}} \big) f \big( k_L \big) + \frac{\lambda^3}{\sqrt{2}}  f \big( k_L + 1 \big) + \big( \lambda^3 + \frac{\lambda^{-3}}{\sqrt{2}} \big)  \bar{f \big( k_L  \big) } + \frac{1}{\sqrt{2}} \bar{f \big( k_L + 1 \big) }   \text{ } \text{ , } \\ \big( P^{\mathrm{AT}} f \big) \big( k_R  \big) \equiv         \frac{\lambda^{-3}}{\sqrt{2}} f \big( k_R - 1 \big) + \big( 1 + \frac{1}{\sqrt{2}} \big) f \big( k_R \big) + \frac{1}{\sqrt{2}} \bar{f \big( k_R - 1 \big) } + \big( \lambda^{-3} + \frac{\lambda^3}{\sqrt{2}} \big) \bar{f \big( k_R \big) }          \text{ } \text{ . } 
\end{align*}

\bigskip

\noindent \textbf{Lemma} \textit{4} (\textit{propagation mechanism for the Ashkin-Teller model below the critical} $\frac{1}{4}\mathrm{log} \big( 3 \big)$ \textit{threshold}). For $J \equiv U \equiv \frac{1}{4} \mathrm{log} \big( 3 \big)$ and some function $f \equiv f^{\mathrm{AT}}$, the propagation mechanism $P^{\mathrm{AT}}_{\frac{1}{4} \mathrm{log} (3 )}: \big( \textbf{R}^2 \big)^{\textbf{I}^{*}} \longrightarrow \big( \textbf{R}^2 \big)^{\textbf{I}^{*}}$ satisfies, for $k \in \textbf{I}^{*} \backslash \big\{ k_L, k_R \big\}$,

\begin{figure}
\begin{align*}
\\ \includegraphics[width=0.23\columnwidth]{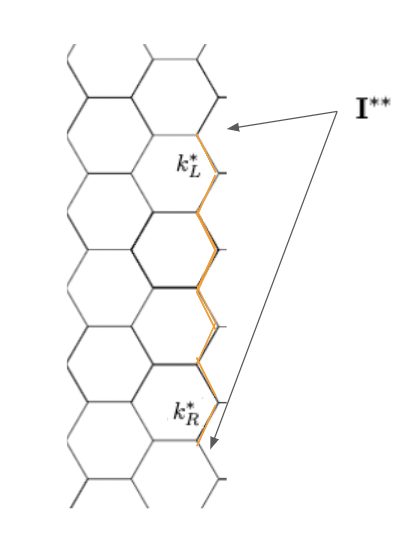}
\end{align*}
\caption{A depiction of the interval $\textbf{I}^{**}$, with endpoints at the middle of edges, $k_L^{*}$ and $k_R^{*}$.}
\end{figure}

\begin{align*}
  \big( P^{\mathrm{AT}}_{\frac{1}{4} \mathrm{log} (3 )} f  \big) \big( k \big) \equiv  \big( \frac{-S-i}{2S} \big)      f \big( k -1 \big) + \frac{C^2}{S}  f \big( k \big) + 
\frac{i - S}{2S} f \big( k +1 \big)   +  \frac{C}{2S} \bar{f \big( k -1 \big)} - C \bar{f \big( k \big)} + \frac{C}{2S} \bar{f \big( k +1 \big) }                  \text{ } \text{ , } \\   \big( P^{\mathrm{AT}}_{\frac{1}{4} \mathrm{log} (3 )} f  \big) \big( k_L \big) \equiv                \frac{\big( S+C \big) C}{2S}   f \big( k_L \big) + \frac{i-S}{2S} f \big( k_L + 1 \big) +                     \bar{f \big( k_L \big) } +    \frac{C}{2S}       \bar{f \big( k_L + 1 \big) }                           \text{ } \text{ , } \\ \big( P^{\mathrm{AT}}_{\frac{1}{4} \mathrm{log} (3 )} f  \big) \big( k_R \big) \equiv  - \frac{S_i}{2S}  f \big( k_R - 1 \big) + 
\frac{\big( S+ C \big) C}{2S} f \big( k_R \big) + \frac{C}{2S} \bar{f \big( k_R - 1 \big)} + \frac{- \big( S+ C\big)S + i \big( C - S \big)}{2S} \bar{f \big( k_R \big)}     \text{ } \text{ , } 
\end{align*}

\noindent with $S \equiv \mathrm{sinh} \big( 2 \beta \big)$ and $C \equiv \mathrm{cosh} \big( 2 \beta \big)$.

\bigskip

\noindent For the propagation mechanisms for the loop model, we make use of a nonempty interval of the hexagonal lattice, which we denote with $\textbf{I}^{**}$ (see the figure above for a visual representation of $\textbf{I}^{**}$). 

\bigskip

\noindent \textbf{Lemma} \textit{5} (\textit{propagation mechanism for the loop model at Nienhuis' critical point}). For $x< x_c \big( n \big)$ and some function $f \equiv f^{\mathrm{loop}}$, the propagation mechanism $P^{\mathrm{loop}}_{\frac{1}{\sqrt{2+\sqrt{2-n}}}} \equiv P^{\mathrm{loop}}_{x_c ( n )} \equiv P^{\mathrm{loop}} : \big( \textbf{H} \big)^{\textbf{I}^{**}} \longrightarrow \big( \textbf{H} \big)^{\textbf{I}^{**}}$ satisfies, for $k \in \textbf{I}^{**} \backslash \big\{ k^{*}_L , k^{*}_R \big\} \equiv \textbf{I}^{**} \backslash \big\{ k_L , k_R \big\}$,

\begin{align*}
  \big( P^{\mathrm{loop}} f \big) \big( k \big) \equiv  \frac{\lambda^{-2}}{\sqrt{3}} f \big( k -1 \big) + 2 f \big( k \big) + \frac{\lambda^2}{\sqrt{3}} f \big( k+1 \big)+ \frac{1}{\sqrt{3}} \bar{f \big( k - 1 \big) } - \sqrt{3} \bar{f \big( k \big)} + \frac{1}{\sqrt{3}} \bar{f \big( k +1 \big) } \text{ } \text{ , } \\ \big( P^{\mathrm{loop}} f \big) \big( k_L \big) \equiv \big( 1 + \frac{1}{\sqrt{3}} \big)     f \big( k_L \big) + \frac{\lambda^2}{\sqrt{3}} f \big( k_L + 1 \big) + \big( \lambda^2 + \frac{\lambda^{-2}}{\sqrt{3}} \big) \bar{f \big( k_L  \big) } + \frac{1}{\sqrt{3}} \bar{f \big( k_L  + 1 \big)}       \text{ } \text{ , } \\   \big( P^{\mathrm{loop}} f \big) \big( k_R \big) \equiv     \frac{\lambda^{-2}}{\sqrt{3}} f \big( k_R - 1 \big) + \big( 1 + \frac{1}{\sqrt{3}} \big) f \big( k_R \big) + \frac{1}{\sqrt{3}} \bar{f \big( k_R - 1 \big)} + \big( \lambda^{-2} + \frac{\lambda^3}{\sqrt{3}} \big) \bar{f \big( k_R \big) }    \text{ } \text{ . } 
\end{align*}

\bigskip

\noindent For the remaining definition of the other loop propagator below Nienhuis' critical point, we make use of the following fact. Independently of the number of loops, recall that the loop probability measure can be expressed possesses a high-temperature expansion, [6],

\begin{align*}
\textbf{P}^{\mathrm{loop}, \xi}_{\Lambda_{\textbf{H}}} \big[ \sigma \big] \overset{\mathrm{HTE}}{\sim }       \frac{n^{k ( \sigma )} x^{e(\sigma)} \mathrm{exp} \big(    h r \big( \sigma \big)  + h^{\prime}    r^{\prime} \big( \sigma \big)     \big) }{Z^{\mathrm{loop},\xi}_{\Lambda_{\textbf{H}}} \big( \sigma \big) }          
\underset{n \equiv 1 , h^{\prime} \equiv 0}{\overset{\beta \equiv \frac{1}{2} | \mathrm{log}  x |}{\longleftrightarrow}} \textbf{P}^{\mathrm{Ising} , \chi}_{\Lambda_{\textbf{Z}^2}} \big[ \sigma^{\mathrm{Ising}} \big]  \equiv \mathrm{O} \big( 1 \big) \text{ } \mathrm{measure} \text{ } \text{ , } 
\end{align*}

\noindent for an Ising model configuration $\sigma^{\mathrm{Ising}}$, under the corresponding probability measure supported over $\textbf{Z}^2$ with boundary conditions $\chi$, for the number of connected components $k \big( \sigma \big)$ of the loop configuration, and also for the factors,

\begin{align*}
  r \big( \sigma \big) \equiv \sum_{u \in G} \sigma_u   \text{ } \text{ , } \\ r^{\prime} \big( \sigma \big) \equiv \sum_{t \equiv \{ u , v , w \}} \textbf{1}_{\{ \sigma_u \equiv \sigma_v \equiv \sigma_w \}} \text{ } \text{ , } 
\end{align*}                

\noindent corresponding to the two external fields, which respectively represent the summation over all spins in finite volume, and the number of monochromatically colored triangles. Hence we can further apply results for the propagator matrix for the loop $\mathrm{O} \big( 1 \big)$ model in the absence of one external field. To distinguish the inverse temperature of the Ising model that is in correspondence with the $\mathrm{O} \big( 1 \big)$ model, denote the following inverse temperature,

\begin{align*}
 \beta^{\mathrm{loop}} \equiv \frac{1}{2} \big| \mathrm{log} \big(  x \big)  \big|    \text{ } \text{ . } 
\end{align*}

\bigskip

\noindent \textbf{Lemma} \textit{6} (\textit{propagation mechanism for the loop model below Nienhuis' critical point}). For $x \equiv x_c \big( n \big)$ and some function $f \equiv f^{\mathrm{loop}}$, the propagation mechanism $P^{\mathrm{loop}}_x : \big( \textbf{H} \big)^{\textbf{I}^{**}} \longrightarrow \big( \textbf{H} \big)^{\textbf{I}^{**}}$ satisfies, for $k \in \textbf{I}^{**} \backslash \big\{ k^{*}_L , k^{*}_R \big\} \equiv \textbf{I}^{**} \backslash \big\{ k_L , k_R \big\}$,

\begin{align*}
       \big( P^{\mathrm{loop}}_{x} f  \big) \big( k \big) \equiv        \frac{-S^{\mathrm{loop}}-i}{2S^{\mathrm{loop}}}              f \big( k -1 \big) + \frac{\big( C^{\mathrm{loop}} \big)^2}{S^{\mathrm{loop}}} f \big( k \big) +  \frac{-S^{\mathrm{loop}} + i }{2S^{\mathrm{loop}}} f \big( k +1 \big) +  \frac{C^{\mathrm{loop}}}{2S^{\mathrm{loop}}} \bar{f \big( k -1 \big)} - \cdots \\ C^{\mathrm{loop}} \bar{f \big( k \big)} + \frac{C^{\mathrm{loop}}}{2S^{\mathrm{loop}}} \bar{f \big( k +1 \big) }     \text{ } \text{ , } \\  \big( P^{\mathrm{loop}}_{x} f \big) \big( k_L \big) \equiv \frac{\big( S^{\mathrm{loop}} + C^{\mathrm{loop}}\big) C^{\mathrm{loop}} }{2S^{\mathrm{loop}}} f \big( k_L \big) + \frac{i-S^{\mathrm{loop}}}{2S^{\mathrm{loop}}} f \big( k_L + 1 \big) + \cdots \\ \big(  \frac{- \big( S^{\mathrm{loop}} +C^{\mathrm{loop}} \big) S^{\mathrm{loop}} + i \big( C^{\mathrm{loop}} - S^{\mathrm{loop}}\big) }{2 S^{\mathrm{loop}}}\big) \bar{f \big( k_L \big)} + \cdots \end{align*}

       \begin{align*}
       \frac{C^{\mathrm{loop}}}{2S^{\mathrm{loop}}} \bar{f \big( k_L + 1 \big) }  \text{ } \text{ , } \\   \big( P^{\mathrm{loop}}_{x} f \big) \big( k_R \big) \equiv \frac{-S^{\mathrm{loop}} - i }{2S^{\mathrm{loop}}} f \big( k_R - 1 \big) + \big( \frac{\big( S^{\mathrm{loop}}+C^{\mathrm{loop}} \big) C^{\mathrm{loop}}}{2S^{\mathrm{loop}}} \big)  f \big( k_R \big) + \frac{C^{\mathrm{loop}}}{2S^{\mathrm{loop}}} \bar{f \big( k_R - 1 \big)} + \cdots \\  \big( \frac{- \big(S^{\mathrm{loop}} + C^{\mathrm{loop}} \big) S^{\mathrm{loop}} + i \big( C^{\mathrm{loop}} - S^{\mathrm{loop}}\big) }{2S^{\mathrm{loop}}} \big) \bar{f \big( k_R\big) } 
 \text{ } \text{ , } 
\end{align*}

\noindent given constants satisfying,

\begin{align*}
            S^{\mathrm{loop}} \equiv \mathrm{sinh} \big( 2 \beta^{\mathrm{loop}} \big)   \text{ } \text{ , } \\        C^{\mathrm{loop}} \equiv \mathrm{cosh}  \big( 2 \beta^{\mathrm{loop}}  \big)      \text{ } \text{ . } 
\end{align*}

\noindent Below, we state results regarding the propagation mechanisms $P^{\mathrm{AT}}$, $P^{\mathrm{AT}}_{\frac{1}{4}\mathrm{log} (3)}$, $P^{\mathrm{loop}}$ and $P^{\mathrm{loop}}_{x}$.

\bigskip

\noindent \textbf{Proposition} \textit{1} (\textit{propagation mechanism of the Ashkin-Teller model at criticality}, \textbf{Proposition } \textit{7}, [13]). The matrix $P^{\mathrm{AT}}_{\frac{1}{4}\mathrm{log}(3)}$ is symmetric, with eigenvalues $\lambda^{\mathrm{AT},\pm}_{\alpha}$, where each $\big\{\lambda_{\alpha}\big\}_{1 \leq \alpha \leq | \textbf{I}^{*}|}$ is distinct and with magnitude strictly larger than $1$.

\bigskip

\noindent \textit{Proof of Proposition 1}. This is an application of the argument for \textbf{Proposition} \textit{7} contained in [13]. \boxed{}

\bigskip

\noindent \textbf{Proposition} \textit{2} (\textit{propagation mechanism of the loop model at Nienhuis' critical point}, \textbf{Proposition} \textit{7}, [13]). The matrix $P^{\mathrm{loop}}_{x_c(n)}$ is symmetric, with eigenvalues $\lambda^{\mathrm{loop},\pm}_{\alpha}$, where each $\big\{ \lambda^{\mathrm{loop}}\big\}_{1 \leq \alpha < |\textbf{I}^{**}|}$ is distinct and with magnitude strictly larger than $1$.

\bigskip

\noindent \textit{Proof of Proposition 2}. We present more dramatic changes to the arguments rather than those provided in \textbf{Proposition} \textit{7} contained in [13]. To begin, observe,

\begin{align*}
      \bigg[  \big( P^{\mathrm{loop}}_{x} f  \big) \big( k \big)\bigg]^{-1} \equiv  \bigg[  \frac{-S^{\mathrm{loop}}-i}{2S^{\mathrm{loop}}}  f \big( k -1 \big) + \frac{\big( C^{\mathrm{loop}} \big)^2}{S^{\mathrm{loop}}} f \big( k \big) + \frac{-S^{\mathrm{loop}} + i }{2S^{\mathrm{loop}}} f \big( k +1 \big) + \frac{C^{\mathrm{loop}}}{2S^{\mathrm{loop}}} \bar{f \big( k -1 \big)} + \cdots \\ C^{\mathrm{loop}} \bar{f \big( k \big)} + \frac{C^{\mathrm{loop}}}{2S^{\mathrm{loop}}} \bar{f \big( k +1 \big) } \bigg]^{-1}  \\ \equiv  \frac{-S^{\mathrm{loop}}-i}{2S^{\mathrm{loop}}} f \big( k -1 \big)^{-1} + \frac{\big( C^{\mathrm{loop}} \big)^2}{S^{\mathrm{loop}}} f\big( k \big)^{-1} + \frac{-S^{\mathrm{loop}}+i }{2 S^{\mathrm{loop}}} f \big(k+1 \big)^{-1} + \frac{{C}^{\mathrm{loop}}  }{2 S^{\mathrm{loop}}} \bar{f \big( k-1\big)}^{-1} + C^{\mathrm{loop}}  \bar{f \big( k \big)}^{-1} + \cdots \\ \frac{C^{\mathrm{loop}}}{2 \mathrm{S}^{\mathrm{loop}}} \bar{f \big( k +1 \big)}^{-1} \\  \equiv  - \frac{i + S^{\mathrm{loop}}}{2 S^{\mathrm{loop}}} f \big( k-2 \big) + \frac{\big( {C}^{\mathrm{loop}}\big)^2 }{S^{\mathrm{loop}}} f \big( k -1 \big) + \frac{-S^{\mathrm{loop}} + i }{2 S^{\mathrm{loop}}} f \big( k \big) +  \frac{C^{\mathrm{loop}}}{2 S^{\mathrm{loop}}} \bar{f \big( k -2 \big)} + C^{\mathrm{loop}} \bar{f \big( k -1 \big)} + \cdots \\ \frac{C^{\mathrm{loop}}}{2S^{\mathrm{loop}} } \bar{f \big( k \big) } \text{ } \text{ , } \end{align*}

      \noindent corresponding to backpropagation from $k$,

      \begin{align*}
      \bigg[ \big( P^{\mathrm{loop}}_{x} f \big) \big( k_L \big) \bigg]^{-1} \equiv \bigg[  \frac{ \big( S^{\mathrm{loop}}+ C^{\mathrm{loop}}\big) C^{\mathrm{loop}} }{2S^{\mathrm{loop}}}         f \big( k_L \big) + \frac{i - S^{\mathrm{loop}}}{2S^{\mathrm{loop}}} f \big( k_L + 1 \big) +     \cdots \\    \big(  \frac{- \big( S^{\mathrm{loop}} +C^{\mathrm{loop}} \big) S^{\mathrm{loop}} + i \big( C^{\mathrm{loop}} - S^{\mathrm{loop}}\big) }{2 S^{\mathrm{loop}}}\big)     \bar{f \big( k_L \big)} +   \frac{\mathrm{C}^{\mathrm{loop}}}{2 S^{\mathrm{loop}}} \bar{f \big( k_L + 1 \big) } \bigg]^{-1}  \\  \equiv        \frac{ \big( S^{\mathrm{loop}} +  C^{\mathrm{loop}} \big) C^{\mathrm{loop}}}{2S^{\mathrm{loop}}}     f \big(k_L\big)^{-1} +         \frac{i - S^{\mathrm{loop}}}{C^{\mathrm{loop}}}           f \big(k_L + 1 \big)^{-1} + 
\cdots \\  \big(  \frac{- \big( S^{\mathrm{loop}} +C^{\mathrm{loop}} \big) S^{\mathrm{loop}} + i \big( C^{\mathrm{loop}} - S^{\mathrm{loop}}\big) }{2 S^{\mathrm{loop}}}\big)  \bar{f \big( k_L \big)}^{-1} + 
 \frac{C^{\mathrm{loop}}}{2 S^{\mathrm{loop}}} \bar{f \big( k_L + 1 \big)}^{-1} \\  \equiv     \frac{ \big( S^{\mathrm{loop}} +  C^{\mathrm{loop}} \big) C^{\mathrm{loop}}}{2S^{\mathrm{loop}}}     f \big( k_L - 1 \big) +        \frac{i - S^{\mathrm{loop}}}{C^{\mathrm{loop}}}     f \big( k_L \big) + \cdots \\   \big(  \frac{- \big( S^{\mathrm{loop}} +C^{\mathrm{loop}} \big) S^{\mathrm{loop}} + i \big( C^{\mathrm{loop}} - S^{\mathrm{loop}}\big) }{2 S^{\mathrm{loop}}}\big)  \bar{f \big( k_L -1 \big)} + \frac{C^{\mathrm{loop}}}{2 S^{\mathrm{loop}}} \bar{f \big( k_L \big) }          \text{ }  \text{ , } \end{align*}
      
     \noindent corresponding to backpropagation from $k_L$, 
      \begin{align*}
      \bigg[  \big( P^{\mathrm{loop}}_{x} f \big) \big( k_R \big) \bigg]^{-1} \equiv \bigg[ - \frac{i + S^{\mathrm{loop}}}{2 S^{\mathrm{loop}}} f \big( k_R - 1 \big) + \big( \frac{\big( S^{\mathrm{loop}} + C^{\mathrm{loop}} \big) C^{\mathrm{loop}}}{2 S^{\mathrm{loop}}}\big)  f \big( k_R \big) + \cdots \end{align*}

      \begin{align*}
      \frac{C^{\mathrm{loop}}}{2 S^{\mathrm{loop}}}  \bar{f \big( k_R - 1 \big)} + \big( \frac{- \big(S^{\mathrm{loop}} + C^{\mathrm{loop}} \big) S^{\mathrm{loop}} + i \big( C^{\mathrm{loop}} - S^{\mathrm{loop}}\big) }{2S^{\mathrm{loop}}} \big)  \bar{f \big( k_R\big) } \bigg]^{-1} \\ \equiv  - \frac{i + S^{\mathrm{loop}}}{2 S^{\mathrm{loop}}}   f \big( k_R -1 \big)^{-1} + \big( \frac{\big( S^{\mathrm{loop}} + C^{\mathrm{loop}} \big) C^{\mathrm{loop}}}{2 S^{\mathrm{loop}}}\big)  
 f \big( k_R \big)^{-1} + \frac{C^{\mathrm{loop}}}{2 S^{\mathrm{loop}}}  \bar{f \big( k_R - 1 \big)}^{-1} + \cdots \\ \big( \frac{- \big(S^{\mathrm{loop}} + C^{\mathrm{loop}} \big) S^{\mathrm{loop}} + i \big( C^{\mathrm{loop}} - S^{\mathrm{loop}}\big) }{2S^{\mathrm{loop}}} \big) \bar{f \big( k_R \big)}^{-1} \\  \equiv  - \frac{i + S^{\mathrm{loop}}}{2 S^{\mathrm{loop}}}  f \big( k_R -2 \big) +  \big( \frac{\big( S^{\mathrm{loop}} + C^{\mathrm{loop}} \big) C^{\mathrm{loop}}}{2 S^{\mathrm{loop}}}\big)    f \big( k_R - 1 \big) + \frac{C^{\mathrm{loop}}}{2 S^{\mathrm{loop}}}  \bar{f \big( k_R \big) } + \cdots \\  \big( \frac{- \big(S^{\mathrm{loop}} + C^{\mathrm{loop}} \big) S^{\mathrm{loop}} + i \big( C^{\mathrm{loop}} - S^{\mathrm{loop}}\big) }{2S^{\mathrm{loop}}} \big)\bar{f \big( k_R - 1 \big) } \text{ } \text{ , } 
\end{align*}

\noindent corresponding to backpropagation from $k_R$, in which the entries of the inverse matrix, $\big( P^{\mathrm{loop}}_{x_c(n)} \big)^{-1}$ are defined from the expressions provided above. Moreover, below Nienhuis' critical point, from \textbf{Definition} \textit{7}, the set of equations about $z_1 , z_2 , z_3 , z_4$ being discretely holomorphic implies the set of relations transform to,

\begin{align*}
   i    \bar{ \big( P^{\mathrm{loop}}_{x} f  \big) \big( k \big)  } \equiv \frac{-S^{\mathrm{loop}}-i}{2S^{\mathrm{loop}}}            \bar{  f \big( k -1 \big) } + \frac{\big( C^{\mathrm{loop}} \big)^2}{S^{\mathrm{loop}}} \bar{ f \big( k \big) } +  \frac{-S^{\mathrm{loop}} + i }{2S^{\mathrm{loop}}}  \bar{f \big( k +1 \big)} +  \frac{C^{\mathrm{loop}}}{2S^{\mathrm{loop}}} \bar{\bar{f \big( k -1 \big)}} - \cdots \\   C^{\mathrm{loop}} \bar{\bar{f \big( k \big)}} + \frac{C^{\mathrm{loop}}}{2S^{\mathrm{loop}}} \bar{\bar{f \big( k +1 \big) }}     \text{ } \text{ , }   \\  i \bar{ \big( P^{\mathrm{loop}}_{x} f  \big) \big( k_L  \big)  }  \equiv  \frac{\big( S^{\mathrm{loop}} + C^{\mathrm{loop}}\big) C^{\mathrm{loop}} }{2S^{\mathrm{loop}}} \bar{f \big( k_L \big)} + \frac{i-S^{\mathrm{loop}}}{2S^{\mathrm{loop}}} \bar{f \big( k_L + 1 \big)} + \cdots \\ \big(  \frac{- \big( S^{\mathrm{loop}} +C^{\mathrm{loop}} \big) S^{\mathrm{loop}} + i \big( C^{\mathrm{loop}} - S^{\mathrm{loop}}\big) }{2 S^{\mathrm{loop}}}\big) \bar{\bar{f \big( k_L \big)}} + \cdots \\   \frac{C^{\mathrm{loop}}}{2S^{\mathrm{loop}}} \bar{\bar{f \big( k_L + 1 \big) }}  \text{ } \text{ , }  \\   \text{ } \text{ , } \\ i \bar{ \big( P^{\mathrm{loop}}_x f  \big) \big( k_R \big)}  \equiv  \frac{-S^{\mathrm{loop}} - i }{2S^{\mathrm{loop}}} \bar{f \big( k_R - 1 \big)} + \big( \frac{\big( S^{\mathrm{loop}}+C^{\mathrm{loop}} \big) C^{\mathrm{loop}}}{2S^{\mathrm{loop}}} \big)  \bar{f \big( k_R \big)} + \frac{C^{\mathrm{loop}}}{2S^{\mathrm{loop}}} \bar{\bar{f \big( k_R - 1 \big)}} + \cdots \\  \big( \frac{- \big(S^{\mathrm{loop}} + C^{\mathrm{loop}} \big) S^{\mathrm{loop}} + i \big( C^{\mathrm{loop}} - S^{\mathrm{loop}}\big) }{2S^{\mathrm{loop}}} \big) \bar{\bar{f \big( k_R\big) }}       \text{ } \text{ , } 
\end{align*}

\noindent under the correspondence,

\begin{align*}
        F^{\mathrm{loop}} \big(      k_L    \big) \underset{\phi_i} {\overset{i \bar{ F^{\mathrm{loop}} ( \cdot)}} {\longleftrightarrow}}    F^{\mathrm{loop}} \big( k_R \big)      \text{ } \text{ . } 
\end{align*}

\noindent Hence,

\begin{align*}
   \underset{\lambda}{\bigcup}  \big\{    \lambda  \big| \lambda \in \mathrm{spec} \big( P^{\mathrm{loop}}_x  \big)   \big\} \equiv               \underset{\lambda}{\bigcup} \big\{  \lambda \big| \lambda \in \mathrm{spec} \big( \big( P^{\mathrm{loop}}_x  \big)^{-1} \big)    \big\} \text{ } \text{ , } 
\end{align*}

\noindent where,

\begin{align*}
   \big( P^{\mathrm{loop}}_x \big)^{-1}  \equiv \textbf{i} \cdot P^{\mathrm{loop}}_x \cdot \textbf{i}^{-1} \text{ } \text{ , } 
\end{align*}

\noindent for the involution,

\begin{align*}
 \textbf{i} : \big( \textbf{H} \big)^{|\textbf{I}^{**} | } \longrightarrow   \big( \textbf{H} \big)^{|\textbf{I}^{**} | }  \text{ } \text{ , } \\  f \mapsto i \bar{f}  \text{ } \text{ . } 
\end{align*}

\noindent Also, the fact that $P^{\mathrm{loop}}_x \equiv \big( P^{\mathrm{loop}}_x  \big)^{\mathrm{T}}$ comes from the observation, from the loop propagator at the critical loop parameter $x_c \big( n \big)$, that,

\begin{align*}
  \bigg[     \big( P^{\mathrm{loop}} f \big) \big( \eta \bar{k} \big)    \bigg]^{\mathrm{T}} \equiv  \bigg[ \frac{\lambda^{-2}}{\sqrt{3}} f \big( \eta_{\bar{k}-1} \big) + 2 f \big( \eta_{\bar{k}} \big) + \frac{\lambda^2}{\sqrt{3}} f \big( \eta_{\bar{k}+1} \big)+ \frac{1}{\sqrt{3}} \bar{f \big( \eta_{k^{\prime}} - 1 \big) } - \sqrt{3} \bar{f \big( \eta_{k^{\prime}} \big)} + \frac{1}{\sqrt{3}} \bar{f \big( \eta_{{k^{\prime}}+1} \big) } \bigg]^{\mathrm{T}} \\ \equiv           \bigg[     \big( P^{\mathrm{loop}} f \big) \big( {k} \big)    \bigg]^{\mathrm{T}}     \text{, }  \\   \bigg[     \big( P^{\mathrm{loop}} f \big) \big( \eta \bar{k_L} \big)    \bigg]^{\mathrm{T}} \equiv   \bigg[   \frac{\lambda^{-2}}{\sqrt{3}} f \big( \eta_{\bar{k_R} - 1} \big) + \big( 1 + \frac{1}{\sqrt{3}} \big) f \big( \eta_{\bar{k_R}} \big) + \frac{1}{\sqrt{3}} \bar{f \big( \eta_{(k_R)^{\prime} - 1)} \big)} + \big( \lambda^{-2} + \frac{\lambda^3}{\sqrt{3}} \big) \bar{f \big( \eta_{(k_R)^{\prime}} \big) }  \bigg]^{\mathrm{T}} \\ \equiv \bigg[     \big( P^{\mathrm{loop}} f \big) \big( {k_L} \big)    \bigg]^{\mathrm{T}}  \text{, }   \\   \bigg[     \big( P^{\mathrm{loop}} f \big) \big( \eta \bar{k_R}  \big)    \bigg]^{\mathrm{T}} \equiv   \bigg[     \frac{\lambda^{-2}}{\sqrt{3}} f \big( \eta_{\bar{k_R} - 1} \big) + \big( 1 + \frac{1}{\sqrt{3}} \big) f \big( \eta_{ \bar{k_R}} \big) + \frac{1}{\sqrt{3}} \bar{f \big( \eta_{(k_R)^{\prime}-1} \big)} + \big( \lambda^{-2} + \frac{\lambda^3}{\sqrt{3}} \big) \bar{f \big( \eta_{(k_R)^{\prime}}  \big) }      \bigg]^{\mathrm{T}} \\   \equiv \bigg[     \big( P^{\mathrm{loop}} f \big) \big( {k_R} \big)    \bigg]^{\mathrm{T}}  \text{, } 
\end{align*}

\noindent for parameters,

\begin{align*}
  \eta_{\bar{k}-1} \equiv \eta \big( \bar{k} - 1 \big)   \text{ } \text{ , } \\ \eta_{\bar{k}}  \equiv  \eta \bar{k} \text{ } \text{ , } \\  \eta_{k^{\prime}} -1  \equiv \eta \big( \bar{k^{\prime}} - 1 \big) \text{ } \text{ , } \\      \eta_{k^{\prime}}    \equiv    \eta \bar{k}  \text{ } \text{ , } \\ \eta_{k^{\prime}+1} \equiv  \eta \big( \bar{k^{\prime}} + 1 \big) \text{ } \text{ , } \\        \eta_{\bar{k_R}-1}         \equiv  \eta \big( \bar{k_R} - 1 \big) \text{ } \text{ , } \\       \eta_{\bar{k_R}}         \equiv  \eta \bar{k_R} \text{ } \text{ , } \\           \eta_{(k_R)^{\prime}-1 } \equiv \eta      \big( \bar{k_R} - 1 \big) \text{ } \text{ , } \\                     \eta_{\bar{k_R}}       \equiv       \eta_{(k_R)^{\prime}}  \text{ } \text{ , }  \\   \eta_{\bar{k_R} - 1}             \equiv    \eta \big( \bar{k_R} - 1 \big)   \text{ } \text{ , }  \\   \eta_{\bar{k_R}}         \equiv  \eta \bar{k_R} \text{ } \text{ , } \\ \eta_{(k_R)^{\prime}-1 }  \equiv \eta \big(  \bar{k_R} - 1 \big) \text{ } \text{ , }  \\  \eta_{(k_R)^{\prime}}\equiv \eta \bar{k_R} \text{ } \text{ . }   
\end{align*}

\noindent Also, observe that two maps, one from $\textbf{I}^{**}_{\frac{1}{2}} \overset{\varphi_1}{\longrightarrow} \textbf{C}$, and the other from $\textbf{C} \overset{\varphi_2}{\longrightarrow} \textbf{I}^{**}_0 \equiv \textbf{I}^{**}$, implies, that the propagator applied to either $\varphi_1$, or to $\varphi_2$, takes the form, 

\begin{align*}
  A \varphi_1  \text{, } \\  A \varphi_2   \text{, }
\end{align*}

\noindent for,

\begin{align*}
  A : \big( \textbf{R}^2 \big)^{\textbf{I}^{*}} \longrightarrow  \big( \textbf{R}^2 \big)^{\textbf{I}^{*}}   \text{ } \text{ . } 
\end{align*}

\noindent From these objects, one has,

\[
\begin{bmatrix} \big( P^{\mathrm{loop}} \big)_{11} & \cdots & \big( P^{\mathrm{loop}} \big)_{n1}\\\vdots & \vdots & \vdots \\  \big( P^{\mathrm{loop}} \big){1n} & \cdots & \big( P^{\mathrm{loop}} \big)_{nn}\end{bmatrix} \begin{bmatrix} \big( P^{\mathrm{loop}} \big)_{11} & \cdots & \big( P^{\mathrm{loop}} \big)_{n1}\\\vdots & \vdots & \vdots \\  \big( P^{\mathrm{loop}} \big){1n} & \cdots & \big( P^{\mathrm{loop}} \big)_{nn}\end{bmatrix} ^{\mathrm{T}} \geq 0 \text{ } \text{ . } 
\]

\noindent Away from the critical point, similar arguments hold which demonstrate that the propagation matrix is continuous in $x_c \big( n \big)$, and hence does not have any eigenvalues in its spectrum which are zero.

\bigskip

\noindent To demonstrate that $1$ cannot be an eigenvalue of $P^{\mathrm{loop}}_x$, or of $P^{\mathrm{loop}}$, argue by contradiction. That is, if $1$ were can eigenvalue of either propagation matrix at, or below the critical parameter $x_c \big( n \big)$, then one would either have that $P^{\mathrm{loop}}_xf= f$, or that $P^{\mathrm{loop}}=f$. For the first case at criticality, to show that one arrives to a contradiction and that $f\equiv 0$, define the extension of the observable with the action, 

\begin{align*}
  i \big[       h \big(   x + \frac{1}{2}   \big) + h \big(  x - \frac{1}{2} \big)   \big] = h \big(  x + \frac{i}{2} \big) + h \big( x - \frac{i}{2}  \big) \equiv 0 \text{ } \text{ , } 
\end{align*}

\noindent from the mapping,

\begin{align*}
  h^{\mathrm{loop}} \equiv h \equiv \varphi_2 \cdot \varphi_1 :  \textbf{I}^{**}_{0,\frac{1}{2},1} \longrightarrow \textbf{C} \text{ } \text{ , } 
\end{align*}

\noindent for the codomain,

\begin{align*}
 \textbf{I}^{**}_{0,\frac{1}{2},1}  \equiv  \textbf{I}^{**}_0 \cup  \textbf{I}^{**}_{\frac{1}{2}} \cup  \textbf{I}^{**}_1   \text{ } \text{ . } 
\end{align*}

\noindent Therefore $h$ is identically zero, and vanishes, implying that $1$ is not an eigenvalue of $P^{\mathrm{loop}}_x$, from the fact that the loop observable $F^{\mathrm{loop}} \equiv 0$ . To show that the same observation holds for $P^{\mathrm{loop}}$ for another contradiction, for $x \neq x_c \big( n \big)$, observe that, for the interval $\textbf{I}^{**}_{\frac{1}{2}}$,

\begin{align*}
     h \big( x + 1 \big)            \equiv  h \big( x \big) + \frac{1}{i} \big[ h \big( x + \frac{i}{2} + \frac{1}{2} \big) - h \big( x - \frac{i}{2} + \frac{1}{2}    \big)  \big]      \equiv h \big( x \big)           \text{ , } 
\end{align*}

\noindent or, similarly, that,

\begin{align*}
        h \big( x \big) \equiv h \big( x - \frac{1}{4} \big)  + \frac{1}{i} \big[  h \big(  x + \frac{i}{2} -  \frac{1}{2}   \big) - h \big(    x - \frac{i}{2} - \frac{1}{2}   \big)              \big] \equiv h \big( x + 1 \big)       \text{ , } 
\end{align*}

\noindent each of which satisfy the relations,

\begin{align*}
       F \big( z_1 \big) + \big( \nu^{\mathrm{loop}}\big)^{-1}\bar{e_1}^{2s} \bar{F \big( z_1 \big) } = \big( \nu^{\mathrm{loop}}\big)^{-1}  F \big( z_2 \big) + \bar{e_1}^{2s} \bar{F \big( z_2 \big) }      \text{ } \text{ , }   \\  
     F \big( z_2 \big) + \big( \nu^{\mathrm{loop}}\big)^{-1}\bar{e_2}^{2s} \bar{F \big( z_2 \big)} =  \big( \nu^{\mathrm{loop}}\big)^{-1}  F \big( z_3 \big) + \bar{e_2}^{2s} \bar{F \big( z_3 \big)} \text{ } \text{ , }  \\    F \big( z_3 \big) + \big( \nu^{\mathrm{loop}}\big)^{-1}\bar{e_3}^{2s} \bar{F \big( z_3 \big)} = \big( \nu^{\mathrm{loop}}\big)^{-1}  F \big( z_4 \big) + \bar{e_3}^{2s} \bar{F \big( z_4 \big) }     \text{ } \text{ , }  \\ F \big( z_4 \big) + \big( \nu^{\mathrm{loop}}\big)^{-1}\bar{e_4}^{2s} \bar{F \big( z_4 \big) } = \big( \nu^{\mathrm{loop}}\big)^{-1} F \big( z_1 \big) + \bar{e_4}^{2s} \bar{F \big( z_1 \big) } \text{ } \text{ , }
\end{align*}

\noindent which implies, from the function $h$, that,

\begin{align*}
     h \big( x + 1 \big) + i \sqrt{\mathrm{Im} \big( \nu \big) } \bar{ h \big( x + 1 \big) } = h \big( x \big) - i  \sqrt{\mathrm{Im} \big( \nu \big) } \bar{h \big( x \big) }       \text{ } \text{ , } 
\end{align*}

\noindent from which we conclude that for the remaining case $h$ is constant, so $1$ is not an eigenvalue of $P^{\mathrm{loop}}$ either, for parameters,

\begin{align*}
  \lambda \text{, } \\ \nu  \text{, } 
\end{align*}

\noindent appearing in the definitions of the propagator matrices at, and away, from $x_c \big( n \big)$.

\bigskip

\noindent With regards to Riemann boundary conditions, we impose boundary conditions at $k_L$, and at $k_R$, with the following. For the leftmost boundary $k_L$, the fact that the function $h$ attains a value of $\mathcal{C}^{\mathrm{loop}} \mathrm{exp} \big( - i \frac{\pi}{3} \big)$, for real $\mathcal{C}^{\mathrm{loop}}$, implies that it also satisfies the relation,

\begin{align*}
          h \big( x_L \big) +    i \sqrt{\mathrm{Im} \big( \nu \big) } \bar{ h \big( x_L \big) } = h \big( x_L -1 \big) - i  \sqrt{\mathrm{Im} \big( \nu \big) } \bar{h \big( x_L -1 \big) } 
 \text{ , } 
\end{align*}

\noindent For the rightmost boundary $k_R$, the fact that the function $h$ attains a value of $\big( \mathcal{C}^{\mathrm{loop}} \big)^{\prime} \mathrm{exp} \big( - i \frac{\pi}{3} \big)$, for another real $\big( \mathcal{C}^{\mathrm{loop}} \big)^{\prime}$, implies that it also satisfies the relation, 

\begin{align*}
         h \big( x_R \big) +    i \sqrt{\mathrm{Im} \big( \nu \big) } \bar{ h \big( x_R \big) } = h \big( x_R -1 \big) - i  \sqrt{\mathrm{Im} \big( \nu \big) } \bar{h \big( x_R -1 \big) }        \text{ . } 
\end{align*}

\noindent To conlude the argument, we show that the eigenvalues of the spectrum are distinct. To show that the eigenspace is one-dimensional, if $f$ denotes an eigenvector over $\big( \textbf{H} \big)^{\textbf{I}^{**}}$, then with a massive s-holomorphic extension $h$ of $f$ (as discussed above), over $\textbf{I}_{\frac{1}{2}} \cup \textbf{I}^{**}_{0,1}$,

\begin{align*}
  \big( P^{\mathrm{loop}} h \big) \big( k \big) \equiv  \frac{\lambda^{-2}}{\sqrt{3}} h \big( k -1 \big) + 2 h \big( k \big) + \frac{\lambda^2}{\sqrt{3}} h \big( k+1 \big)+ \frac{1}{\sqrt{3}} \bar{h \big( k - 1 \big) } - \sqrt{3} \bar{h \big( k \big)} + \frac{1}{\sqrt{3}} \bar{h \big( k +1 \big) } \text{ } \text{ , } \\ \big( P^{\mathrm{loop}} h \big) \big( k_L \big) \equiv \big( 1 + \frac{1}{\sqrt{3}} \big)     h \big( k_L \big) + \frac{\lambda^2}{\sqrt{3}} h \big( k_L + 1 \big) + \big( \lambda^2 + \frac{\lambda^{-2}}{\sqrt{3}} \big) \bar{h \big( k_L  \big) } + \frac{1}{\sqrt{3}} \bar{h \big( k_L  + 1 \big)}       \text{ } \text{ , } \\   \big( P^{\mathrm{loop}} h \big) \big( k_R \big) \equiv     \frac{\lambda^{-2}}{\sqrt{3}} h \big( k_R - 1 \big) + \big( 1 + \frac{1}{\sqrt{3}} \big) h \big( k_R \big) + \frac{1}{\sqrt{3}} \bar{h \big( k_R - 1 \big)} + \big( \lambda^{-2} + \frac{\lambda^3}{\sqrt{3}} \big) \bar{h \big( k_R \big) }    \text{ } \text{ , } 
\end{align*}

\noindent one can solve for $h \big( x+1 \big)$. From the expression for $h$ obtained from the massive s-holomorphic equations that are propagated with $P^{\mathrm{loop}}$ above, we conclude that the eigenspace is one-dimensional, from which we conclude the argument. \boxed{}

\subsection{Transfer matrices}

\subsubsection{Ashkin-Teller model}

\noindent \underline{Generator relations}. For the loop $\mathrm{O} \big( n \big)$ and Ashkin-Teller models, to make use of the operator formalism one needs to introduce transfer matrices over each state space of the model. In order to formulate the transfer matrices for these other models, with respective state spaces $\Omega^{\mathrm{AT}}$ and $\Omega^{\mathrm{O} ( n)}$ for the Ashkin-Teller and loop $\mathrm{O} \big( n \big)$ models, it suffices to demonstrate the existence of several bases for the Ashkin-Teller model, the first of which takes the form,

\begin{align*}
\mathcal{S}_1 \equiv  \mathrm{span} \big\{       e_{\tau ( i )} \text{ }        \big| \text{ }           \tau \big( i \big)  \equiv + 1        \big\}    \text{ } \text{ , } 
\end{align*}

\noindent corresponding to the spins at side $i$ for which $\tau \big( i \big) \equiv +1$, the second of which takes the form,

\begin{align*}
\mathcal{S}_2 \equiv \mathrm{span} \big\{       e_{\tau ( j )} \text{ }        \big| \text{ }           \tau \big( j \big)  \equiv + 1        \big\}    \text{ } \text{ , } 
\end{align*}

\noindent corresponding to the spins at site $j$ for which $\tau \big( j \big) \equiv +1$, the third of which takes the form,

\begin{align*}
\mathcal{S}_3 \equiv \mathrm{span} \big\{       e_{\tau^{\prime} ( i )} \text{ }        \big| \text{ }           \tau^{\prime} \big( i \big)  \equiv + 1        \big\}    \text{ } \text{ , } 
\end{align*}

\noindent corresponding to the spins at site $i$ for which $\tau^{\prime} \big( i \big) \equiv +1$, and finally, the fourth of which takes the form,

\begin{align*}
  \mathcal{S}_4 \equiv \mathrm{span} \big\{       e_{\tau^{\prime} ( j )} \text{ }        \big| \text{ }           \tau^{\prime} \big( j \big)  \equiv + 1        \big\}   \text{ } \text{ . } 
\end{align*}

\noindent corresponding to the spins at site $j$ for which $\tau^{\prime}\big( j \big) \equiv +1$. From $\mathcal{S}_1 , \mathcal{S}_2 , \mathcal{S}_3 , \mathcal{S}_4$ above, observe,

\begin{align*}
 \mathcal{S}_1 \cap \mathcal{S}_2 \equiv \mathrm{span} \big\{ e_{\tau(i)}  , e_{\tau(j)} \text{ } | \text{ } \tau \big( i \big) \tau \big( j \big) \equiv + 1 \big\} \text{ } \text{ , }  \\    \mathcal{S}_3 \cap \mathcal{S}_4 \equiv \mathrm{span} \big\{ e_{\tau(i)}  , e_{\tau(j)} \text{ } | \text{ } \tau^{\prime} \big( i \big) \tau^{\prime} \big( j \big) \equiv + 1 \big\}    \text{ } \text{ , }   \\   \mathcal{S}_1 \cap \mathcal{S}_2 \cap \mathcal{S}_3 \cap \mathcal{S}_4 \equiv \mathrm{span} \big\{ e_{\tau ( i ) } , e_{\tau(j)} , e_{\tau^{\prime}(i) } , e_{\tau^{\prime} ( j )}     \text{ } | \text{ }       \tau \big( i \big) \tau \big( j \big) \tau^{\prime}  \big( i \big) \tau^{\prime} \big( j \big) \equiv + 1        \big\}       \text{ } \text{ . } 
\end{align*}

\noindent Moreover,

\begin{align*}
    \mathrm{dim}_{\Omega^{\mathrm{AT}}} \big( \mathcal{S}_1 \big) =  \mathrm{dim}_{\Omega^{\mathrm{AT}}} \big( \mathcal{S}_2 \big) = \mathrm{dim}_{\Omega^{\mathrm{AT}}} \big( \mathcal{S}_3 \big) = \mathrm{dim}_{\Omega^{\mathrm{AT}}} \big( \mathcal{S}_4 \big) = 2^{|\textbf{I}|} = 2^{b-a} \text{ } \text{ , } \end{align*}

    \begin{align*}
    \mathrm{dim}_{\Omega^{\mathrm{AT}}} \bigg[ \underset{1 \leq i \leq 4}{\bigcap} \mathcal{S}_i  \bigg] = 2^{b-a+1} \text{ } \text{ . } 
\end{align*}

\bigskip

\noindent We define similar basies $\mathcal{S}^{\prime}_1 , \cdots , \mathcal{S}^{\prime}_4$ spanning the remaining color of the Potts models in the next section.

\bigskip

\noindent As a result, in a similar way to how the transfer matrix $V$ is formulated for the Ising model, for two coupled Potts models in the Ashkin-Teller model, the transfer matrix takes the form,

\[
    V^{\mathrm{AT}} \equiv \big( V^{\mathrm{AT},h} \big)^{\frac{1}{2}}  V^{\mathrm{AT},V}  \big( V^{\mathrm{AT},h} \big)^{\frac{1}{2}}   \text{ } \text{ , } 
\]

\bigskip

\noindent under the assignment, for $\tau \big( i \big) \equiv \tau_i$, $\tau \big( i+1\big) \equiv \tau_{i+1}$, $\tau^{\prime} \big( i \big) \equiv \tau^{\prime}_i$ and $\tau^{\prime} \big( j \big) \equiv \tau^{\prime}_j$,

\[
V^{\mathrm{AT},V}_{\tau , \tau^{\prime}} \equiv V^{\mathrm{AT},V}  \equiv  \text{ } 
\left\{\!\begin{array}{ll@{}>{{}}l}   \underset{a \equiv i_0 \sim \cdots \sim i_{n-1} \sim i_n \equiv b}{\sum} \big[     J \big( \tau_i \tau_{i+1} + \tau^{\prime}_i \tau^{\prime}_{i+1} \big) + U \big( \tau_i \tau_{i+1} \tau^{\prime}_{i} \tau^{\prime}_{i+1} \big)   \big]      & \text{ } \text{if} \text{ } \tau_{i_0} \equiv  a \text{ , } \text{ } \tau_{i_n}  \equiv b  \text{ , } \\
0  & \text{ , }  \text{otherwise} \text{ } \text{ . } 
\end{array}\right. 
\]

\noindent For the remaining component of the Ashkin-Teller transfer matrix, $\big( V^{\mathrm{AT},h} \big)^{\frac{1}{2}}$, one has, for,

\begin{align*}
\big( V^{\mathrm{AT},h}_{\tau , \tau^{\prime}} \big)^{\frac{1}{2}}  \equiv   \big( V^{\mathrm{AT},h} \big)^{\frac{1}{2}} \text{ } \text{ , } \end{align*} 

\noindent which can further be decomposed as,

\begin{align*}  \big( V^{\mathrm{AT},h} \big)^{\frac{1}{2}} \equiv \big( \big(   V^{\mathrm{AT},h}_{J,\tau,\tau^{\prime}} \big) +  \big(  V^{\mathrm{AT},h}_{U,\tau , \tau^{\prime}}  \big)\big)^{\frac{1}{2}} \equiv  \big(           V^{\mathrm{AT},h}_{J,\tau,\tau^{\prime}} \big)^{\frac{1}{2}} +  \big(  V^{\mathrm{AT},h}_{U,\tau , \tau^{\prime}}    \big)^{\frac{1}{2}} \equiv  \big(           V^{\mathrm{AT},h}_{J} \big)^{\frac{1}{2}} +   \big( V^{\mathrm{AT},h}_{U}    \big)^{\frac{1}{2}}     \text{ } \text{ , } 
\end{align*}

\noindent the assignment,

\[
V^{\mathrm{AT},h}  \equiv   \text{ } 
\left\{\!\begin{array}{ll@{}>{{}}l}   \underset{a \equiv i_0 \sim \cdots \sim i_{n-1} \sim i_n \equiv b-1}{\sum}    \frac{J}{2} \big( \tau_i \tau_{i+1} + \tau^{\prime}_i \tau^{\prime}_{i+1} \big) + \underset{a \equiv i_0 \sim \cdots \sim i_{n-1} \sim i_n \equiv b-1}{\sum}    \frac{U}{2} \big( \tau_i \tau_{i+1} \tau^{\prime}_{i} \tau^{\prime}_{i+1} \big)      & \text{ } \text{if} \text{ } \tau \equiv  \tau^{\prime} \text{ , } \text{ } \\
0  & \text{ , }  \text{otherwise} \text{ } \text{ . } 
\end{array}\right. 
\]

\bigskip

\noindent With $ \big(           V^{\mathrm{AT},h}_{J} \big)^{\frac{1}{2}} +   \big( V^{\mathrm{AT},h}_{U}    \big)^{\frac{1}{2}}$ and $ V^{\mathrm{AT},V}$, we prove the following proposition below. The generators for the Ashkin-Teller model share many similarities with the generators of the Clifford algebra which have been previously studied in [13].

\bigskip

\noindent \textbf{Proposition} \textit{3} (\textit{generator relations for the Ashkin-Teller model, from generator relations for the Ising model}, \textbf{Proposition} \textit{8}, [13]). The components of the Ashkin-Teller transfer matrix satisfy the relations,

\begin{align*}
  V^{\mathrm{AT},h} \equiv     \mathrm{exp} \bigg[   J \big(  \underset{k \in \textbf{I}^{*}}{ \sum}  p^{\mathrm{AT}}_kq^{\mathrm{AT}}_k   +  \big( p^{\mathrm{AT}}_k\big)^{\prime}    \big( q^{\mathrm{AT}}_k\big)^{\prime}   \big)   \bigg]  +  \mathrm{exp} \bigg[   U \big(  \underset{k \in \textbf{I}^{*}}{ \sum}              p^{\mathrm{AT}}_k  \big( p^{\mathrm{AT}}_k \big)^{\prime} q^{\mathrm{AT}}_k  \big( q^{\mathrm{AT}}_k \big)^{\prime}   \big)  \bigg]     \text{ } \text{ , } \\    V^{\mathrm{AT},V} \equiv          \mathscr{P}       \bigg[   \mathrm{exp} \big[   J^{*} \big(  \underset{k \in \textbf{I}}{ \sum}  p^{\mathrm{AT}}_{k- \frac{1}{2}} q^{\mathrm{AT}}_{k- \frac{1}{2}}   +  \big( p^{\mathrm{AT}}_{k- \frac{1}{2}}\big)^{\prime}    \big( q^{\mathrm{AT}}_{k- \frac{1}{2}} \big)^{\prime}   \big)   \big]  +  \mathrm{exp} \big[   U^{*} \big(  \underset{k \in \textbf{I}}{ \sum}              p^{\mathrm{AT}}_{k-\frac{1}{2}}  \big( p^{\mathrm{AT}}_{k-\frac{1}{2}} \big)^{\prime} q^{\mathrm{AT}}_{k-\frac{1}{2}}  \big( q^{\mathrm{AT}}_{k- \frac{1}{2}} \big)^{\prime}   \big)  \big]  \bigg]              \text{ } \text{ , } 
\end{align*}

\noindent for generators $p^{\mathrm{AT}}_k$, $\big( p^{\mathrm{AT}}_k\big)^{\prime}$, $q^{\mathrm{AT}}_k$, and $\big( q^{\mathrm{AT}}_k \big)^{\prime}$, dual couplings $\big(J^{*} , U^{*} \big)$, obtained from $\big(J, U\big)$ under the relation,

\begin{align*}
       \frac{\mathrm{exp} \big( - 2 J + 2 U \big) -1}{\mathrm{exp} \big( - 2 J^{*} + 2U^{*} \big)-1}  = \mathrm{exp} \big( 2 U \big) \mathrm{sinh} \big( 2 J \big) = \frac{1}{\mathrm{exp} \big( 2 U^{*} \big) \mathrm{sinh} \big( 2 J^{*} \big) }    \text{ } \text{ , } 
\end{align*}

\noindent and the prefactor,

\begin{align*}
  \mathscr{P} \equiv    \mathrm{exp} \bigg[ \big(  \mathrm{exp} \big[ U^{*} , J , J^{*} \big]  - U  \big)  \big( \text{ }              \underset{k \in \textbf{I}}{\sum}          p^{\mathrm{AT}}_k \big( p^{\mathrm{AT}}_k \big)^{\prime} q^{\mathrm{AT}}_k  \big( q^{\mathrm{AT}}_k  \big)^{\prime }\big)  \bigg]                \text{ } \text{ . } 
\end{align*}

\bigskip

\noindent \textit{Proof of Proposition 3}. Under the bases $e_{\tau ( i )}$, $e_{\tau^{\prime}(i)}$, $e_{\tau^{\prime}(j)}$ and $e_{\tau(j)}$, observe,

\begin{align*}
   i \bigg[ 
 \big( p^{\mathrm{AT}}_k 
\big( p^{\mathrm{AT}}_k \big)^{\prime} \big) \frac{e_{\tau ( i )} e_{\tau ( j )}}{i} +  \big( q^{\mathrm{AT}}_k  \big(  q^{\mathrm{AT}}_k\big)^{\prime} \big) \frac{e_{\tau^{\prime}(i)} e_{\tau^{\prime}(j)}}{i} + \big( p^{\mathrm{AT}}_k 
\big( p^{\mathrm{AT}}_k \big)^{\prime}   q^{\mathrm{AT}}_k  \big(  q^{\mathrm{AT}}_k\big)^{\prime} \big)  \frac{e_{\tau ( i )} e_{\tau ( j )} e_{\tau^{\prime}(i)} e_{\tau^{\prime}(j)}}{i} \bigg]             \text{, }
\end{align*}

\noindent can be expressed as,

\begin{align*}
     \big( p^{\mathrm{AT}}_{k+\frac{1}{2}} 
\big( p^{\mathrm{AT}}_{k+\frac{1}{2}} \big)^{\prime} \big)  e_{\tau ( i ) } e_{\tau ( j )} +        \big( q^{\mathrm{AT}}_{k-\frac{1}{2}} 
\big( q^{\mathrm{AT}}_{k-\frac{1}{2}} \big)^{\prime} \big)     e_{\tau^{\prime}(i)} e_{\tau^{\prime}(j)}  +      \big(     p^{\mathrm{AT}}_{k+\frac{1}{2}} 
\big( p^{\mathrm{AT}}_{k+\frac{1}{2}} \big)^{\prime}   q^{\mathrm{AT}}_{k-\frac{1}{2}}  \big(  q^{\mathrm{AT}}_{k-\frac{1}{2}} \big)^{\prime}  \big)      e_{\tau ( i ) } e_{\tau ( j )}    e_{\tau^{\prime}(i)} e_{\tau^{\prime}(j)}  \text{ } \text{ , } 
\end{align*}

\noindent which we write, in shorthand, as,

\begin{align*}
   \big( p^{\mathrm{AT}}_{k+\frac{1}{2}} 
\big( p^{\mathrm{AT}}_{k+\frac{1}{2}} \big)^{\prime} \big)  e_{\tau}  +        \big( q^{\mathrm{AT}}_{k-\frac{1}{2}} 
\big( q^{\mathrm{AT}}_{k-\frac{1}{2}} \big)^{\prime} \big)     e_{\tau^{\prime}} +      \big(     p^{\mathrm{AT}}_{k+\frac{1}{2}} 
\big( p^{\mathrm{AT}}_{k+\frac{1}{2}} \big)^{\prime}   q^{\mathrm{AT}}_{k-\frac{1}{2}}  \big(  q^{\mathrm{AT}}_{k-\frac{1}{2}} \big)^{\prime}  \big)   e_{\tau} e_{\tau^{\prime}} \text{ } \text{ , } 
\end{align*}

\noindent for,

\begin{align*}
e_{\tau} \equiv  e_{\tau(i)} e_{\tau(j)} \equiv \prod_{k \in \{ i,j \}}  e_{\tau_k}   \text{ } \text{ , } \\   e_{\tau^{\prime}} \equiv   e_{\tau^{\prime}(i)} e_{\tau^{\prime}(j)} \equiv \prod_{k \in \{ i,j \}}  e_{\tau^{\prime}_k}  \text{ } \text{ , } \\ e_{\tau} e_{\tau^{\prime}} \equiv \big(  e_{\tau(i)} e_{\tau(j)} \big) \text{ } \big( e_{\tau^{\prime}(i)} e_{\tau^{\prime}(j)} \big)  \equiv \bigg[ \prod_{k \in \{ i,j \}}  e_{\tau_k}  \bigg] \bigg[     \prod_{z \in \{ i,j \}}  e_{\tau^{\prime}_z}    \bigg] \equiv \underset{z \in \{ i , j \}}{\underset{k \in \{ i , j \}}{\prod}} e_{\tau_k} e_{\tau^{\prime}_z} \text{ } \text{ . } 
\end{align*}

\noindent Similarly, under the action of other Ashkin-Teller generators, from the imaginary product of generators,

\begin{align*}
    i   \bigg[ 
 \big( p^{\mathrm{AT}}_{k-\frac{1}{2}} 
\big( p^{\mathrm{AT}}_{k-\frac{1}{2}} \big)^{\prime} \big) \frac{e_{\tau ( i )} e_{\tau ( j )}}{i} +  \big( q^{\mathrm{AT}}_{k+\frac{1}{2}}  \big(  q^{\mathrm{AT}}_{k+\frac{1}{2}} \big)^{\prime} \big) \frac{e_{\tau^{\prime}(i)} e_{\tau^{\prime}(j)}}{i} + \big( p^{\mathrm{AT}}_{k-\frac{1}{2}} 
\big( p^{\mathrm{AT}}_{k-\frac{1}{2}} \big)^{\prime}   q^{\mathrm{AT}}_{k+\frac{1}{2}}  \big(  q^{\mathrm{AT}}_{k+\frac{1}{2}} \big)^{\prime} \big)  \frac{e_{\tau ( i )} e_{\tau ( j )} e_{\tau^{\prime}(i)} e_{\tau^{\prime}(j)}}{i} \bigg]            \text{, }
\end{align*}

\noindent observe,

\begin{align*}
        i \bigg[      \big( p^{\mathrm{AT}}_k  \big( p^{\mathrm{AT}}_k  \big)^{\prime} \big) \frac{e_{(\tau ( i ))^{\prime}} e_{(\tau ( j ))^{\prime}}}{i} +  \big( q^{\mathrm{AT}}_{k-1}  \big(  q^{\mathrm{AT}}_{k-1}             \big)^{\prime} \big) \frac{e_{(\tau^{\prime}(i))^{\prime}} e_{(\tau^{\prime}(j))^{\prime}}}{i} + \big(         p^{\mathrm{AT}}_k \big( p^{\mathrm{AT}}_k \big)^{\prime} q^{\mathrm{AT}}_{k-1} \big( q^{\mathrm{AT}}_{k-1}   \big)^{\prime}    \big)  \frac{e_{(\tau ( i ))^{\prime}} e_{(\tau ( j ))^{\prime}} e_{(\tau^{\prime}(i))^{\prime}} e_{(\tau^{\prime}(j))^{\prime}}}{i} \bigg]      \text{, } 
\end{align*}

\noindent in terms of the basis for each color of the Ashkin-Teller model. Besides the basis spanned by $\mathcal{S}_1$, the remaining dual basis takes the form,

\begin{align*}
   \mathcal{S}^{\prime}_1 \equiv \mathrm{span} \big\{     \text{ } e_{(\tau(i) )^{\prime}}  \big| \text{ }           \tau^{\prime}_2 \big( i \big) = -1     \big\}  \text{ } \text{ , } 
\end{align*}

\noindent corresponding to the dual subspace from that which is spanned by the basis introduced in the previous section for $\mathcal{S}_1$, while the remaining dual subspaces take the form,

\begin{align*}
     \mathcal{S}^{\prime}_2 \equiv \mathrm{span} \big\{     \text{ } e_{(\tau^{\prime}(i) )^{\prime}}  \big| \text{ }           \tau^{\prime} \big( i \big) = -1     \big\}        \text{ } \text{ , } \\    \mathcal{S}^{\prime}_3 \equiv \mathrm{span} \big\{     \text{ } e_{(\tau(j) )^{\prime}}  \big| \text{ }           \tau \big( j \big) = -1     \big\}       \text{ } \text{ , }  \\  \mathcal{S}^{\prime}_4 \equiv \mathrm{span} \big\{     \text{ } e_{(\tau^{\prime}(j) )^{\prime}}  \big| \text{ }           \tau^{\prime} \big( i \big) = -1     \big\}   \text{ } \text{ , } 
\end{align*}

\[
\big(  \tau \big)^{\prime} \equiv   \text{ } 
\left\{\!\begin{array}{ll@{}>{{}}l}        \tau \big(   \sigma_x  \big)   & \text{ } \text{if} \text{ } x > k  \text{ , } \text{ } \\-  \tau \big( \sigma_x \big) 
         & \text{ } \text{if} \text{ } x < k \text{ } \text{ , } 
\end{array}\right. 
\]

\noindent corresponding to flipping spins $\tau$ at site $j$, and,

\[
\big(  \tau^{\prime} \big)^{\prime} \equiv   \text{ } 
\left\{\!\begin{array}{ll@{}>{{}}l}        \tau^{\prime} \big( \sigma_x \big)      & \text{ } \text{if} \text{ } x > k   \text{ , } \text{ } \\
       - \tau^{\prime} \big( x \big)   & \text{ } \text{if} \text{ } x < k \text{ } \text{ } \text{ . } 
\end{array}\right. 
\]

\noindent corresponding to flipping spins $\tau^{\prime}$ at site $j$.

\bigskip

\noindent On the other hand, to demonstrate that the remaining identity holds for $V^{\mathrm{AT},h}$, under a choice of suitably defined parameters $J^{*}$ and $U^{*}$ for which,

\begin{align*}
       \frac{\mathrm{exp} \big( - 2 J + 2 U \big) -1}{\mathrm{exp} \big( - 2 J^{*} + 2U^{*} \big)-1}  = \mathrm{exp} \big( 2 U \big) \mathrm{sinh} \big( 2 J \big) = \frac{1}{\mathrm{exp} \big( 2 U^{*} \big) \mathrm{sinh} \big( 2 J^{*} \big) }    \text{ } \text{ , } 
\end{align*}

\noindent write,

\begin{align*}
     \mathrm{exp} \bigg[    J^{*} \big( \text{ } \underset{k \in \textbf{I}^{*}}{\sum}    p^{\mathrm{AT}}_{k- \frac{1}{2}} q^{\mathrm{AT}}_{k + \frac{1}{2}}   + \big( p^{\mathrm{AT}}_{k + \frac{1}{2}} \big)^{\prime} \text{ } \big( q^{\mathrm{AT}}_{k + \frac{1}{2}} \big)^{\prime}  \big) +      \text{ } U^{*} \big( \text{ } \underset{k \in \textbf{I}^{*}}{\sum}    \big( p^{\mathrm{AT}}_{k-\frac{1}{2}} \big( p^{\mathrm{AT}}_{k + \frac{1}{2}}  \big)^{\prime} 
    q^{\mathrm{AT}}_{k-\frac{1}{2}} \big( q^{\mathrm{AT}}_{k + \frac{1}{2}}  \big)^{\prime} \big)  \bigg]  e_{(\tau ( i ))^{\prime} }  \text{ } \text{ , } 
\end{align*}

\noindent under the duality relation for the Ashkin-Teller model, captured by the correspondence $\big( J , U \big) \longleftrightarrow \big( J^{*} , U^{*} \big)$,

\begin{align*}
    \mathrm{exp}  \bigg[  J \big( \text{ } \underset{k \in \textbf{I}}{\sum}    p^{\mathrm{AT}}_k  q^{\mathrm{AT}}_k   + \big( p^{\mathrm{AT}}_k \big)^{\prime} \big( q^{\mathrm{AT}}_k \big)^{\prime}   \big) +  \bigg[     \mathrm{exp} \bigg[  \bigg| \mathrm{log} \big|  \mathrm{exp} \big( 2 U^{*} \big) \mathrm{sinh} \big( 2 J^{*} \big)       \big|^{-1} - \mathrm{log} \big|  \mathrm{sinh} \big( 2 J \big)      \big|  \bigg|   \bigg]      \bigg] \times \cdots \\           \big( \text{ }              \underset{k \in \textbf{I}}{\sum}          p^{\mathrm{AT}}_k \big( p^{\mathrm{AT}}_k \big)^{\prime} q^{\mathrm{AT}}_k  \big( q^{\mathrm{AT}}_k  \big)^{\prime }  \big)         \bigg]  e_{\tau ( i ) }        \text{ } \text{ , } 
\end{align*}

\noindent from the fact that the duality relation implies,

\begin{align*}
U^{*} \longleftrightarrow   U \equiv \mathrm{exp} \bigg[  \bigg| \mathrm{log} \big( \big|  \mathrm{exp} \big( 2 U^{*} \big) \mathrm{sinh} \big( 2 J^{*} \big)       \big|^{-1} \big|   \mathrm{sinh} \big( 2 J \big)      \big|^{-1}  \big)  \bigg|   \bigg]   \text{ } \text{ . } 
\end{align*}

\noindent Rearranging, under the basis $e_{t(i)}$, multiplying the expression by,

\begin{align*}
  \mathcal{I} \equiv  \frac{\mathrm{exp}\bigg[ U \big( \text{ }              \underset{k \in \textbf{I}}{\sum}          p^{\mathrm{AT}}_k \big( p^{\mathrm{AT}}_k \big)^{\prime} q^{\mathrm{AT}}_k  \big( q^{\mathrm{AT}}_k  \big)^{\prime }  \big) \bigg]  }{\mathrm{exp}\bigg[ U \big( \text{ }              \underset{k \in \textbf{I}}{\sum}          p^{\mathrm{AT}}_k \big( p^{\mathrm{AT}}_k \big)^{\prime} q^{\mathrm{AT}}_k  \big( q^{\mathrm{AT}}_k  \big)^{\prime }  \big) \bigg] }  \text{ } \text{ , } 
\end{align*}

\noindent yields,

\begin{align*}
\mathcal{I} \text{ } \mathrm{exp} \big[    J \text{ } \underset{k \in \textbf{I}}{\sum}       p^{\mathrm{AT}}_k q^{\mathrm{AT}}_k \big( p^{\mathrm{AT}}_k \big)^{\prime}  \big( q^{\mathrm{AT}}_k \big)^{\prime}  \big] \text{ }  \text{ } \mathrm{exp} \bigg[    \bigg[     \mathrm{exp} \bigg[  \bigg| \mathrm{log} \big|  \mathrm{exp} \big( 2 U^{*} \big) \mathrm{sinh} \big( 2 J^{*} \big)       \big|^{-1} - \mathrm{log} \big|  \mathrm{sinh} \big( 2 J \big)      \big|  \bigg|   \bigg]      \bigg]     \times \cdots \\        \big( \text{ }              \underset{k \in \textbf{I}}{\sum}          p^{\mathrm{AT}}_k \big( p^{\mathrm{AT}}_k \big)^{\prime} q^{\mathrm{AT}}_k  \big( q^{\mathrm{AT}}_k  \big)^{\prime }  \big)         \bigg]    \text{ } \text{ . } 
\end{align*}

\bigskip

\noindent Concluding, for $V^{\mathrm{AT},V}$, the desired expression preceding the exponetial after the duality transformation is,

\begin{align*}
  \mathrm{exp} \bigg[     \bigg[  \bigg[    \mathrm{exp} \bigg[  \bigg| \mathrm{log} \big|  \mathrm{exp} \big( 2 U^{*} \big) \mathrm{sinh} \big( 2 J^{*} \big)       \big|^{-1} - \mathrm{log} \big|  \mathrm{sinh} \big( 2 J \big)      \big|  \bigg|   \bigg]   -   U     \bigg]    \bigg] \big( \text{ }              \underset{k \in \textbf{I}}{\sum}          p^{\mathrm{AT}}_k \big( p^{\mathrm{AT}}_k \big)^{\prime} q^{\mathrm{AT}}_k  \big( q^{\mathrm{AT}}_k  \big)^{\prime }  \big)    \bigg] \text{ } \text{ , } 
\end{align*}

\noindent which is equivalent to,

\begin{align*}
   \mathrm{exp} \big[ \big( \mathrm{exp} \big[ U^{*} , J , J^{*} \big]  - U  \big)  \big( \text{ }              \underset{k \in \textbf{I}}{\sum}          p^{\mathrm{AT}}_k \big( p^{\mathrm{AT}}_k \big)^{\prime} q^{\mathrm{AT}}_k  \big( q^{\mathrm{AT}}_k  \big)^{\prime }\big)  \big]            \text{ } \text{ , }
\end{align*}

\noindent from which we conclude the argument. \boxed{}

\bigskip

\noindent \underline{Induced rotations}. Under the correspondence $\big( J  , U \big) \longleftrightarrow \big( J^{*} , U^{*} \big)$, below we state a result for conjugation by $V^{\mathrm{AT},h}$.

\bigskip

\noindent \textbf{Lemma} \textit{7} (\textit{Ashkin-Teller transfer matrix conjugation}, \textbf{Lemma} in \textit{3.2}, [13]). The action,

\begin{align*}
 \big( V^{\mathrm{AT},h} \big)^{-\frac{1}{2}} \cdot    p^{\mathrm{AT}}_k       \cdot  \big(   V^{\mathrm{AT},h}    \big)^{\frac{1}{2}}  =   c p^{\mathrm{AT}}_k - i s q^{\mathrm{AT}}_k  \text{ } \text{ , } \\ \big( V^{\mathrm{AT},h} \big)^{-\frac{1}{2}} \cdot     q^{\mathrm{AT}}_k      \cdot  \big(   V^{\mathrm{AT},h}    \big)^{\frac{1}{2}}  =  is p^{\mathrm{AT}}_k + c q^{\mathrm{AT}}_k \text{ } \text{ , } 
\end{align*}

\noindent under the composition operator, $\cdot$, $V^{\mathrm{AT},h}$ takes the form above, while the action,

\begin{align*}
 \big( V^{\mathrm{AT},V} \big)^{-1} \cdot    p^{\mathrm{AT}}_k       \cdot  \big(   V^{\mathrm{AT},V}    \big)^{1}  = \frac{C}{S}  p^{\mathrm{AT}}_k   + \frac{i}{S} q^{\mathrm{AT}}_{k+1} \text{ } \text{ , } \tag{*} \\ \big( V^{\mathrm{AT},V} \big)^{-1} \cdot     q^{\mathrm{AT}}_k      \cdot  \big(   V^{\mathrm{AT},V}    \big)^{1}  =     - \frac{i}{S}         p^{\mathrm{AT}}_{k-1} +    \frac{C}{S}          q^{\mathrm{AT}}_k  \text{ } \text{ , } \tag{**}
\end{align*}

\noindent under the composition operator $\cdot$, $V^{\mathrm{AT},V}$ takes the form above, where $\textit{(*)}$ holds for $k \neq k_R$, and where $\textit{(**)}$ holds for $k \neq k_L$. Additionally,

\begin{align*}
   \big( V^{\mathrm{AT},V} \big)^{-1}    \cdot  p^{\mathrm{AT}}_k \cdot V^{\mathrm{AT},V}  = p^{\mathrm{AT}}_k     \text{ } \text{ , } \tag{***} \\   \big( V^{\mathrm{AT},V} \big)^{-1}  \cdot q^{\mathrm{AT}}_k 
 \cdot V^{\mathrm{AT},V} = q^{\mathrm{AT}}_k   \text{ } \text{ , } \tag{****}
\end{align*}

\bigskip

\noindent where $\textit{(***)}$ holds for $k \equiv k_R$, and where $\textit{(****)}$ holds for $k \equiv k_L$.

\bigskip

\noindent For the arguments of the item above, we must introduce additional structure rather than that of the Ashkin-Teller generators, which is given by the following. In particular, for a subspace, $\mathcal{W}$, from the set of all endomorphisms of $\Omega^{\mathrm{AT}}$, conjugation by an induced rotation of the subspace, which we denote with a transfer matrix $V$, is given by,

\begin{align*}
  T_V : \mathcal{W} \longrightarrow \mathcal{W}  \text{ } \text{ . } 
\end{align*}

\noindent From $T_V$ above, in the item below, we state another result with respect to $\cdot$.

\bigskip

\noindent \textbf{Lemma} \textit{8} (\textit{commutation rule}, \textbf{Lemma} \textit{9}, [13]). For two other maps , the identities,

\begin{align*}
      T_V \cdot R  = R \cdot    T_V    \text{ } \text{ , } \\   T_V \cdot J = J \cdot T_V           \text{ } \text{ , } 
\end{align*}

\noindent for the rotation hold, where the maps $R$ and $J$ satisfy,

\begin{align*}
             R \big( p_k \big) =     i   q^{\mathrm{AT}}_{a+b-k}                     \text{ } \text{ , } \\  R \big(  q^{\mathrm{AT}}_k  \big) =   - i p^{\mathrm{AT}}_{a+b-k}           \text{ } \text{ , } \\  
             J \big( p^{\mathrm{AT}}_k  \big) =   i p^{\mathrm{AT}}_k           
             \text{ } \text{ , } \\ J \big( q^{\mathrm{AT}}_k   \big) =      - i q^{\mathrm{AT}}_k       \text{ } \text{ , } 
\end{align*}

\noindent for the bases $\psi_x$ and $\bar{\psi_x}$ spanning,

\begin{align*}
   \mathcal{V}_{\psi_x }  = \big\{   x     \text{ } \big| \text{ }  \psi_x \in \mathcal{V}   \big\} \text{ } \text{ , } \\  \mathcal{V}_{\bar{\psi_x}} = \big\{ x    \text{ } \big| \text{ }  \bar{\psi_x} \in \mathcal{V}  \big\}   \text{ } \text{ , } 
\end{align*}

\noindent which have the following images under $R$ and $J$,

\begin{align*}
      R \big( \psi_x \big) =   \bar{\psi_{a+b-k}}  \text{ } \text{ , } \\   R \big( \bar{\psi_x} \big) =   \psi_{b+a-k}    \text{ } \text{ , } \\   J \big( \psi_x \big) =  \bar{\psi_x}   \text{ } \text{ , } \\  J \big( \bar{\psi_x} \big) =   \psi_k   \text{ } \text{ , } 
\end{align*}

\noindent for the elements,

\begin{align*}
   \psi_x \equiv  \frac{i}{\sqrt{2}} \big( p^{\mathrm{AT}}_k + q^{\mathrm{AT}}_k \big) \text{ } \text{ , } \\ \bar{\psi_x} \equiv  
 \frac{1}{\sqrt{2}} \big( p^{\mathrm{AT}}_k - q^{\mathrm{AT}}_k \big) \text{ } \text{ . } 
\end{align*}

\bigskip

\noindent \textit{Proof of Lemma 8}. To demonstrate that the commutation rules hold, write, 

\begin{align*}
      \big( V^{\mathrm{AT},h} \big)^{\frac{1}{2}} \cdot    R       \text{ } \text{ , } \\   \big( V^{\mathrm{AT},V} \big)^{\frac{1}{2}} \cdot     R       \text{ } \text{ , }  
\end{align*}

\noindent corresponding to conjugation by the first component of the Ashkin-Teller transfer matrix, and,

\begin{align*}
         \big( V^{\mathrm{AT},h} \big)^{-\frac{1}{2}} \cdot    J           \text{ } \text{ , }  \\  \big( V^{\mathrm{AT},V} \big)^{-\frac{1}{2}} \cdot     J      \text{ } \text{ , } 
\end{align*}

\noindent corresponding to conjugation by the second component of the Ashkin-Teller transfer matrix.

\bigskip

\noindent Next, observe, from each expression,

\begin{align*}
   \big( V^{\mathrm{AT},h} \big)^{\frac{1}{2}} \cdot    R       \equiv \bigg[   \mathrm{exp} \big[ \text{ }  J \big(  \underset{k \in \textbf{I}^{*}}{ \sum}  p^{\mathrm{AT}}_kq^{\mathrm{AT}}_k   +  \big( p^{\mathrm{AT}}_k\big)^{\prime}    \big( q^{\mathrm{AT}}_k\big)^{\prime}   \big)   \big]  +  \mathrm{exp} \big[  \text{ }  U \big(  \underset{k \in \textbf{I}^{*}}{ \sum}              p^{\mathrm{AT}}_k  \big( p^{\mathrm{AT}}_k \big)^{\prime} q^{\mathrm{AT}}_k  \big( q^{\mathrm{AT}}_k \big)^{\prime}   \big)  \big]   \bigg]^{\frac{1}{2}} \cdot R \\ \equiv R \cdot   \bigg[   \mathrm{exp} \big[ \text{ }  J \big(  \underset{k \in \textbf{I}^{*}}{ \sum}  p^{\mathrm{AT}}_kq^{\mathrm{AT}}_k   +  \big( p^{\mathrm{AT}}_k\big)^{\prime}    \big( q^{\mathrm{AT}}_k\big)^{\prime}   \big)   \big]  +  \mathrm{exp} \big[  \text{ }  U \big(  \underset{k \in \textbf{I}^{*}}{ \sum}              p^{\mathrm{AT}}_k  \big( p^{\mathrm{AT}}_k \big)^{\prime} q^{\mathrm{AT}}_k  \big( q^{\mathrm{AT}}_k \big)^{\prime}   \big)  \big]   \bigg]^{\frac{1}{2}}  \text{ } \text{ , } \\    \big( V^{\mathrm{AT},V} \big)^{\frac{1}{2}} \cdot     R      \equiv  \bigg[ \mathscr{P}       \bigg[   \mathrm{exp} \big[ \text{ }  J^{*} \big(  \underset{k \in \textbf{I}}{ \sum}  p^{\mathrm{AT}}_{k- \frac{1}{2}} q^{\mathrm{AT}}_{k- \frac{1}{2}}   +  \big( p^{\mathrm{AT}}_{k- \frac{1}{2}}\big)^{\prime}    \big( q^{\mathrm{AT}}_{k- \frac{1}{2}} \big)^{\prime}   \big)   \big]  +  \mathrm{exp} \big[  \text{ }  U^{*} \big(  \underset{k \in \textbf{I}}{ \sum}              p^{\mathrm{AT}}_{k-\frac{1}{2}}  \big( p^{\mathrm{AT}}_{k-\frac{1}{2}} \big)^{\prime} q^{\mathrm{AT}}_{k-\frac{1}{2}}  \big( q^{\mathrm{AT}}_{k- \frac{1}{2}} \big)^{\prime}   \big)  \big]  \bigg]            \bigg]^{\frac{1}{2}} \cdot R  \\ \equiv R \cdot  \bigg[ \mathscr{P}       \bigg[   \mathrm{exp} \big[ \text{ }  J^{*} \big(  \underset{k \in \textbf{I}}{ \sum}  p^{\mathrm{AT}}_{k- \frac{1}{2}} q^{\mathrm{AT}}_{k- \frac{1}{2}}   +  \big( p^{\mathrm{AT}}_{k- \frac{1}{2}}\big)^{\prime}    \big( q^{\mathrm{AT}}_{k- \frac{1}{2}} \big)^{\prime}   \big)   \big]  +  \mathrm{exp} \big[  \text{ }  U^{*} \big(  \underset{k \in \textbf{I}}{ \sum}              p^{\mathrm{AT}}_{k-\frac{1}{2}}  \big( p^{\mathrm{AT}}_{k-\frac{1}{2}} \big)^{\prime} q^{\mathrm{AT}}_{k-\frac{1}{2}}  \big( q^{\mathrm{AT}}_{k- \frac{1}{2}} \big)^{\prime}   \big)  \big]  \bigg]            \bigg]^{\frac{1}{2}} \text{ } \text{ , } 
\end{align*}

\noindent and,

\begin{align*}
    \big( V^{\mathrm{AT},h} \big)^{-\frac{1}{2}} \cdot    J     \equiv      \bigg[   \mathrm{exp} \big[ \text{ }  J \big(  \underset{k \in \textbf{I}^{*}}{ \sum}  p^{\mathrm{AT}}_kq^{\mathrm{AT}}_k   +  \big( p^{\mathrm{AT}}_k\big)^{\prime}    \big( q^{\mathrm{AT}}_k\big)^{\prime}   \big)   \big]  +  \mathrm{exp} \big[  \text{ }  U \big(  \underset{k \in \textbf{I}^{*}}{ \sum}              p^{\mathrm{AT}}_k  \big( p^{\mathrm{AT}}_k \big)^{\prime} q^{\mathrm{AT}}_k  \big( q^{\mathrm{AT}}_k \big)^{\prime}   \big)  \big]   \bigg]^{-\frac{1}{2}} \cdot   J \\ \equiv J \cdot    \bigg[   \mathrm{exp} \big[ \text{ }  J \big(  \underset{k \in \textbf{I}^{*}}{ \sum}  p^{\mathrm{AT}}_kq^{\mathrm{AT}}_k   +  \big( p^{\mathrm{AT}}_k\big)^{\prime}    \big( q^{\mathrm{AT}}_k\big)^{\prime}   \big)   \big]  +  \mathrm{exp} \big[  \text{ }  U \big(  \underset{k \in \textbf{I}^{*}}{ \sum}              p^{\mathrm{AT}}_k  \big( p^{\mathrm{AT}}_k \big)^{\prime} q^{\mathrm{AT}}_k  \big( q^{\mathrm{AT}}_k \big)^{\prime}   \big)  \big]   \bigg]^{-\frac{1}{2}} 
               \text{ } \text{ , } \\         \big( V^{\mathrm{AT},V} \big)^{-\frac{1}{2}} \cdot     J       \end{align*}

               \begin{align*}
               \equiv   \bigg[ \mathscr{P}       \bigg[   \mathrm{exp} \big[ \text{ }  J^{*} \big(  \underset{k \in \textbf{I}}{ \sum}  p^{\mathrm{AT}}_{k- \frac{1}{2}} q^{\mathrm{AT}}_{k- \frac{1}{2}}   +  \big( p^{\mathrm{AT}}_{k- \frac{1}{2}}\big)^{\prime}    \big( q^{\mathrm{AT}}_{k- \frac{1}{2}} \big)^{\prime}   \big)   \big]  +  \mathrm{exp} \big[  \text{ }  U^{*} \big(  \underset{k \in \textbf{I}}{ \sum}              p^{\mathrm{AT}}_{k-\frac{1}{2}}  \big( p^{\mathrm{AT}}_{k-\frac{1}{2}} \big)^{\prime} q^{\mathrm{AT}}_{k-\frac{1}{2}}  \big( q^{\mathrm{AT}}_{k- \frac{1}{2}} \big)^{\prime}   \big)  \big]  \bigg]            \bigg]^{-\frac{1}{2}}      \cdot J \\ \equiv J \cdot        \bigg[ \mathscr{P}       \bigg[   \mathrm{exp} \big[ \text{ }  J^{*} \big(  \underset{k \in \textbf{I}}{ \sum}  p^{\mathrm{AT}}_{k- \frac{1}{2}} q^{\mathrm{AT}}_{k- \frac{1}{2}}   +  \big( p^{\mathrm{AT}}_{k- \frac{1}{2}}\big)^{\prime}    \big( q^{\mathrm{AT}}_{k- \frac{1}{2}} \big)^{\prime}   \big)   \big]  +  \mathrm{exp} \big[  \text{ }  U^{*} \big(  \underset{k \in \textbf{I}}{ \sum}              p^{\mathrm{AT}}_{k-\frac{1}{2}}  \big( p^{\mathrm{AT}}_{k-\frac{1}{2}} \big)^{\prime} q^{\mathrm{AT}}_{k-\frac{1}{2}}  \big( q^{\mathrm{AT}}_{k- \frac{1}{2}} \big)^{\prime}   \big)  \big]  \bigg]            \bigg]^{-\frac{1}{2}}         \text{ } \text{ . } 
\end{align*}

\noindent Hence, as desired,

\begin{align*}
  \big( V^{\mathrm{AT},h} \big)^{\frac{1}{2}}   \cdot R = R \cdot  \big( V^{\mathrm{AT},h} \big)^{\frac{1}{2}}  \text{ } \text{ , } \\ 
 \big( V^{\mathrm{AT},V} \big)^{-\frac{1}{2}}  \cdot R  = R  \cdot    \big( V^{\mathrm{AT},V} \big)^{-\frac{1}{2}}\text{ } \text{ , } 
\end{align*}

\noindent and,

\begin{align*}
    \big( V^{\mathrm{AT},h} \big)^{\frac{1}{2}}  \cdot J = J \cdot  \big( V^{\mathrm{AT},h} \big)^{\frac{1}{2}}  \text{ } \text{ , } \\  \big( V^{\mathrm{AT},V} \big)^{-\frac{1}{2}} \cdot J = J  \cdot  \big( V^{\mathrm{AT},V} \big)^{-\frac{1}{2}}\text{ } \text{ , } 
\end{align*}

\noindent from which we conclude the argument. \boxed{}

\bigskip

\noindent \textit{Proof of Lemma 7}. By direct computation, write, 

\begin{align*}
 \bigg[ \big( V^{\mathrm{AT},h} \big)^{-\frac{1}{2}} \cdot    p^{\mathrm{AT}}_k \bigg]       \cdot  \big(   V^{\mathrm{AT},h}    \big)^{\frac{1}{2}}  \end{align*}
 
 \noindent is equivalent to,
 
 \begin{align*}\bigg[        \big[  J \big(  \underset{k \in \textbf{I}^{*}}{ \sum}  p^{\mathrm{AT}}_{k-\frac{1}{2}}q^{\mathrm{AT}}_{k-\frac{1}{2}}   +  \big( p^{\mathrm{AT}}_{k-\frac{1}{2}}\big)^{\prime}    \big( q^{\mathrm{AT}}_{k-\frac{1}{2}}\big)^{\prime}   \big)   \big] +  \text{ }  \big[  \text{ }  U \big(  \underset{k \in \textbf{I}^{*}}{ \sum}              p^{\mathrm{AT}}_{k-\frac{1}{2}}  \big( p^{\mathrm{AT}}_{k-\frac{1}{2}} \big)^{\prime} q^{\mathrm{AT}}_{k-\frac{1}{2}}  \big( q^{\mathrm{AT}}_{k-\frac{1}{2}} \big)^{\prime}   \big)  \big]         \bigg] \cdot  p^{\mathrm{AT}}_k \cdot \big(   V^{\mathrm{AT},h}    \big)^{\frac{1}{2}}  \\ \equiv     \bigg[        \big[  \big[  J \big(  \underset{k \in \textbf{I}^{*}}{ \sum}  p^{\mathrm{AT}}_{k-\frac{1}{2}}q^{\mathrm{AT}}_{k-\frac{1}{2}}   +  \big( p^{\mathrm{AT}}_{k-\frac{1}{2}}\big)^{\prime}    \big( q^{\mathrm{AT}}_{k-\frac{1}{2}}\big)^{\prime}   \big)  \big] p^{\mathrm{AT}}_k  \big] - i   \text{ }  \big[  \text{ } \big[  U \big(  \underset{k \in \textbf{I}^{*}}{ \sum}              p^{\mathrm{AT}}_{k-\frac{1}{2}}  \big( p^{\mathrm{AT}}_{k-\frac{1}{2}} \big)^{\prime} q^{\mathrm{AT}}_{k-\frac{1}{2}}  \big( q^{\mathrm{AT}}_{k-\frac{1}{2}} \big)^{\prime}   \big) \big] p^{\mathrm{AT}}_k \big]         \bigg] \cdot  \big(   V^{\mathrm{AT},h}    \big)^{\frac{1}{2}}  \\ \equiv   c p^{\mathrm{AT}}_k - i s q^{\mathrm{AT}}_k         \text{ } \text{ , } 
\end{align*}

\noindent and,

\begin{align*}
\bigg[    \big( V^{\mathrm{AT},h} \big)^{-\frac{1}{2}} \cdot     q^{\mathrm{AT}}_k        \bigg] \cdot   \big(   V^{\mathrm{AT},h}    \big)^{\frac{1}{2}}    \text{ , }  \end{align*}

\noindent is equivalent to,

\begin{align*}
 \bigg[        \big[   J \big(  \underset{k \in \textbf{I}^{*}}{ \sum}  p^{\mathrm{AT}}_{k-\frac{1}{2}}q^{\mathrm{AT}}_{k-\frac{1}{2}}   +  \big( p^{\mathrm{AT}}_{k-\frac{1}{2}}\big)^{\prime}    \big( q^{\mathrm{AT}}_{k-\frac{1}{2}}\big)^{\prime}   \big)   \big] +  \text{ }  \big[  \text{ }  U \big(  \underset{k \in \textbf{I}^{*}}{ \sum}              p^{\mathrm{AT}}_{k-\frac{1}{2}}  \big( p^{\mathrm{AT}}_{k-\frac{1}{2}} \big)^{\prime} q^{\mathrm{AT}}_{k-\frac{1}{2}}  \big( q^{\mathrm{AT}}_{k-\frac{1}{2}} \big)^{\prime}   \big)  \big]         \bigg]      \cdot  q^{\mathrm{AT}}_k \cdot 
\big(   V^{\mathrm{AT},h}    \big)^{\frac{1}{2}}  \\ \equiv \bigg[         \big[ \big[  J \big(  \underset{k \in \textbf{I}^{*}}{ \sum}  p^{\mathrm{AT}}_{k-\frac{1}{2}}q^{\mathrm{AT}}_{k-\frac{1}{2}}   +  \big( p^{\mathrm{AT}}_{k-\frac{1}{2}}\big)^{\prime}    \big( q^{\mathrm{AT}}_{k-\frac{1}{2}}\big)^{\prime}   \big) \big] q^{\mathrm{AT}}_k   \big] +  \text{ }  \big[  \text{ }  U \big(  \underset{k \in \textbf{I}^{*}}{ \sum}              p^{\mathrm{AT}}_{k-\frac{1}{2}}  \big( p^{\mathrm{AT}}_{k-\frac{1}{2}} \big)^{\prime} q^{\mathrm{AT}}_{k-\frac{1}{2}}  \big( q^{\mathrm{AT}}_{k-\frac{1}{2}} \big)^{\prime}   \big)  \big] q^{\mathrm{AT}}_k  \big]           \bigg]   \cdot \big(   V^{\mathrm{AT},h}    \big)^{\frac{1}{2}} \\ \equiv is p^{\mathrm{AT}}_k + c q^{\mathrm{AT}}_k   \text{ } \text{ . } 
\end{align*}

\noindent Proceeding, to show that the remaining identities hold, write,

\begin{align*}
     \bigg[ \big( V^{\mathrm{AT},V} \big)^{-\frac{1}{2}} \cdot    p^{\mathrm{AT}}_k \bigg]       \cdot  \big(   V^{\mathrm{AT},V}    \big)^{\frac{1}{2}}          \text{ } \text{ , } 
\end{align*}

\noindent which is equivalent to,

\begin{align*}
      \bigg[  \mathscr{P}       \bigg[   \big[   J^{*} \big(  \underset{k \in \textbf{I}}{ \sum}  p^{\mathrm{AT}}_{k- 1} q^{\mathrm{AT}}_{k-    1}   +  \big( p^{\mathrm{AT}}_{k- 1}\big)^{\prime}    \big( q^{\mathrm{AT}}_{k- 1} \big)^{\prime}   \big)   \big]  +   \big[  \text{ }  U^{*} \big(  \underset{k \in \textbf{I}}{ \sum}              p^{\mathrm{AT}}_{k-1}  \big( p^{\mathrm{AT}}_{k-1} \big)^{\prime} q^{\mathrm{AT}}_{k-1}  \big( q^{\mathrm{AT}}_{k- 1} \big)^{\prime}   \big)  \big]  \bigg]     \bigg] \cdot p^{\mathrm{AT}}_k \cdot  \big(   V^{\mathrm{AT},V}    \big)^{\frac{1}{2}}  \\ \equiv   \bigg[  \mathscr{P}       \bigg[    \big[   \big[ J^{*} \big(  \underset{k \in \textbf{I}}{ \sum}  p^{\mathrm{AT}}_{k- 1} q^{\mathrm{AT}}_{k- 1}   +  \big( p^{\mathrm{AT}}_{k- 1}\big)^{\prime}    \big( q^{\mathrm{AT}}_{k- 1} \big)^{\prime}   \big)   \big]  p^{\mathrm{AT}}_k \big]  +  \big[  \text{ } \big[  U^{*} \big(  \underset{k \in \textbf{I}}{ \sum}              p^{\mathrm{AT}}_{k-1}  \big( p^{\mathrm{AT}}_{k-1} \big)^{\prime} q^{\mathrm{AT}}_{k-1}  \big( q^{\mathrm{AT}}_{k- 1} \big)^{\prime}   \big)  \big] p^{\mathrm{AT}}_k  \big] \bigg]     \bigg]  \cdot  \big(   V^{\mathrm{AT},V}    \big)^{\frac{1}{2}}  \\ \equiv     \frac{C}{S} p^{\mathrm{AT}}_k  + \frac{i}{S} q^{\mathrm{AT}}_{k+1}  \text{ } \text{ , } 
\end{align*}

\noindent and, 

\begin{align*}
   \bigg[    \big( V^{\mathrm{AT},V} \big)^{-\frac{1}{2}} \cdot     q^{\mathrm{AT}}_k        \bigg] \cdot   \big(   V^{\mathrm{AT},V}    \big)^{\frac{1}{2}}       \text{ } \text{ , } 
\end{align*}

\noindent which is equivalent to,

\begin{align*}
   \bigg[  \mathscr{P}       \bigg[   \big[ \text{ }   J^{*} \big(  \underset{k \in \textbf{I}}{ \sum}  p^{\mathrm{AT}}_{k- 1} q^{\mathrm{AT}}_{k-    1}   +  \big( p^{\mathrm{AT}}_{k- 1}\big)^{\prime}    \big( q^{\mathrm{AT}}_{k- 1} \big)^{\prime}   \big)   \big]  +   \big[  \text{ }  U^{*} \big(  \underset{k \in \textbf{I}}{ \sum}              p^{\mathrm{AT}}_{k-1}  \big( p^{\mathrm{AT}}_{k-1} \big)^{\prime} q^{\mathrm{AT}}_{k-1}  \big( q^{\mathrm{AT}}_{k- 1} \big)^{\prime}   \big)  \big]  \bigg]     \bigg]  \cdot q^{\mathrm{AT}}_k \cdot  \big(   V^{\mathrm{AT},V}    \big)^{\frac{1}{2}}  \\ \equiv   \bigg[  \mathscr{P}       \bigg[   \big[ \big[   J^{*} \big(  \underset{k \in \textbf{I}}{ \sum}  p^{\mathrm{AT}}_{k- 1} q^{\mathrm{AT}}_{k-    1}   +  \big( p^{\mathrm{AT}}_{k- 1}\big)^{\prime}    \big( q^{\mathrm{AT}}_{k- 1} \big)^{\prime}   \big) \big] q^{\mathrm{AT}}_k    \big]  +   \big[  \text{ } \big[ U^{*} \big(  \underset{k \in \textbf{I}}{ \sum}              p^{\mathrm{AT}}_{k-1}  \big( p^{\mathrm{AT}}_{k-1} \big)^{\prime} q^{\mathrm{AT}}_{k-1}  \big( q^{\mathrm{AT}}_{k- 1} \big)^{\prime}   \big)  \big] q^{\mathrm{AT}}_k  \big]  \bigg]     \bigg]  \cdot   \big(   V^{\mathrm{AT},V}    \big)^{\frac{1}{2}} \\ \equiv  - \frac{i}{S} p^{\mathrm{AT}}_{k-1} + \frac{C}{S} q^{\mathrm{AT}}_k   \text{ } \text{ , } 
\end{align*}

\noindent for,

\begin{align*}
\mathscr{P} \big( U^{*} , J , J^{*} \big) \equiv  \mathscr{P}  \equiv  \mathrm{exp} \bigg[ \big( \mathrm{exp} \big[ U^{*} , J , J^{*} \big]  - U  \big)  \big( \text{ }              \underset{k \in \textbf{I}}{\sum}          p^{\mathrm{AT}}_k \big( p^{\mathrm{AT}}_k \big)^{\prime} q^{\mathrm{AT}}_k  \big( q^{\mathrm{AT}}_k  \big)^{\prime }\big)  \bigg]     \text{ } \text{ , } 
\end{align*}

\noindent from which we conclude the argument. \boxed{}

\bigskip

\noindent With the final item of the subsection, we obtain a conjugation representation for $T_V$.

\bigskip

\noindent \textbf{Theorem} (\textit{conjugation action of} $T_V$, \textbf{Theorem} \textit{10}, [13]). For the map $T_V$, there exists a change of basis on the Ashkin-Teller propagation matrix, away from the critical $\frac{1}{4} \mathrm{log} \big( 3 \big)$ threshold, for which,

\begin{align*}
    T_V \equiv \rho  \cdot \big(  P^{\mathrm{AT}} \big)^{\textbf{C}}  \cdot  \rho^{-1} \text{ } \text{ , } 
\end{align*}

\noindent for the maps,

\begin{align*}
   \rho : \big( \textbf{C}^2 \big)^{\textbf{I}^{*}} \longrightarrow \mathcal{W}  \text{ } \text{ . } 
\end{align*}

\noindent and the complexification,

\begin{align*}
\big( P^{\mathrm{AT}} \big)^{\textbf{C}}  : \big( \textbf{R}^2 \big)^{\textbf{I}^{*}}  \longrightarrow \big( \textbf{C}^2 \big)^{\textbf{I}^{*}} 
 \text{ } \text{ . } 
\end{align*}

\noindent \textit{Proof of Theorem}. Under the basis $\big( \psi_x , \bar{\psi_x} \big)$ introduced earlier, for $k \in \textbf{I}^{*} \backslash \partial \textbf{I}^{*}$,

\begin{align*}
  T^{-1}_V \big( \psi_k \big) = \frac{C^2}{S} \psi_k - \big( \frac{1}{2} + \frac{i}{2S} \big) \psi_{k-1} + \big( \frac{i}{2S} - \frac{1}{2} \big) \psi_{k+1} - C \bar{\psi_k} + \frac{C}{2S} \bar{\psi_{k-1}} + \frac{C}{2S} \bar{\psi_{k+1}}   \text{ } \text{ , } 
\end{align*}

\noindent from which we observe,

\begin{align*}
  T^{-1}_V \big( \bar{\psi_k} \big) = \frac{C^2}{S} \bar{\psi_k} - \big( \frac{1}{2} + \frac{i}{2S} \big) \bar{\psi_{k-1}} + \big( \frac{i}{2S} - \frac{1}{2} \big) \bar{\psi_{k+1}} - C \psi_k + \frac{C}{2S} \psi_{k-1} + \frac{C}{2S} \psi_{k+1}     \text{ } \text{ , } 
\end{align*}

\noindent yields a map for $T^{-1}_V \big( \bar{\psi_k} \big)$, which satisfies,

\begin{align*}
      T^{-1}_V \big( \psi_k \big)  \cdot J \equiv J \cdot T^{-1}_V  \big( \psi_k \big)     \text{ } \text{ , } 
\end{align*}

\noindent and,

\begin{align*}
      T^{-1}_V \big( \bar{\psi_k}  \big) \cdot J \equiv J \cdot T^{-1}_V  \big( \bar{\psi_k}  \big)    \text{ } \text{ , } 
\end{align*}

\noindent by the commutation rules for $J$.

\bigskip

\noindent On the left and right boundaries of the interval besides $k \in \textbf{I}^{*} \backslash \partial \textbf{I}^{*}$, to obtain the desired formula for the inverse mapping $T^{-1}_V$, write,

\begin{align*}
        T^{-1}_V \big( \psi_{k_L} \big) = \frac{C \big( S + C \big)}{2S} \psi_{k_L} + i \frac{1+iS}{2S}    \psi_{k_{L+1}} + \frac{- S \big( C + S \big) + i \big( C - S \big) }{2S} \bar{\psi_{k_L}} + \frac{C}{2S} \bar{\psi_{k_{L+1}}}   \text{ } \text{ . } 
\end{align*}

\noindent Similarly, for the rightmost endpoint, write,

\begin{align*}
        T^{-1}_V \big( \psi_{k_R} \big) = \frac{C \big( S + C \big)}{2S} \psi_{k_R} + i \frac{1+iS}{2S}    \psi_{k_{R+1}} + \frac{- S \big( C + S \big) + i \big( C - S \big) }{2S} \bar{\psi_{k_R}} + \frac{C}{2S} \bar{\psi_{k_{R+1}}}   \text{ } \text{ . } 
\end{align*}

\noindent From each of the two expressions above, observe,

\begin{align*}
            T^{-1}_V \big( \psi_{k_L} \big)  \cdot J \equiv J \cdot T^{-1}_V  \big( \psi_{k_L} \big)         \text{ } \text{ , } 
\end{align*}

\noindent 

\begin{align*}
         T^{-1}_V \big( \psi_{k_R} \big)  \cdot J \equiv J \cdot T^{-1}_V  \big( \psi_{k_R}  \big)     \text{ } \text{ , } 
\end{align*}

\noindent by the commutation rules for $J$. 

\bigskip

\noindent Altogether, by inspection the basis elements of the map $T^{-1}_V$ analyzed above coincide with the coefficients of the complexification of the Ashkin-Teller propagation matrix away from $\frac{1}{4} \mathrm{log} \big( 3 \big)$. Hence $P^{\mathrm{AT}}$ and $\big( P^{\mathrm{AT}} \big)^{-1}$ are conjugate, from which we conclude the argument. \boxed{}

\bigskip

\noindent From the results of this subsection, we now transition back to the loop model.

\subsubsection{Loop model}

\noindent For the loop $\mathrm{O} \big( n \big)$ model, the desired basis for the transfer matrix, from a construction due to Blöte, and Nienhuis, [3], takes the form,

\begin{align*}
     Z \equiv \underset{\mathrm{configurations}}{\sum}    K^{\epsilon} n^L       \text{ } \text{ , } 
\end{align*}

\noindent corresponding to the product of the number of edges for each $K$, for real $K$, raised to some arbitrary, strictly positive $\epsilon$, and the number of loops, for some strictly positive, real $L$, following an expansion of the partition function as given above. For the number of occupied vertices in a loop configuration, the partition above is equivalent to,

\begin{align*}
     Z    \equiv  \underset{\mathrm{configurations}}{\sum}    K^{N_V} n^L        \text{ } \text{ , } 
\end{align*}

\noindent corresponding to the number of vertices, $N_V$, which are occupied, under the assumption,

\begin{align*}
 \epsilon  \approx N_V \text{ } \text{ , } 
\end{align*}

\noindent To consolidate vertices which could belong to the same loop of the configuration, if two lines, $k$ and $m$, both of which intersect with the finite volume, then we say that $k$ and $m$ are connected. The same observation can be applied to more than $k$ lines, in particular an arbitrary number of lines, to say that more than two lines of the loop configuration are connected. Iteratively, one can define the partition function of the finite volume over $\textbf{H}$, with side length $N+1$, in terms of the partition function of the finite volume over $\textbf{H}$, with side length $N$, with,

\begin{align*}
       Z^{N+1}_{\alpha} \equiv \underset{ \beta \in \textbf{N}}{\sum} T_{\alpha \beta} Z^N_{\beta}   \text{ } \text{ , } 
\end{align*}

\noindent for a connectivity parameter $\alpha$. For the parameters,

\begin{align*}
  N^{\prime}_v \equiv N_v + n_v  \text{ } \text{ , } \\ N^{\prime}_l \equiv N_l + n_l  \text{ } \text{ , } 
\end{align*}

\noindent from the expansion of the partition function over $\beta$,

\begin{align*}
   Z^{(N)} \equiv \underset{\beta \in \textbf{N}}{\sum}   Z^N_{\beta}   \text{ } \text{ , } 
\end{align*}

\noindent one can write,

\begin{align*}
     Z^{(N)}_{\alpha}  \equiv \underset{S ( N+1 )  \in G}{\sum}     \delta_{\alpha , \phi ( S ( N+1))}     K^{N^{\prime}_v} n^{N^{\prime}_l}      \text{ } \text{ , } 
\end{align*}

\noindent and also,

\begin{align*}
    Z^{(N)}_{\alpha}  \equiv       \underset{S ( N  )  \in G}{\sum}   K^{n_v} n^{N_l} \bigg[      \underset{g_{N} | S ( N ) }{\sum}    \delta_{\alpha , \phi ( S ( N))} K^{n_v} n^{n_l}  \bigg]            \text{ } \text{ , } 
\end{align*}

\noindent corresponding to the delta function $\delta$ for the parameters $\alpha$ and the connectivity $\phi \big( S \big( N +1 \big) \big)$ of line $N$, where the index of summation $g_N | S(N)$ denotes the sum of graphs over line $N$.

\bigskip

\noindent Hence, the partition function of the loop model can be expressed in terms of the transfer matrix, with,

\begin{align*}
   Z^{(N+1)}_{\alpha}  \equiv    \underset{ \beta \in \textbf{N}}{\sum}  \bigg[   \bigg[        \underset{g_{N+1} |\beta }{\sum}    \delta_{\alpha , \psi ( S ( N))} K^{n_v} n^{n_l}      \bigg]    \underset{S ( N  )  \in G}{\sum}   K^{n_v} n^{N_l} \bigg[     \underset{g_{N} | S ( N ) }{\sum}    \delta_{\beta , \phi ( S ( N))} K^{n_v} n^{n_l}  \bigg]    \bigg]   \text{ } \text{ , } 
\end{align*}

\noindent from the expansion,

\begin{align*}
 Z^{(N+1)}_{\alpha} \equiv \underset{ \beta \in \textbf{N}}{\sum}    \bigg[        \underset{g_{N+1} |\beta }{\sum}    \delta_{\alpha , \psi ( S ( N))} K^{n_v} n^{n_l}      \bigg] Z^{(N)}_{\beta}       \text{ } \text{ . } 
\end{align*}

\noindent \underline{Generator relations}. Equipped with the loop transfer matrix construction, we establish the following results, akin to \textbf{Proposition} \textit{3}, \textbf{Lemma} \textit{8}, and \textbf{Theorem} of \textit{2.2.1}. Again, we return to the fact of the HTE for the loop $\mathrm{O} \big( 1 \big)$ model, in which,

\begin{align*}
\textbf{P}^{\mathrm{loop}, \xi}_{\Lambda_{\textbf{H}}} \big[ \sigma \big] \overset{\mathrm{HTE}}{\sim }         \frac{n^{k ( \sigma )} x^{e(\sigma)} \mathrm{exp} \big(    h r \big( \sigma \big)  + h^{\prime}    r^{\prime} \big( \sigma \big)     \big) }{Z^{\mathrm{loop},\xi}_{\Lambda_{\textbf{H}}}\big( \sigma \big)  }          
\underset{n \equiv 1 , h^{\prime} \equiv 0}{\overset{\beta^{\mathrm{loop}} \equiv \frac{1}{2} | \mathrm{log}  x |}{\longleftrightarrow}} \textbf{P}^{\mathrm{Ising} , \chi}_{\Lambda_{\textbf{Z}^2}} \big[ \sigma^{\mathrm{Ising}} \big] \equiv \mathrm{O} \big( 1 \big) \text{ } \mathrm{measure}  \text{ } \text{ , } 
\end{align*}

\noindent introduce loop generators,

\begin{align*}
  \frac{p^{\mathrm{loop}}_u}{\beta} \equiv  h r \big( \sigma \big)   \text{ } \text{ , } \\ \frac{q^{\mathrm{loop}}_u}{\beta} \equiv   \mathrm{exp} \big( \mathrm{log} \big( x^{e(\sigma)} \big) \big)      \text{ } \text{ , } 
\end{align*}

\noindent because,

\begin{align*}
   x^{e(\sigma)}  \equiv \mathrm{exp} \big( \mathrm{log} \big( x^{e(\sigma)} \big) \big) \Rightarrow   \mathrm{exp} \big[ \mathrm{log} \big( x^{e(\sigma)}\big)  + h^{\prime} r^{\prime} \big( \sigma \big) \big]         \text{ } \text{ . } 
\end{align*}

\noindent With the Blöte-Nienhuis construction of the transfer matrix, state the following result.

\bigskip

\noindent \textbf{Lemma}  $8^{*}$ ($\textit{analogue of the Ashkin-Teller generator relation for the loop model}$). For $p^{\mathrm{loop}}_k$ and $q^{\mathrm{loop}}_k$, one has the relations, 

\begin{align*}
  V^{\mathrm{loop},h} \equiv \mathrm{exp} \bigg[ \beta^{\mathrm{loop}}  p^{\mathrm{loop}}_u + \beta^{\mathrm{loop}} q^{\mathrm{loop}}_u \bigg]  \equiv \mathrm{exp} \bigg[ \beta^{\mathrm{loop}}  \big( p^{\mathrm{loop}}_u +  q^{\mathrm{loop}}_u \big)  \bigg] \text{ } \text{ , } 
\end{align*}

\noindent and,

\begin{align*}
   V^{\mathrm{loop},V} \equiv \bigg[             \mathrm{exp} \big[ \big( \beta^{\mathrm{loop}} - \big( \beta^{\mathrm{loop}}\big)^{*} \big)  \big( p^{\mathrm{loop}}_k +  q^{\mathrm{loop}}_k \big)  \big]    \bigg] \mathrm{exp} \bigg[ \big( \beta^{\mathrm{loop}} \big)^{*}  p^{\mathrm{loop}}_{u+\frac{1}{2}} + \big( \beta^{\mathrm{loop}} \big)^{*}q^{\mathrm{loop}}_{u+\frac{1}{2}} \bigg]  \text{ } \text{ , } 
\end{align*}

\noindent where the dual loop temperature $\big( \beta^{\mathrm{loop}}\big)^{*}$ is obtained from the correspondence,

\begin{align*}
   \big( \beta^{\mathrm{loop}}\big)^{*} \longleftrightarrow \beta^{\mathrm{loop}}  \text{ } \text{ , } 
\end{align*}

\noindent with the relations,

\begin{align*}
  \mathrm{tanh} \big( \big( \beta^{\mathrm{loop}} \big)^{*} \big) = \mathrm{exp} \big( - 2 \beta^{\mathrm{loop}} \big) \text{ } \text{ , } \\ S^{\mathrm{loop}} \equiv \mathrm{sinh} \big( 2  \big( \beta^{\mathrm{loop}} \big)^{*} 
 \big) \text{ } \text{ . } 
\end{align*}

\bigskip

\noindent To prove the \textbf{Lemma} above, one needs to introduce several subspaces, similar to those defined in \textit{2.2.1}. We denote these subspaces with similar notation, namely $\mathcal{S}^{\mathrm{loop}}_1 \equiv \mathcal{S}_1 , \cdots , \mathcal{S}^{\mathrm{loop}}_4 \equiv \mathcal{S}_4$. The dual subspaces to each of the four loop subspaces are given by, respectively, $\mathcal{S}^{\prime}_1 , \cdots , \mathcal{S}^{\prime}_4$.

\bigskip

\noindent $\textit{Proof of Lemma}\text{ }  8^{*}$.  By direct computation, from the loop generators,

\begin{align*}
        i \bigg[      p^{\mathrm{loop}}_u \frac{e_{\tau(u)}}{i}  + q^{\mathrm{loop}}_u \frac{e_{\tau^{\prime}(u)} }{i}   \bigg]   \equiv p^{\mathrm{loop}}_{u+\frac{1}{2}} e_{\tau^{\prime}(u)} + q^{\mathrm{loop}}_{u+\frac{1}{2}}    e_{(\tau^{\prime}(u))^{\prime}}  \text{ } \text{ , } 
\end{align*}

\noindent corresponding to the generator relation for $V^{\mathrm{loop},h}$, and,

\begin{align*}
      i \bigg[     p^{\mathrm{loop}}_{u-\frac{1}{2}} \frac{e_{\tau^{\prime}(u)}}{i}  +  q^{\mathrm{loop}}_{u+\frac{1}{2}}  \frac{(e_{\tau^{\prime}(u)})^{\prime}}{i} \bigg] \equiv          p^{\mathrm{loop}}_u         e_{\tau(u)} +               q^{\mathrm{loop}}_{u-1}       e_{\tau^{\prime}(u)}   \text{ } \text{ , } 
\end{align*}

\noindent corresponding to the generator relation for $V^{\mathrm{loop},V}$. The bases for each loop subspace with its dual is obtained from the correspondence,

\begin{align*}
    \mathcal{S}_1 \underset{e_{\tau(u)} \longleftrightarrow (e_{\tau(u))^{\prime}}}{\overset{(*)}{\longleftrightarrow}} \mathcal{S}^{\prime}_1 \text{ } 
\end{align*}

\noindent Furthermore, along the lines of the argument of the previous section for \textbf{Proposition} \textit{3}, we work towards concluding the argument by observing,

\begin{align*}
  \mathrm{exp} \bigg[ \big( \beta^{\mathrm{loop}} \big)^{*}  p^{\mathrm{loop}}_{u+\frac{1}{2}} + q^{\mathrm{loop}}_{u+\frac{1}{2}}\big)  \bigg] e_{\tau} = \mathrm{cosh} \big( \big( \beta^{*} \big)^{\mathrm{loop}} \big) e_{\tau} +  \mathrm{cosh} \big( \big( \beta^{*} \big)^{\mathrm{loop}} \big) e_{(\tau)^{\prime}}   \text{ } \text{ , } 
\end{align*}

\noindent for the first identity for $V^{\mathrm{loop},h}$, and also that the prefactor for the second identity takes the form,

\begin{align*}
  \mathrm{exp} \bigg[ \big( \beta^{\mathrm{loop}} - \big( \beta^{\mathrm{loop}}\big)^{*} \big)  \big( p^{\mathrm{loop}}_k +  q^{\mathrm{loop}}_k \big)  \bigg]     \text{ } \text{ , } 
\end{align*}

\noindent corresponding to the second identity for $V^{\mathrm{loop},V}$ after performing similar computations as given in \textit{2.2.1}, from which we conclude the argument. \boxed{}

\bigskip

\noindent \underline{Induced rotations}. From the correspondence between $\beta^{\mathrm{loop}}$ and the dual inverse temperature $\big(\beta^{\mathrm{loop}}\big)^{*}$, below we state a conjugation result.

\bigskip

\noindent $\textbf{Lemma} \text{ } 7^{*}$ (\textit{loop transfer matrix conjugation, from Ashkin-Teller transfer matrix conjugation}). The action,

\begin{align*}
      \big( V^{\mathrm{loop},h}\big)^{-\frac{1}{2}}  \cdot p^{\mathrm{loop}}_u \cdot   \big( V^{\mathrm{loop},h} \big)^{\frac{1}{2}}    = c p^{\mathrm{loop}}_u - i s q^{\mathrm{loop}}_u \text{ } \text{ , } \\           \big( V^{\mathrm{loop},V}\big)^{-\frac{1}{2}}      \cdot q^{\mathrm{loop}}_u 
 \cdot    \big( V^{\mathrm{loop},V} \big)^{\frac{1}{2}}    = is p^{\mathrm{loop}}_u + c q^{\mathrm{loop}}_u   \text{ } \text{ , } 
\end{align*}

\noindent takes the form above, under the composition $\cdot$, for the Blöte-Nienhuis construction of the loop transfer matrix,

\begin{align*}
  T_{\alpha \beta} \equiv  \underset{g_{N+1} |\beta }{\sum}    \delta_{\alpha , \psi ( S ( N))} K^{n_v} n^{n_l}      \text{ } \text{ . } 
\end{align*}

\noindent \textit{Proof of Lemma $7^{*}$}. Perform similar computations to those used for arguing that \textbf{Lemma} \textit{7} holds, from \textit{2.2.1}. \boxed{}

\bigskip

\noindent $\textbf{Lemma} \text{ } 8^{*}$ (\textit{loop commutation rule}). For the same maps $R,J$ as defined in $\textbf{Lemma} \text{ } \textit{8}$, the same identities hold, for the induced rotations.

\bigskip

\noindent \textit{Proof of Lemma $8^{*}$}. Perform similar computations to those used for arguing that \textbf{Lemma} \textit{7} holds, from \textit{2.2.1}. \boxed{}

\bigskip

\noindent $\textbf{Theorem}^{*}$ (\textit{loop conjugation action of} $T_V$). For the loop propagator matrix, away from $x_c \big( n \big)$, the map $T_V$ admits the representation,

\begin{align*}
  T_V \equiv \rho \cdot \big( P^{\mathrm{loop}} \big)^{\textbf{C}} \cdot \rho^{-1}  \text{ } \text{ , } 
\end{align*}

\noindent for the maps, 

\begin{align*}
  \rho : \big( \textbf{C}^2 \big)^{\textbf{I}^{**}} \longrightarrow \mathcal{W}   \text{ } \text{ , } 
\end{align*}

\noindent and the complexification,

\begin{align*}
\big( P^{\mathrm{loop}} \big)^{\textbf{C}}  : \big( \textbf{R}^2 \big)^{\textbf{I}^{*}}  \longrightarrow \big( \textbf{C}^2 \big)^{\textbf{I}^{*}} 
 \text{ } \text{ . } 
\end{align*}

\noindent $\textit{Proof of}\text{ }  Theorem^{*}$. Perform similar computations to those used for arguing that \textbf{Theorem} holds, from \textit{2.2.1}. \boxed{}

\bigskip

\noindent In the next section, we discuss Fock representations.

\subsection{Fock representation}

\noindent For both models, we introduce a decomposition of the algebra into a Fock decomposition, from a decomposition of the generators of the algebra, with,

\begin{align*}
  \mathcal{W} \equiv \mathcal{W}_{\mathrm{cr}}          \oplus \mathcal{W}_{\mathrm{ann}}   \text{ } \text{ , } 
\end{align*}

\noindent and with,

\begin{align*}
  \mathcal{W}^{\prime} \equiv \mathcal{W}^{\prime}_{\mathrm{cr}}          \oplus \mathcal{W}^{\prime}_{\mathrm{ann}}   \text{ } \text{ , } 
\end{align*}

\noindent where in the direct sum decomposition above, the subspace, with its bilinear form $\big( \cdot , \cdot\big)$, has creation and annihilation operators. These operators satisfy the identities,

\begin{align*}
  \big( w_{\mathrm{cr}} , w^{\prime}_{\mathrm{cr}} \big) = 0   \text{ } \text{ , } \\  \big( w_{\mathrm{ann}} , w^{\prime}_{\mathrm{ann}} \big) = 0  \text{ } \text{ , } 
\end{align*}

\noindent for $w_{\mathrm{cr}} , w^{\prime}_{\mathrm{cr}} \in \mathcal{W}_{\mathrm{cr}}$, and $w_{\mathrm{ann}} , w^{\prime}_{\mathrm{ann}} \in \mathcal{W}_{\mathrm{ann}}$. Analogous properties hold for the subspaces $\mathcal{W}^{\prime}_{\mathrm{cr}}$, and $\mathcal{W}^{\prime}_{\mathrm{ann}}$. As a result of the bilinear form over $\mathcal{W}$ being nondegenerate, the subspaaces of $\mathcal{W}$ spanned by the creation and annihilation generators have dimension that is preciselty half of that of the dimension of $\mathcal{W}$. From the subspaces $\mathcal{W}_{\mathrm{cr}}$ and $\mathcal{W}_{\mathrm{ann}}$, one can introduce the decomposition of the exterior algebra of $W_{\mathrm{cr}}$, with,

\begin{align*}
 \bigwedge \mathcal{W}_{\mathrm{cr}} =  \underset{0 \leq n \leq |\textbf{I}^{*}|}{\bigoplus}  \big( \wedge^n \mathcal{W}_{\mathrm{cr}} \big)   \text{ } \text{ , } 
\end{align*}

\noindent into components for the $n$ th wedge product of $\mathcal{W}_{\mathrm{cr}}$ for the Ashkin-Teller model, or similarly, 

\begin{align*}
 \bigwedge \mathcal{W}^{\prime}_{\mathrm{cr}} =  \underset{0 \leq n \leq |\textbf{I}^{**}|}{\bigoplus}  \big( \wedge^n \mathcal{W}^{\prime}_{\mathrm{cr}} \big)   \text{ } \text{ , } 
\end{align*}

\noindent for the loop model. From these wedge product decompositions over a direct sum, for bases $\big\{ a^{\dagger}_{\alpha} \big\}_{1 \leq \alpha \leq | \textbf{I}^{*}|}$ and $\big\{ a^{\prime,\dagger}_{\alpha} \big\}_{1 \leq \alpha \leq | \textbf{I}^{**}|}$ of $\mathcal{W}_{\mathrm{cr}}$ and $\mathcal{W}^{\prime}_{\mathrm{cr}}$, the action on the Fock spaces, $ \bigwedge \mathcal{W}_{\mathrm{cr}}$ and  $\bigwedge \mathcal{W}^{\prime}_{\mathrm{cr}}$, are given by,

\begin{align*}
   a^{\dagger}_{\alpha} \cdot   \big(     \underset{1 \leq i \leq n}{\wedge}  \alpha^{\dagger}_{\beta_i}  \big) \equiv a^{\dagger}_{\alpha} \big(  \underset{1 \leq i \leq n}{\wedge}  \alpha^{\dagger}_{\beta_i} \big)   \text{ } \text{ , } \\  a^{\dagger}_{\alpha}\cdot   \big(   \underset{1 \leq i \leq n}{\wedge}  \alpha^{\dagger}_{\beta_i}     \big)   \equiv \underset{1 \leq j \leq n}{\sum}     \big( -1 \big)^{j-i} \delta_{\alpha, \beta_j}  \big( \underset{1 \leq i \leq n}{\wedge}  \alpha^{\dagger}_{\beta_i}   \big)     \equiv  \underset{1 \leq j \leq n}{\sum}             \big( -1 \big)^{j-i} \delta_{\alpha, \beta_j}  \big(  \alpha^{\dagger}_{\beta_1} \wedge \cdots \wedge \alpha^{\dagger}_{\beta_{j-1}} \wedge \alpha^{\dagger}_{\beta_j} \wedge \cdots \wedge \alpha^{\dagger}_{\beta_n} \big)               \text{ } \text{ , } 
\end{align*}

\noindent corresponding to identities for the product of the Ashkin-Teller basis elements $a^{\dagger}_{\alpha}$, while for the loop basis elements $a^{\prime,\dagger}_{\alpha}$, with the wedge product over $n$ other basis elements, and,

\begin{align*}
   a^{\dagger}_{\alpha} \cdot   \big(     \underset{1 \leq i \leq n}{\wedge}  \alpha^{\dagger}_{\beta^{\prime}_i}  \big) \equiv a^{\dagger}_{\alpha} \big(  \underset{1 \leq i \leq n}{\wedge}  \alpha^{\dagger}_{\beta^{\prime}_i} \big)   \text{ } \text{ , } \\  a^{\prime,\dagger}_{\alpha}\cdot   \big(   \underset{1 \leq i \leq n}{\wedge}  \alpha^{\dagger}_{\beta^{\prime}_i}     \big)   \equiv \underset{1 \leq j \leq n}{\sum}     \big( -1 \big)^{j-i} \delta_{\alpha, \beta^{\prime}_j}  \big( \underset{1 \leq i \leq n}{\wedge}  \alpha^{\dagger}_{\beta_i}   \big)     \equiv  \underset{1 \leq j \leq n}{\sum}             \big( -1 \big)^{j-i} \delta_{\alpha, \beta^{\prime}_j}  \big(  \alpha^{\dagger}_{\beta^{\prime}_1} \wedge \cdots \wedge \alpha^{\dagger}_{\beta^{\prime}_{j-1}} \wedge \alpha^{\dagger}_{\beta^{\prime}_j} \wedge \cdots \wedge \alpha^{\dagger}_{\beta^{\prime}_n} \big)               \text{ } \text{ , } 
\end{align*}

\noindent corresponding to identities for the product of the basis element with the complex conjugate of $n$ basis elements. Finally, define,

\begin{align*}
  \mathrm{Pf} \big[ A \big] \equiv        \frac{1}{2^m m!} \underset{\pi \in S_{2m}}{\sum}    \mathrm{sgn} \big( \pi \big)     \underset{1 \leq s \leq m}{\prod}    A_{\pi ( 2s-1 ) , \pi (  2s) }       \text{ } \text{ , } 
\end{align*}

\noindent corresponding to the Pfaffian of an antisymmetric matrix, for permutations $\pi$ belonging to the $2m$ letter symmetric group.

\subsubsection{Ashkin-Teller model}

\noindent To formalize the Fock representation, introduce the following.

\bigskip

\noindent \textbf{Lemma} \textit{10} (\textit{Ashkin-Teller polarization}, \textbf{Lemma} \textit{11}, [13]). For a polarization of $\mathcal{W}$, one has the following correspondence between irreducible representations of the algebra, and the Fock space representation of the creation subspace of $\mathcal{W}$, in which,

\begin{align*}
    \mathrm{Irr} \big( \mathcal{W} \big) \longleftrightarrow \bigwedge \mathcal{W}_{\mathrm{cr}} \text{ } \text{ . } 
\end{align*}

\noindent The correspondence above implies the existence of an embedding into the irreducible representation space into $\mathcal{W}$ for which,

\begin{align*}
      \underset{1 \leq j \leq n}{\bigwedge}   \alpha^{\dagger}_{\beta_j }         \mapsto \bigg[ \underset{1 \leq j \leq n}{\prod}  \alpha^{\dagger}_{\beta_j }    \bigg] v_{\mathrm{vac}}  \text{ } \text{ , } 
\end{align*}

\noindent where $v_{\mathrm{vac}}$ denotes an element of $\mathcal{W} \subsetneq \mathcal{V}$ for which $\mathcal{W}_{\mathrm{ann}} v_{\mathrm{vac}} \equiv 0$.

\bigskip

\noindent \textit{Proof of Lemma 10}. Begin by considering some representation of the algebra. Irrespective of irreducibility, fix some nonzero $v^{(0)} \in \mathcal{V}$. Introduce the following action on basis elements of the Fock representation, with,

\begin{align*}
    a_{\alpha}  v^{(\alpha)}    \neq 0    \text{ } \text{ , } 
\end{align*}

\noindent and, otherwise,

\begin{align*}
         v^{(\alpha)} \equiv 0     \text{ } \text{ , } 
\end{align*}

\noindent for $1 \leq \alpha \leq \big| \textbf{I}^{*}\big|$. From the construction of the vector above, the vacuum vector stated in the \textbf{Lemma} is nonzero, and equals $v^{(|\textbf{I}^{*}|}$, from which one has that $\mathcal{W}_{\mathrm{ann}} v_{\mathrm{vac}} \equiv 0$. One can extend this argument to other subrepresentations of the Fock representation $\bigwedge \mathcal{W}_{\mathrm{cr}}$, in which any subrepresentation of the Fock representation, satisfying,

\begin{align*}
   \bigwedge  \big( \mathcal{W}_{\mathrm{cr}} \big)^{\prime} \neq \emptyset \text{ } \text{ , } 
\end{align*}

\noindent for $\big( \mathcal{W}_{\mathrm{cr}}\big)^{\prime} \subsetneq \mathcal{W}_{\mathrm{cr}}$, contains a vacuum vector which we denote with $1 \in \bigwedge \mathcal{W}_{\mathrm{cr}}$. Hence, the existence of the $1$ vector in the Fock representation implies irreducibility of the representation. Besides demonstrating irreducibility of the Fock representation, the correspondence given in the statement of the \textbf{Lemma} above defines an intertwining operation which is an embedding by the  irreducibility of the Fock representation. Furthermore, the intertwining operation between the irreducible representations of the algebra and the Fock representation implies that the operation, which is an embedding, must be surjective, and hence an isomorphism, from which we conclude the argument. \boxed{}

\bigskip

\noindent For the next result, we make use of the Pfaffian.

\bigskip

\noindent \textbf{Lemma} \textit{11} (\textit{Ashkin-Teller Wick's Formula}, \textbf{Lemma} \textit{12}, [13]). For a polarization of $\mathcal{W}$ and Fock representation of the creation subspace $\mathcal{W}_{\mathrm{cr}}$, fix a vacuum vector $v_{\mathrm{vac}} \in \bigwedge \big( \mathcal{W}_{\mathrm{cr}} \big)^{\prime}$, and a dual vacuum vector $v^{*}_{\mathrm{vac}} \in \bigwedge \big( \mathcal{W}^{*}_{\mathrm{cr}} \big)^{\prime}$. Under the choice of these basis elements such that $< v_{\mathrm{vac}} , v^{*}_{\mathrm{vac}} > \equiv 1$, for elements $\phi_1 , \cdots , \phi_n \in \mathcal{W}$,

\begin{align*}
       < v^{*}_{\mathrm{vac}} \big(  \underset{1 \leq i \leq n}{\prod} \phi_i  \big) v_{\mathrm{vac}} > \text{ } \equiv  \underset{1 \leq i,j \leq n }{\bigcup}\mathrm{Pf} \bigg[   <  v^{*}_{\mathrm{vac}} \big(  \underset{k \in \{ i,j \} }{\prod}  \phi_k    \big) v_{\mathrm{vac}}    >          \bigg]   \text{ } \text{ . }
\end{align*}

\noindent \textit{Proof of Lemma 11}. To demonstrate that the desired expression for the inner product between the dual vacuum vector, product of $\phi_i$, and the vacuum vector can be expressed as a union over Pfaffians, express $\phi_i$ as,

\begin{align*}
  \phi_i \equiv   \alpha_i  w_{\mathrm{cr}} +       \beta_i  w_{\mathrm{ann}}  \text{ } \text{ , } 
\end{align*}

\noindent for real $\alpha_i, \beta_i$, with $w_{\mathrm{cr}} \in \mathcal{W}_{\mathrm{cr}}$ and $w_{\mathrm{ann}} \in \mathcal{W}_{\mathrm{ann}}$. One can similarly define linear combinations for the remaining fields $\phi_2 , \cdots , \phi_n$. From the linear combination above. in the creation and annihilation bases, write,

\begin{align*}
       < v^{*}_{\mathrm{vac}} \big(  \underset{1 \leq i \leq n}{\prod} \big(  \alpha_i  w_{\mathrm{cr}} +       \beta_i  w_{\mathrm{ann}}  \big)  \big) v_{\mathrm{vac}} >          \text{ } \text{ , } 
\end{align*}

\noindent which is equivalent to, after distribtuing terms to the product over $1 \leq i \leq n$,

\begin{align*}
       <   \underset{1 \leq i \leq n}{\prod} \big(   v^{*}_{\mathrm{vac}} \alpha_i  w_{\mathrm{cr}} v_{\mathrm{vac}}  +      v^{*}_{\mathrm{vac}}  \beta_i  w_{\mathrm{ann}}  v_{\mathrm{vac}} \big)  >     \text{ } \text{ , } 
\end{align*}

\noindent from which we perform the following anticommutation calculation, (AC), which yields,

\begin{align*}
     <   \underset{1 \leq i \leq n}{\prod} \big(  \big( \alpha_i  w_{\mathrm{cr}} v_{\mathrm{vac}} 
 \big)^{\dagger }    v^{*}_{\mathrm{vac}}    +   \big( \beta_i  w_{\mathrm{ann}}  v_{\mathrm{vac}}  \big)^{\dagger}   v^{*}_{\mathrm{vac}}   \big)  >      \text{ } \text{ , } 
\end{align*}

\noindent from which we obtain the desired expression, from the union of Pfaffians,

\begin{align*}
 \underset{1 \leq i , j \leq n }{\bigcup}  \mathrm{Pf} \bigg[   \text{ }     \underset{1 \leq i , j \leq n}{\prod}          \bigg(   <    v^{*}_{\mathrm{vac}} \alpha_i  w_{\mathrm{cr}} v_{\mathrm{vac}}  +      v^{*}_{\mathrm{vac}}  \beta_i  w_{\mathrm{ann}}  v_{\mathrm{vac}}   >         \bigg)     \bigg] \text{ } \text{ , } \end{align*}

 \noindent which we rewrite as,

 \begin{align*}
 \underset{1 \leq i , j \leq n }{\bigcup}  \mathrm{Pf} \bigg[   <   \underset{1 \leq i \leq n}{\prod} \big(   v^{*}_{\mathrm{vac}} \alpha_i  w_{\mathrm{cr}} v_{\mathrm{vac}}  +      v^{*}_{\mathrm{vac}}  \beta_i  w_{\mathrm{ann}}  v_{\mathrm{vac}} \big)  >      \bigg]  \\ \overset{(\mathrm{AC})}{\equiv}    \underset{1 \leq i , j \leq n }{\bigcup}           \mathrm{Pf} \bigg[   <   \underset{1 \leq i \leq n}{\prod} \big(  \big( \alpha_i  w_{\mathrm{cr}} v_{\mathrm{vac}} 
 \big)^{\dagger }    v^{*}_{\mathrm{vac}}    +   \big( \beta_i  w_{\mathrm{ann}}  v_{\mathrm{vac}}  \big)^{\dagger}   v^{*}_{\mathrm{vac}}   \big)  >   \bigg]       \\ \equiv  \underset{1 \leq i , j \leq n }{\bigcup}   \mathrm{Pf} \bigg[  <  \underset{1 \leq i \leq n}{\prod}  \bigg(   \big( \alpha_i  w_{\mathrm{cr}} v_{\mathrm{vac}} 
 \big)^{\dagger }    v^{*}_{\mathrm{vac}}   \bigg)  + \underset{1 \leq i \leq n}{\prod}  \bigg(   \big( \beta_i  w_{\mathrm{ann}}  v_{\mathrm{vac}}  \big)^{\dagger}   v^{*}_{\mathrm{vac}}   \bigg)  >  \bigg] \\ \equiv  \underset{1 \leq i , j \leq n}{\bigcup} \mathrm{Pf} \bigg[  <  \underset{1 \leq i \leq n}{\prod} \bigg(    \big( \alpha_i  w_{\mathrm{cr}} v_{\mathrm{vac}} 
 \big)^{\dagger }    v^{*}_{\mathrm{vac}}  \bigg)  > + < \underset{1 \leq i \leq n}{\prod}    \bigg( \big( \beta_i  w_{\mathrm{ann}}  v_{\mathrm{vac}}  \big)^{\dagger}   v^{*}_{\mathrm{vac}}    \bigg) >          \bigg] \\ \equiv \underset{1 \leq i , j \leq n}{\bigcup} \bigg[ \text{ } \mathrm{Pf} \bigg[          <   \underset{1 \leq i \leq n}{\prod}   \bigg(   \alpha^{\dagger}_i w^{\dagger}_{\mathrm{cr}} v^{\dagger}_{\mathrm{cr}} v^{*}_{\mathrm{vac}}   \bigg)           >            \bigg] + \mathrm{Pf} \bigg[      <  \underset{1 \leq i \leq n}{\prod}        \bigg(  \beta^{\dagger}_i          w^{\dagger}_{\mathrm{cr}} w^{\dagger}_{\mathrm{ann}} v^{*}_{\mathrm{vac}}         \bigg)     >                \bigg] \text{ }  \bigg] \end{align*}

 \noindent which, for the final steps, implies,
 
 \begin{align*} \underset{1 \leq i , j \leq n}{\bigcup} \mathrm{Pf} \bigg[        <   v^{*}_{\mathrm{vac}}      \bigg(    \big( \underset{1 \leq i \leq n}{\prod}  \alpha_i  \big)  w_{\mathrm{cr}} \bigg)  >        \bigg]  + \underset{1 \leq i , j \leq n}{\bigcup} \mathrm{Pf} \bigg[      <    \bigg(   \big( \underset{1 \leq i \leq n}{\prod} \beta_i   \big)    w_{\mathrm{ann}}  \bigg)  v_{\mathrm{vac}}  >          \bigg]  \\ \equiv \underset{1 \leq i , j \leq n}{\bigcup} \mathrm{Pf} \bigg[    <   v^{*}_{\mathrm{vac}}      \bigg(    \big( \underset{1 \leq i \leq n}{\prod}  \alpha_i  \big)  w_{\mathrm{cr}} \bigg) 
  > + < \bigg(   \big( \underset{1 \leq i \leq n}{\prod} \beta_i   \big)    w_{\mathrm{ann}}  \bigg)  v_{\mathrm{vac}}  >            \bigg]  \\ \equiv \underset{1 \leq i , j \leq n}{\bigcup} \mathrm{Pf} \bigg[    <   v^{*}_{\mathrm{vac}}      \bigg(    \big(  \alpha_1 \times \cdots \times \alpha_n   \big)  w_{\mathrm{cr}} \bigg) 
 +  \bigg(   \big(  \beta_1 \times \cdots \times \beta_n   \big)    w_{\mathrm{ann}}  \bigg)  v_{\mathrm{vac}}  >            \bigg]   \\ \equiv \underset{1 \leq i,j \leq n }{\bigcup}\mathrm{Pf} \bigg[   <  v^{*}_{\mathrm{vac}} \big(  \underset{k \in \{ i,j \} }{\prod}  \phi_k    \big) v_{\mathrm{vac}}    >          \bigg]   \text{ } \text{ , } 
\end{align*}

\noindent from which we conclude the argument. \boxed{}

\bigskip

\noindent We introduce the final item of the subsection below, for polarization at low temperatures.

\bigskip

\noindent \textbf{Lemma} \textit{12} (\textit{Ashkin-Teller polarization at low temperature}, \textbf{Lemma} \textit{13}, [13]). One has, for the Ashkin-Teller polarization, that,

\begin{align*}
  \mathcal{W}^{(+)}_{\mathrm{cr}} \equiv \mathrm{span} \big\{    p^{\mathrm{AT}}_k - i q^{\mathrm{AT}}_k \text{ } \big| \text{ } k \in \textbf{I}^{*}      \big\}   \text{ } \text{ , } \\ \mathcal{W}^{(+)}_{\mathrm{ann}} \equiv \mathrm{span} \big\{      p^{\mathrm{AT}}_k + i q^{\mathrm{AT}}_k \text{ } \big| \text{ } k \in \textbf{I}^{*}      \big\} 
 \text{ } \text{ . } 
\end{align*}

\bigskip

\noindent \textit{Proof of Lemma 12}. It is straightforward to check, by direct computation, that the creation and annihilation subspaces are spanned by the sets of vectors given above. \boxed{}

\bigskip

\noindent \textbf{Lemma} \textit{13} (\textit{Ashkin-Teller polarization in the vanishing temperature limit}, \textbf{Lemma} \textit{15}, [13]). For $N \in \textbf{N}$, a polarization is defined from the two subspaces,

\begin{align*}
          \mathcal{W}^{(+),N}_{\mathrm{cr}}  \equiv \mathrm{span}  \big\{  V^{-N}   \big(    p^{\mathrm{AT}}_k - i q^{\mathrm{AT}}_k       \big) V^N \text{ } 
 \big| \text{ }    k \in \textbf{I}^{*} 
 \big\}   \text{ } \text{ , } \\    \mathcal{W}^{(+),N}_{\mathrm{ann}}     \equiv \mathrm{span}  \big\{      p^{\mathrm{AT}}_k + i q^{\mathrm{AT}}_k      \text{ }  \big| \text{ }   k \in \textbf{I}^{*}   \big\}    \text{ } \text{ , } 
\end{align*}

\noindent for countably many inverse temperatures $\beta$, given a wedge product,

\begin{align*}
  \bigwedge \mathcal{W}^{(+),N}_{\mathrm{cr}}  \text{ } \text{ , } 
\end{align*}

\noindent and vacuum vectors,

\begin{align*}
     v^{(+),N}_{\mathrm{vac}} \equiv          e_{(+)}     \text{ }\text{ , }  \\ \big( v^{(+),N}_{\mathrm{vac}}\big)^{*}  \equiv       \big( e_{(+)} \big)^{*} \equiv    \frac{e^{\mathrm{T}}_{(+)} V^N }{e^{\mathrm{T}}_{(+)} V^N e_{(+)}}    \text{ } \text{ . } 
\end{align*}

\bigskip

\noindent $\textit{Proof of Lemma 13}$. It is straightforward to check, by direct computation, that the creation and annihilation subspaces are spanned by the sets of vectors given above, in which $ \mathcal{S}_{(+)} \overset{\sim}{=} \wedge \mathcal{W}_{\mathrm{cr}}^{(+),N}$, and $\mathcal{S}_{(-)} \overset{\sim}{=} \wedge \mathcal{W}_{\mathrm{ann}}^{(+),N}$ (refer to the arguments for \textbf{Lemma} \textit{15} in [13]). \boxed{}

\bigskip

\noindent In comparison to the polarization above which establishes an isomorphism between the direct sum of the annihilator and creation subspaces and the Fock representation, we introduce the physical polarization below which is determined by the eigenvectors of $T_V$.

\bigskip

\noindent \textbf{Lemma} \textit{14} (\textit{physical polarization for the Ashkin-Teller model}, \textbf{Lemma} \textit{16}, [13]). For the subspace $\mathcal{W}^{\mathrm{phys}}_{\mathrm{cr}} \subsetneq \mathcal{W}$ spanned by eigenvectors of $T_V$ with eigenvalues $<1$, and the subspace $\mathcal{W}^{\mathrm{phys}}_{\mathrm{ann}} \subsetneq \mathcal{W}$ spanned by eigenvectors $T_V$ with eigenvectors $>1$, we say that the direct sum $\mathcal{W}^{\mathrm{phys}}_{\mathrm{cr}} \oplus \mathcal{W}^{\mathrm{phys}}_{\mathrm{ann}}$ is a physical polarization of the Ashkin-Teller model. Hence, $\mathcal{S}_{+} \overset{\sim}{=} \wedge \mathcal{W}^{\mathrm{phys}}_{\mathrm{cr}}$.

\bigskip

\noindent \textit{Proof of Lemma 14}. To exhibit that the form of the physical polarization holds, recall from previous arguments that if $1$ is not an eigenvalue, then for vectors $u,v \in \mathcal{W}$, from the fact that $\big( T_V u , T_V v \big) \equiv \big( u , v \big)$, the tuple is nonzero iff $\big( u , v \big) \neq 0 \Longleftrightarrow u \equiv v^{-1} \Longleftrightarrow u^{-1} \equiv u$. Hence the binlinear form $\big( \cdot , \cdot \big)$ only vanishes when it is defined over $\mathcal{W}^{\mathrm{phys}}_{\mathrm{cr}}$, or over $\mathcal{W}^{\mathrm{phys}}_{\mathrm{ann}}$. As a result $T_V$ is diagonalizable, with eigenvalues that are real, and neither of which equal $1$, from which we conclude the argument. \boxed{}

\bigskip

\noindent Besides the polarization and its physical counterpart, we list the properties below.

\bigskip

\noindent \textbf{Proposition} \textit{AT 1} (\textit{properties of the Ashkin-Teller polarization}, \textbf{Proposition} \textit{17}, [13]). Given a physical polarization, introduce the basis,

\begin{align*}
   \mathrm{span} \big\{   \alpha_{\alpha}^{|\textbf{I}^{*}|} \text{ } \big| \text{ }   1 \leq \alpha \leq \big| \textbf{I}^{*} \big|     \big\}  \text{ } \text{ , } 
\end{align*}

\noindent for $\mathcal{W}^{\mathrm{phys}}_{\mathrm{ann}}$, which under the induced rotation $T_V$ yields a basis of the form,

\begin{align*}
  T_V \bigg[  \mathrm{span} \big\{   \alpha_{\alpha}^{|\textbf{I}^{*}|} \text{ } \big| \text{ }   1 \leq \alpha \leq \big| \textbf{I}^{*} \big|     \big\} \bigg] = \mathrm{span} \big\{ \lambda_{\alpha} a_{\alpha} \text{ } \big| \text{ }  1 \leq \alpha \leq \big| \textbf{I}^{*} \big|  \big\}        \text{ } \text{ , } 
\end{align*}

\noindent while the basis for the creation subspace takes the form,

\begin{align*}
  \bigg[  \mathrm{span} \big\{   \alpha_{\alpha}^{|\textbf{I}^{*}|} \text{ } \big| \text{ }   1 \leq \alpha \leq \big| \textbf{I}^{*} \big|     \big\} \bigg]^{\dagger} \equiv    \mathrm{span} \big\{   \big( \alpha_{\alpha}^{|\textbf{I}^{*}|}\big)^{\dagger} \text{ } \big| \text{ }   1 \leq \alpha \leq \big| \textbf{I}^{*} \big|     \big\}    \text{ } \text{ , } 
\end{align*}

\noindent which under the induced rotation $T_V$ yields a basis of the form,

\begin{align*}
  T_V \bigg[     \mathrm{span} \big\{   \big( \alpha_{\alpha}^{|\textbf{I}^{*}|}\big)^{\dagger} \text{ } \big| \text{ }   1 \leq \alpha \leq \big| \textbf{I}^{*} \big|     \big\}     \bigg]  =     \mathrm{span} \big\{    \lambda^{-1}_{\alpha} a^{\dagger}_{\alpha}    \text{ } \big| \text{ }      1 \leq \alpha \leq \big| \textbf{I}^{*} \big|       \big\}     \text{ } \text{ . }
\end{align*}

\noindent With the four bases above, it is possible, for some $v \in \mathcal{S}$, that:

\begin{itemize}
    \item[$\bullet$] \underline{Case one}: The product $a^{\dagger}_{\alpha} v \equiv 0 \in \mathcal{S}$.

    \item[$\bullet$] \underline{Case two}: The product $a^{\dagger}_{\alpha} v \neq 0 \in \mathcal{S}$. 
\end{itemize}

\noindent For the remaining possibility, we also have:

\begin{itemize}
    \item[$\bullet$] \underline{Case one}: The product $a_{\alpha} v \equiv 0 \in \mathcal{S}$.

    \item[$\bullet$] \underline{Case two}: The product $a_{\alpha} v \neq 0 \in \mathcal{S}$. 
\end{itemize}

\noindent From the first two cases above, it is possible that the eigenvectors living in $\mathcal{S}$ take the form:

\begin{itemize}
    \item[$\bullet$] \underline{Case one}: From the product $a^{\dagger}_{\alpha} v$ given in \underline{Case one} above, the eigenvalue $\lambda^{\dagger}_{\alpha} \equiv 0$.

    \item[$\bullet$] \underline{Case two}: From the product $a^{\dagger}_{\alpha} v$ given in \underline{Case two} above, the eigenvalue $\lambda^{\dagger}_{\alpha} \neq 0$.
\end{itemize}

\noindent From the second two cases above, it is possible that the eigenvectors living in $\mathcal{S}$ take the form:

\begin{itemize}
    \item[$\bullet$] \underline{Case one}: From the product $a_{\alpha} v$ given in \underline{Case one} above, the eigenvalue $\lambda_{\alpha} \equiv 0$.

    \item[$\bullet$] \underline{Case two}: From the product $a_{\alpha} v$ given in \underline{Case two} above, the eigenvalue $\lambda_{\alpha} \neq 0$.
\end{itemize}

\noindent From the four possible cases described above, if we denote $\Lambda_0$ as the largest eigenvalue in the spectrum of $\mathcal{V}$, with corresponding eigenvector $v^{\mathrm{phys}}_{\mathrm{vac}} \in \mathcal{S}_{+}$, then $\mathrm{span} \big\{          \prod_{1 \leq i \leq n} \alpha^{\dagger}_{\alpha_i} v_{\mathrm{vac}}  \text{ } \big| \text{ } v_{\mathrm{vac}} \in \mathcal{V}  \big\}$ forms a basis of $\mathcal{S}_{+}$, with spectrum,

\begin{align*}
   \mathrm{spec} \big( \mathcal{S}_{+} \big) =  \Lambda_0 \underset{1 \leq s \leq n}{\prod} \lambda^{-1}_{\alpha_s }\text{ } \text{ . } 
\end{align*}

\noindent \textit{Proof of Proposition AT 1}. For an eigenvector $v \in \mathcal{S}$, which satisfies the eigenvalue-eigenvector equation $V v = \Lambda v$, for the basis elements, $\big\{ a_{\alpha} \big\}$, $v a_{\alpha} v \equiv \big( V a_{\alpha} V^{-1} \big) V v \equiv T_V \big( a_{\alpha} \big) V v \equiv \lambda_{\alpha} \Lambda a_{\alpha} v$, and similarly, for the dual basis elements $\big\{ a^{\dagger}_{\alpha} \big\}$, one has, $v a^{\dagger}_{\alpha} v \equiv \big( V a^{\dagger}_{\alpha} V^{-1} \big) V v \equiv T_V \big( a^{\dagger}_{\alpha} \big) V v \equiv \lambda^{\dagger}_{\alpha} \Lambda a^{\dagger}_{\alpha} v$. To demonstrate that the remaining item holds, observe that the fact that $v^{\mathrm{phys}}_{\mathrm{vac}}\equiv 0$ when multiplied by any vector of the annihilator subspace $\mathcal{W}_{\mathrm{ann}}$, because,

\begin{align*}
     \lambda_{\alpha} \Lambda_0 > \underset{v}{\mathrm{sup}} \big\{ v : v \in \mathrm{spec} \big( V \big) \big\}     \text{ } \text{ , } 
\end{align*}

\noindent from which we conclude the argument. \boxed{}

 \bigskip

 \noindent We conclude the subsection with the following result. 

 \bigskip

 \noindent \textbf{Theorem} \textit{2} (\textit{isomorphism between the complexified propagation operator and a direct sum of wedge products}, \textbf{Theorem} \textit{18}, [13]). From the complexified Ashkin-Teller propagator, $\big( P^{\mathrm{AT}} \big)^{\textbf{C}}$, denote the subspace $\mathcal{W}_{\cdot}$ as the one spanned by eigenvectors of the complexified propagator with magnitude $<1$. From the exterior algebra,

 \begin{align*}
    \bigwedge \mathcal{W}_{\cdot} = \underset{0 \leq n \leq | \textbf{I}^{*}|}{\bigoplus}     \wedge^n    \mathcal{W}_{\cdot}      \text{ } \text{ , }
 \end{align*}

 \noindent introduce,

\begin{align*}
  \rho_{+} : \mathcal{S}_{+} \longrightarrow \bigwedge \mathcal{W}_{\cdot }  \text{ } \text{ , } 
\end{align*}

 \noindent so that,

\begin{align*}
  \rho_{+} \cdot V \cdot \rho_{+} = \mathcal{C} \big( \rho_{+} , V \big)  \Gamma \bigg[ \big( P^{\mathrm{AT}} \big)^{\textbf{C}} \bigg]  \text{ } \text{ , } 
\end{align*}

\noindent for some constant $\mathcal{C} \big( \rho_{+} , V \big) \equiv \mathcal{C}$, and,

\begin{align*}
  \Gamma \big( \big( P^{\mathrm{AT}} \big)^{\textbf{C}} \big) = \underset{0 \leq n \leq | \textbf{I}^{*}|} {\bigoplus}  \bigg[ \big( P^{\mathrm{AT}} \big)^{\textbf{C}}|_{\mathcal{W}_{\cdot}} 
 \bigg]^{\otimes n} \text{ } \text{ , } 
\end{align*}

\noindent for parameters in the complexified propagator which are below the $\frac{1}{4} \mathrm{log} \big( 3 \big)$ threshold.

\bigskip

\noindent \textit{Proof of Theorem 2}. To demonstrate that the result above holds, observe, from a previous result, that the isomorphism,

\begin{align*}
    \mathcal{S}_{+}    \overset{\sim}{=}     \bigwedge \bigg[ \mathcal{W}^{\mathrm{phys}}_{\mathrm{cr}} \oplus     \mathcal{W}^{\mathrm{phys}}_{\mathrm{ann}}\bigg]   \text{ } \text{ , } 
\end{align*}

\noindent holds. From the previous result, the transfer matrix of the Ashkin-Teller model, in the basis,

\begin{align*}
\mathrm{span} \big\{   \underset{1 \leq i \leq n}{\wedge} a^{\dagger}_{\alpha_i}  \text{ } \big| \text{ } 1 \leq a \leq \big| \textbf{I}^{*} \big|  \big\} \text{ } \text{ , } 
\end{align*}
 
\noindent has the same eigenvalues as the direct sum of transfer matrices raised to the $n$ th power of the tensor product,

\begin{align*}
  \underset{0 \leq n \leq | \textbf{I}^{*}|}{\bigoplus} \bigg[            T_V |_{W^{\mathrm{phys} }_{\mathrm{cr}} \otimes W^{\mathrm{phys} }_{\mathrm{ann}} }     \bigg]^{\otimes n} \text{ } \text{ . } 
\end{align*}

\noindent Hence, the induced rotation,

\begin{align*}
  T_V : \mathcal{W} \longrightarrow \mathcal{W}  \text{ } \text{ , } 
\end{align*}

\noindent satisfies,

\begin{align*}
  T_V \overset{\sim}{=}     \big( P^{\mathrm{loop}} \big)^{\textbf{C}}          \text{ } \text{ . } 
\end{align*}

\noindent Similarly, two other isomorphisms hold, in which,

\begin{align*}
    T|_{W^{\mathrm{phys}}_{\mathrm{cr}}} \overset{\sim}{=}     \big( P^{\mathrm{loop}} |_{(\mathcal{W} \cdot)_{\mathrm{cr}}} \big)^{\textbf{C}}        \text{ } \text{ , } \\  T|_{W^{\mathrm{phys}}_{\mathrm{ann}}} \overset{\sim}{=}     \big( P^{\mathrm{loop}} |_{(\mathcal{W} \cdot)_{\mathrm{ann}}} \big)^{\textbf{C}}      \text{ } \text{ . } 
\end{align*}

\noindent Hence,

\begin{align*}
  T|_{\mathcal{W}^{\mathrm{phys}}_{\mathrm{cr}}} \bigoplus     T|_{\mathcal{W}^{\mathrm{phys}}_{\mathrm{ann}}}  \overset{\sim}{=}   \big( P^{\mathrm{loop}} |_{(\mathcal{W} \cdot)_{\mathrm{cr}}} \big)^{\textbf{C}}  \bigoplus       \big( P^{\mathrm{loop}} |_{(\mathcal{W} \cdot)_{\mathrm{ann}}} \big)^{\textbf{C}} \overset{\sim}{=}  \underset{(\mathcal{W}^{\prime} \cdot)_{\cdot} \in \{ (\mathcal{W}_{\cdot})_{\mathrm{cr}} , (\mathcal{W}_{\cdot})_{\mathrm{ann}}  \}}{\bigoplus}   \big( P^{\mathrm{loop}} |_{(\mathcal{W}^{\prime} \cdot)_{\cdot}} \big)^{\textbf{C}} \overset{\sim}{=}   \big( P^{\mathrm{loop}} |_{(\mathcal{W} \cdot)} \big)^{\textbf{C}}  \text{ , }
\end{align*}

\noindent from which we conclude the argument. \boxed{}

\subsubsection{Loop model}

\noindent To formalize the Fock representation, introduce the following.

\bigskip

\noindent \textbf{Lemma} $10^{*}$ (\textit{loop polarization}). For a polarization of $\mathcal{W}^{\prime}$, one has the following correspondence between irreducible representations of the algebra, and the Fock space representation of the creation subspace of $\mathcal{W}^{\prime}$, in which,

\begin{align*}
    \mathrm{Irr} \big( \mathcal{W}^{\prime} \big) \longleftrightarrow \bigwedge \mathcal{W}^{\prime}_{\mathrm{cr}} \text{ } \text{ . } 
\end{align*}

\noindent The correspondence above implies the existence of an embedding into the irreducible representation space into $\mathcal{W}^{\prime}$ for which,

\begin{align*}
      \underset{1 \leq j \leq n}{\bigwedge}   \alpha^{\dagger}_{\beta^{\prime}_j }         \mapsto \bigg[ \underset{1 \leq j \leq n}{\prod}  \alpha^{\dagger}_{\beta^{\prime}_j }    \bigg] v_{\mathrm{vac}}  \text{ } \text{ , } 
\end{align*}

\noindent where $v_{\mathrm{vac}}$ denotes an element of $\mathcal{W}^{\prime} \subsetneq \mathcal{V}$ for which $\mathcal{W}^{\prime}_{\mathrm{ann}} v_{\mathrm{vac}} \equiv 0$.

\bigskip

\noindent \textit{Proof of Lemma} $10^{*}$. Begin by considering some representation of the algebra. Irrespective of irreducibility, fix some nonzero $v^{(0)} \in \mathcal{V}^{\prime}$. Introduce the following action on basis elements of the Fock representation, with,

\begin{align*}
    a_{\alpha}  v^{(\alpha)}    \neq 0    \text{ } \text{ , } 
\end{align*}

\noindent and, otherwise,

\begin{align*}
         v^{(\alpha)} \equiv 0     \text{ } \text{ , } 
\end{align*}

\noindent for $1 \leq \alpha \leq \big| \textbf{I}^{**} \big|$. With this recursive definition of the vector above, the same components provided in \textbf{Lemma} \textit{10} follow, from which we conclude the argument. \boxed{}

\bigskip

\noindent For the next result, we make use of the Pfaffian.

\bigskip

\noindent \textbf{Lemma} $11^{*}$ (\textit{Loop Wick's Formula}). For a polarization of $\mathcal{W}^{\prime}$ and Fock representation of the creation subspace $\mathcal{W}^{\prime}_{\mathrm{cr}}$, fix a vacuum vector $v_{\mathrm{vac}} \in \bigwedge \big( \mathcal{W}^{\prime}_{\mathrm{cr}} \big)^{\prime}$, and a dual vacuum vector $v^{*}_{\mathrm{vac}} \in \bigwedge \big( \mathcal{W}^{\prime,*}_{\mathrm{cr}} \big)^{\prime}$. Under the choice of these basis elements such that $< v_{\mathrm{vac}} , v^{*}_{\mathrm{vac}} > \equiv 1$, for elements $\phi_1 , \cdots , \phi_n \in \mathcal{W}^{\prime}$,

\begin{align*}
       < v^{*}_{\mathrm{vac}} \big(  \underset{1 \leq i \leq n}{\prod} \phi_i  \big) v_{\mathrm{vac}} > \text{ } 
     \equiv  \underset{1 \leq i,j \leq n }{\bigcup}\mathrm{Pf} \bigg[   <  v^{*}_{\mathrm{vac}} \big(  \underset{k \in \{ i,j \} }{\prod}  \phi_k    \big) v_{\mathrm{vac}}    >          \bigg]   \text{ } \text{ . }
\end{align*}

\noindent \textit{Proof of Lemma} $11^{*}$. Perform computations along the lines of those provided for \textbf{Lemma} \textit{11} in the previous subsection. \boxed{}

\bigskip

\noindent We introduce the final item of the subsection below, for polarization at low temperatures.

\bigskip

\noindent \textbf{Lemma} $12^{*} $(\textit{loop polarization at low temperature}). One has, for the loop polarization, that,

\begin{align*}
  \mathcal{W}^{\prime,(+)}_{\mathrm{cr}} \equiv \mathrm{span} \big\{    p^{\mathrm{loop}}_u - i q^{\mathrm{loop}}_u \text{ } \big| \text{ } u \in \textbf{I}^{**}      \big\}   \text{ } \text{ , } \\ \mathcal{W}^{\prime,(+)}_{\mathrm{ann}} \equiv \mathrm{span} \big\{      p^{\mathrm{loop}}_u + i q^{\mathrm{loop}}_u \text{ } \big| \text{ } u \in \textbf{I}^{**}      \big\} 
 \text{ } \text{ . } 
\end{align*}

\bigskip

\noindent \textit{Proof of Lemma} $12^{*}$. It is straightforward to check, by direct computation, that the creation and annihilation subspaces are spanned by the sets of vectors given above. \boxed{}

\bigskip

\noindent \textbf{Lemma} $13^{*}$ (\textit{loop polarization in the vanishing temperature limit}, \textbf{Lemma} \textit{15}, [13]). For $N \in \textbf{N}$, a polarization is defined from the two subspaces,

\begin{align*}
          \mathcal{W}^{(+),N}_{\mathrm{cr}}  \equiv \mathrm{span}  \big\{  V^{-N}   \big(    p^{\mathrm{loop}}_u - i q^{\mathrm{loop}}_u       \big) V^N \text{ } 
 \big| \text{ }    u \in \textbf{I}^{**} 
 \big\}   \text{ } \text{ , } \\    \mathcal{W}^{(+),N}_{\mathrm{ann}}     \equiv \mathrm{span}  \big\{      p^{\mathrm{loop}}_u + i q^{\mathrm{loop}}_u      \text{ }  \big| \text{ }  u \in \textbf{I}^{**}   \big\}    \text{ } \text{ , } 
\end{align*}

\noindent for countably many inverse temperatures $\beta^{\mathrm{loop}}$, given a wedge product,

\begin{align*}
  \bigwedge \mathcal{W}^{\prime,(+),N}_{\mathrm{cr}}  \text{ } \text{ , } 
\end{align*}

\noindent and vacuum vectors,

\begin{align*}
     v^{\prime,(+),N}_{\mathrm{vac}} \equiv             e^{\prime}_{(+)}  \text{ }\text{ , }  \\ \big( v^{\prime,(+),N}_{\mathrm{vac}}\big)^{*}  \equiv   \frac{e^{\prime,\mathrm{T}}_{(+)} V^N }{e^{\prime,\mathrm{T}}_{(+)} V^N e^{\prime}_{(+)}}      \text{ } \text{ . } 
\end{align*}
 
\bigskip

\noindent \textit{Proof of Lemma} $13^{*}$. It is straightforward to check, by direct computation, that the creation and annihilation subspaces are spanned by the sets of vectors given above, in which $ \mathcal{S}^{\prime}_{(+)} \overset{\sim}{=} \wedge \mathcal{W}_{\mathrm{cr}}^{\prime,(+),N}$, and $\mathcal{S}^{\prime}_{(-)} \overset{\sim}{=} \wedge \mathcal{W}_{\mathrm{ann}}^{\prime,(+),N}$ (refer to the arguments for \textbf{Lemma} \textit{15} in [13]). \boxed{}

\bigskip

\noindent In comparison to the polarization above which establishes an isomorphism between the direct sum of the annihilator and creation subspaces and the Fock representation, we introduce the physical polarization below which is determined by the eigenvectors of $T_V$.

\bigskip

\noindent \textbf{Lemma} $14^{*}$ (\textit{physical polarization for the loop model}). For the subspace $\mathcal{W}^{\prime,\mathrm{phys}}_{\mathrm{cr}} \subsetneq \mathcal{W}^{\prime}$ spanned by eigenvectors of $T_V$ with eigenvalues $<1$, and the subspace $\mathcal{W}^{\prime,\mathrm{phys}}_{\mathrm{ann}} \subsetneq \mathcal{W}$ spanned by eigenvectors $T_V$ with eigenvectors $>1$, we say that the direct sum $\mathcal{W}^{\prime,\mathrm{phys}}_{\mathrm{cr}} \oplus \mathcal{W}^{\prime,\mathrm{phys}}_{\mathrm{ann}}$ is a physical polarization of the Ashkin-Teller model. Hence, $\mathcal{S}^{\prime}_{+} \overset{\sim}{=} \wedge \mathcal{W}^{\prime,\mathrm{phys}}_{\mathrm{cr}}$.

\bigskip

\noindent \textit{Proof of Lemma} $14^{*}$. To exhibit that the form of the physical polarization holds, appeal to very similar arguments to those provided for \textbf{Lemma} \textit{14} of the previous subsection. \boxed{}  

\bigskip

\noindent Besides the polarization and its physical counterpart, we list the properties below.

\bigskip

\noindent \textbf{Proposition} \textit{Loop 1} (\textit{properties of the loop polarization}). Given a physical polarization, introduce the basis,

\begin{align*}
   \mathrm{span} \big\{   \alpha_{\alpha}^{|\textbf{I}^{**}|} \text{ } \big| \text{ }   1 \leq \alpha \leq \big| \textbf{I}^{**} \big|     \big\}  \text{ } \text{ , } 
\end{align*}

\noindent for $\mathcal{W}^{\mathrm{phys}}_{\mathrm{ann}}$, which under the induced rotation $T_V$ yields a basis of the form,

\begin{align*}
  T_V \bigg[  \mathrm{span} \big\{   \alpha_{\alpha}^{|\textbf{I}^{**}|} \text{ } \big| \text{ }   1 \leq \alpha \leq \big| \textbf{I}^{**} \big|     \big\} \bigg] = \mathrm{span} \big\{ \lambda_{\alpha} a_{\alpha} \text{ } \big| \text{ }  1 \leq \alpha \leq \big| \textbf{I}^{**} \big|  \big\}        \text{ } \text{ , } 
\end{align*}

\noindent while the basis for the creation subspace takes the form,

\begin{align*}
  \bigg[  \mathrm{span} \big\{   \alpha_{\alpha}^{|\textbf{I}^{**}|} \text{ } \big| \text{ }   1 \leq \alpha \leq \big| \textbf{I}^{**} \big|     \big\} \bigg]^{\dagger} \equiv    \mathrm{span} \big\{   \big( \alpha_{\alpha}^{|\textbf{I}^{**}|}\big)^{\dagger} \text{ } \big| \text{ }   1 \leq \alpha \leq \big| \textbf{I}^{**} \big|     \big\}    \text{ } \text{ , } 
\end{align*}

\noindent which under the induced rotation $T_V$ yields a basis of the form,

\begin{align*}
  T_V \bigg[     \mathrm{span} \big\{   \big( \alpha_{\alpha}^{|\textbf{I}^{**}|}\big)^{\dagger} \text{ } \big| \text{ }   1 \leq \alpha \leq \big| \textbf{I}^{**} \big|     \big\}     \bigg]  =     \mathrm{span} \big\{    \lambda^{-1}_{\alpha} a^{\dagger}_{\alpha}    \text{ } \big| \text{ }      1 \leq \alpha \leq \big| \textbf{I}^{**} \big|       \big\}     \text{ } \text{ . }
\end{align*}

\noindent With the four bases above, it is possible, for some $v \in \mathcal{S}^{\prime}$, that:

\begin{itemize}
    \item[$\bullet$] \underline{Case one}: The product $a^{\prime,\dagger}_{\alpha} v \equiv 0 \in \mathcal{S}^{\prime}$.

    \item[$\bullet$] \underline{Case two}: The product $a^{\prime,\dagger}_{\alpha} v \neq 0 \in \mathcal{S}^{\prime}$. 
\end{itemize}

\noindent For the remaining possibility, we also have:

\begin{itemize}
    \item[$\bullet$] \underline{Case one}: The product $a^{\prime}_{\alpha} v \equiv 0 \in \mathcal{S}^{\prime}$.

    \item[$\bullet$] \underline{Case two}: The product $a^{\prime}_{\alpha} v \neq 0 \in \mathcal{S}^{\prime}$. 
\end{itemize}

\noindent From the first two cases above, it is possible that the eigenvectors living in $\mathcal{S}^{\prime}$ take the form:

\begin{itemize}
    \item[$\bullet$] \underline{Case one}: From the product $a^{\prime,\dagger}_{\alpha} v$ given in \underline{Case one} above, the eigenvalue $\lambda^{\prime,\dagger}_{\alpha} \equiv 0$.

    \item[$\bullet$] \underline{Case two}: From the product $a^{\prime,\dagger}_{\alpha} v$ given in \underline{Case two} above, the eigenvalue $\lambda^{\prime,\dagger}_{\alpha} \neq 0$.
\end{itemize}

\noindent From the second two cases above, it is possible that the eigenvectors living in $\mathcal{S}$ take the form:

\begin{itemize}
    \item[$\bullet$] \underline{Case one}: From the product $a^{\prime}_{\alpha} v$ given in \underline{Case one} above, the eigenvalue $\lambda^{\prime}_{\alpha} \equiv 0$.

    \item[$\bullet$] \underline{Case two}: From the product $a^{\prime}_{\alpha} v$ given in \underline{Case two} above, the eigenvalue $\lambda^{\prime}_{\alpha} \neq 0$.
\end{itemize}

\noindent From the four possible cases described above, if we denote $\Lambda_0$ as the largest eigenvalue in the spectrum of $\mathcal{V}^{\prime}$, with corresponding eigenvector $v^{\mathrm{phys}}_{\mathrm{vac}} \in \mathcal{S}^{\prime}_{+}$, then $\mathrm{span} \big\{          \prod_{1 \leq i \leq n} \alpha^{\dagger}_{\alpha_i} v^{\prime}_{\mathrm{vac}}  \text{ } \big| \text{ } v_{\mathrm{vac}} \in \mathcal{V}^{\prime}  \big\}$ forms a basis of $\mathcal{S}^{\prime}_{+}$, with spectrum,

\begin{align*}
   \mathrm{spec} \big( \mathcal{S}^{\prime}_{+} \big) =  \Lambda_0 \underset{1 \leq s \leq n}{\prod} \lambda^{\prime,-1}_{\alpha_s }\text{ } \text{ . } 
\end{align*}

\noindent \textit{Proof of Proposition Loop 1}. Identical results stated above for the loop model follow from arguments given in the previous subsection for the Ashkin-Teller model, from \textbf{Proposition} \textit{AT 1}. \boxed{}

 \bigskip

 \noindent We conclude the subsection with the following result.

 \bigskip

\noindent \textbf{Theorem} $2^{*}$ (\textit{isomorphism between the complexified propagation operator and a direct sum of wedge products}). From the complexified loop propagator, $\big( P^{\mathrm{loop}} \big)^{\textbf{C}}$, denote the subspace $\mathcal{W}_{\cdot}$ as the one spanned by eigenvectors of the complexified propagator with magnitude $<1$. From the exterior algebra,

 \begin{align*}
    \bigwedge \mathcal{W}^{\prime}_{\cdot} = \underset{0 \leq n \leq | \textbf{I}^{*}|}{\bigoplus}     \wedge^n    \mathcal{W}_{\cdot}      \text{ } \text{ , }
 \end{align*}

 \noindent introduce,

\begin{align*}
  \rho^{\prime}_{+} : \mathcal{S}_{+} \longrightarrow \bigwedge \mathcal{W}_{\cdot }  \text{ } \text{ , } 
\end{align*}

 \noindent so that,

\begin{align*}
  \rho^{\prime}_{+} \cdot V \cdot \rho^{\prime}_{+} = \mathcal{C} \big( \rho^{\prime}_{+} , V \big)  \Gamma \bigg[ \big( P^{\mathrm{loop}} \big)^{\textbf{C}} \bigg]  \text{ } \text{ , } 
\end{align*}

\noindent for some constant $\mathcal{C} \big( \rho^{\prime}_{+} , V \big) \equiv \mathcal{C}$, and,

\begin{align*}
  \Gamma \big( \big( P^{\mathrm{loop}} \big)^{\textbf{C}} \big) = \underset{0 \leq n \leq | \textbf{I}^{**}|} {\bigoplus}  \bigg[ \big( P^{\mathrm{loop}} \big)^{\textbf{C}}|_{\mathcal{W}^{\prime}_{\cdot}} 
 \bigg]^{\otimes n} \text{ } \text{ , } 
\end{align*}

\noindent for parameters in the complexified propagator below Nienhuis' critical point.

\bigskip

\noindent \textit{Proof of Theorem} $2^{*}$. Appeal to similar arguments as those for \textbf{Theorem} \textit{2} in the previous subsection. \boxed{}

\subsection{Operator correlations and observables}

\noindent We analyze correlations and observables for each model in the subsections below.

\subsubsection{Quantum correspondence with the Ashkin-Teller model}

\noindent We define the following. Over a finite volume of the square lattice which can be expressed with a Cartesian product of rows and columns, $\mathscr{I} \times \mathscr{J}$, introduce, over the state space $\big\{ \pm 1 \big\}^{\mathscr{I}\times \mathscr{J}}$,

\begin{align*}
           \mathscr{P}^{\xi}_{\Lambda} \big( \sigma \big) \equiv \frac{\mathrm{exp}\big[ \mathscr{H}^{\mathrm{AT}} \big( \tau , \tau^{\prime} \big)  \big]}{\mathscr{Z}^{\xi}_{\Lambda} \big( \sigma \big) }  \equiv    \frac{\mathrm{exp}\bigg[ \underset{i \sim j}{\sum} \big[  J \big(  \tau \big( i \big) \tau \big( j \big) + \tau^{\prime} \big( i \big) \tau^{\prime} \big( j \big) \big) + U \big(  \tau \big( i \big) \tau \big( j \big) \tau^{\prime} \big( i \big) \tau^{\prime} \big( j \big)  \big)     \big]   \bigg]}{\mathscr{Z}^{\xi}_{\Lambda} \big( \sigma \big) }                     \text{ } \text{ , } 
\end{align*}

\noindent corresponding to the Ashkin-Teller measure, where the partition function can be placed into correspondence with the quantum state,

\begin{align*}
            \mathscr{Z}^{\xi}_{\Lambda} \big( \sigma \big) \equiv \underset{\tau , \tau^{\prime} \in \{ \pm \}^{\mathscr{I} \times \mathscr{J}}}{\sum}   \mathrm{exp} \big[ \mathscr{H}^{\mathrm{AT}}  \big( \tau , \tau^{\prime} \big) \big]  \longleftrightarrow      \bra{f} V^{\mathrm{AT},N}  \ket{i}     \text{ } \text{ , } 
\end{align*}

\noindent where the states $\bra{f}$ and $\ket{i}$ respectively denote the spins of the Ashkin-Teller configuration at the leftmost, and rightmost, boundaries, which satisfy,

\begin{align*}
    f \equiv i \equiv \big( V^{\mathrm{AT}, h} \big)^{\frac{1}{2}} e_{(+)} \equiv \mathrm{exp} \bigg[ J \big( \frac{|\textbf{I}^{*}|}{2}+ \frac{|\textbf{I}^{*}|}{2} \big) + U \big( \frac{|\textbf{I}^{*}|}{2} \big)^2  \bigg] \equiv \mathrm{exp} \bigg[     \frac{|\textbf{I}^{*} |}{2}    \big(  2 J  +          \frac{U}{2}           \big)      \bigg]      \text{ } \text{ . } 
\end{align*}

\noindent Under the stipulation that boundary spins are $\equiv +1$, this means that for boundary conditions $\xi$, we take,

\begin{align*}
  \xi \equiv \big\{          \sigma^{\mathrm{AT}} \equiv \big( \tau , \tau^{\prime} \big) \in \big\{ \pm 1 \big\}^{\mathscr{I} \times \mathscr{J}}  \times \big\{ \pm 1 \big\}^{\mathscr{I} \times \mathscr{J}}  \text{ } \big| \text{ }      \bra{f } \equiv + 1 , \ket{i} \equiv + 1    \big\}   \text{ } \text{ . } 
\end{align*}

\bigskip

\noindent Furthermore, there exists a spin operator, $\hat{\tau}_j : \big\{ \pm 1 \big\}^{\mathscr{I}\times \mathscr{J}} \longrightarrow \big\{ \pm 1 \big\}^{\mathscr{I}\times \mathscr{J}}$ with the action,

\begin{align*}
  \hat{\tau_j} \big( \tau \big( i \big) \big) \mapsto -\tau \big( i \big)   \text{ } \text{ , } \\  \hat{\tau_j} \big( \tau^{\prime} \big( i \big) \big) \mapsto - \tau^{\prime} \big( i \big)  \text{ } \text{ , } 
\end{align*}

\noindent in which the spins of either Potts model at site $i$ of the lattice are reversed. With respect to the measure $\mathscr{P}^{\xi}_{\Lambda} \big( \cdot \big)$, the expectatioon for the Ashkin-Teller model can be written as,

\begin{align*}
  \mathscr{E}^{\xi}_{\Lambda} \big( \sigma_z \big) =    \underset{\tau_z , \tau^{\prime}_z \in \{ \pm 1 \}^{\mathscr{I} \times \mathscr{J}}}{\sum}   \sigma_z      \text{ }  \mathrm{d}  \mathscr{P}^{\xi}_{\Lambda} \big( \tau_z , \tau^{\prime}_z \big) \equiv    \frac{\underset{\tau , \tau^{\prime} \in \{ \pm 1 \}^{\mathscr{I} \times \mathscr{J}}}{\sum}  \sigma_z   \mathrm{exp} \big[ \mathscr{H}^{\mathrm{AT}}  \big( \tau_z , \tau^{\prime}_z \big) \big] }{\mathscr{Z}^{\xi}_{\Lambda} \big( \sigma_z \big) }  \equiv \frac{\underset{\tau , \tau^{\prime} \in \{ \pm 1 \}^{\mathscr{I} \times \mathscr{J}}}{\sum}  \sigma_z   \mathrm{exp} \big[ \mathscr{H}^{\mathrm{AT}}  \big( \tau_z , \tau^{\prime}_z \big) \big] }{ \underset{\tau_z,\tau^{\prime}_z \in \{ \pm 1 \}^{\mathscr{I} \times \mathscr{J}}}{\sum} \mathrm{exp} \big[ \mathscr{H}^{\mathrm{AT}}  \big( \tau_z , \tau^{\prime}_z \big) \big]}   \text{ , } 
\end{align*}

\noindent for $z = x + i y$. The final expression above is equivalent to the state,

\begin{align*}
  \frac{\bra{f} V^{\mathrm{AT}, N -y }   \big(    \hat{\tau \big( i \big)} \hat{\tau \big( j \big) }       \big)     V^{\mathrm{AT}, y }          \ket{i }}{\bra{f } V^{\mathrm{AT},N} \ket{i}}  \text{ } \text{ . } 
\end{align*}

\noindent Generalizing the expectation above to a product over spins instead of a single spin at $x$ yields the expression,

\begin{align*}
  \mathscr{E}^{\xi}_{\Lambda} \bigg[ \text{ } \underset{x \in \textbf{C}}{\prod} \sigma_x \text{ } \bigg] =  \underset{x \in \textbf{C}}{\prod} 
  \bigg[ \text{ } \frac{\underset{\tau , \tau^{\prime} \in \{ \pm 1 \}^{\mathscr{I} \times \mathscr{J}}}{\sum}  \sigma_x   \mathrm{exp} \big[ \mathscr{H}^{\mathrm{AT}}  \big( \tau_z , \tau^{\prime}_z \big) \big] }{ \underset{\tau_z,\tau^{\prime}_z \in \{ \pm 1 \}^{\mathscr{I} \times \mathscr{J}}}{\sum} \mathrm{exp} \big[ \mathscr{H}^{\mathrm{AT}}  \big( \tau , \tau^{\prime} \big) \big]} \bigg] \equiv \frac{    \underset{x \in \textbf{C}}{\prod}    \text{ }    \underset{\tau_z , \tau^{\prime}_z \in \{ \pm 1 \}^{\mathscr{I} \times \mathscr{J}}}{\sum}  \sigma_x   \mathrm{exp} \big[ \mathscr{H}^{\mathrm{AT}}  \big( \tau , \tau^{\prime} \big) \big] }{ \underset{\tau,\tau^{\prime} \in \{ \pm 1 \}^{\mathscr{I} \times \mathscr{J}}}{\sum} \mathrm{exp} \big[ \mathscr{H}^{\mathrm{AT}}  \big( \tau_z , \tau^{\prime}_z \big) \big]}   \text{ } \text{ , }
\end{align*}

\noindent which is equivalent to,

\begin{align*}
  \frac{\bra{e_{(+)}} V^{\mathrm{AT},N}  \underset{x \in \textbf{C}}{\underset{1 \leq i \leq r}{\prod}} \big( \hat{\sigma_i} \big)_x \ket{e_{(+)}}}{\bra{e_{(+)}} V^{\mathrm{AT},N} \ket{e_{(+)}}}  \text{ } \text{ . } 
\end{align*}

\noindent Moreover, from the expressions above, observe,

\begin{align*}
 \hat{\tau \big( x + iy  \big) \tau \big( x+i y^{\prime} \big) } \equiv \hat{ \tau \big( x + iy \big)} \hat{\tau \big( x+i y^{\prime} \big)}  = \bigg[ V^{\mathrm{AT},-y} \big(  \hat{\tau_x}    \big) V^{\mathrm{AT},y}  \bigg] \bigg[ V^{\mathrm{AT},-y^{\prime}}              \big( \hat{\tau_x} \big) V^{\mathrm{AT},y^{\prime}}  \bigg] \\ \equiv  V^{\mathrm{AT},-y} \big(  \hat{\tau_x}    \big) \bigg[ \text{ } \underset{i \in \{ y , -y^{\prime}\}}{\prod} V^{\mathrm{AT},i} \text{ } \bigg] \big( \hat{\tau_x} \big) V^{\mathrm{AT},y^{\prime}} 
 \text{ } \text{ . } 
\end{align*}

\noindent From the generators $p^{\mathrm{AT}}_k$ and $q^{\mathrm{AT}}_k$ analyzed in \textit{2.2.1}, over the complex plane,

\begin{align*}
       \psi \big( k + iy \big) =  V^{\mathrm{AT},-y} \psi_k V^{\mathrm{AT},y} \text{ } \text{ , } \\    \bar{\psi} \big( k + iy \big) =    V^{\mathrm{AT},-y} \bar{\psi_k }V^{\mathrm{AT},y}  \text{ } \text{ , } 
\end{align*}

\noindent for $k \in \textbf{I}^{*}$. Recall,

\begin{align*}
       \psi \big( k + i y \big) = \frac{i}{\sqrt{2}} \big( p^{\mathrm{AT}}_k + q^{\mathrm{AT}}_k \big)    \text{ } \text{ , } \\     \bar{\psi} \big( k + i y \big) = \frac{1}{\sqrt{2}} \big( p^{\mathrm{AT}}_k - q^{\mathrm{AT}}_k \big)    \text{ } \text{ . } 
\end{align*}

\noindent Over $\mathscr{I} \times \mathscr{J}$, for,

\begin{align*}
  \psi^{\mathrm{AT}} \equiv \psi   \text{ } \text{ , } 
\end{align*}

\noindent the multipoint correlation function for the fermion operator takes the form, for $\big\{ z_i \big\}_{i \in \textbf{C}}$,

\begin{align*}
   <                         \psi^{(1)} \big( z_1 \big) \times  \cdots  \times \psi^{(n)} \big( z_n \big)          >^{\xi}_{\mathscr{I}\times \mathscr{J}}  =      \frac{\bra{e_{(+)}} V^{\mathrm{AT},N}  \psi^{(1)} \big( z_1 \big) \times  \cdots  \times \psi^{(n)} \big( z_n \big)   \ket{e_{(+)}}}{\bra{e_{(+)}} V^{\mathrm{AT},N} \ket{e_{(+)}}}              \text{ } \text{ . } 
\end{align*}

\noindent Under the induced rotation introduced in \textit{2.2.1}, the basis functions $\phi \big( z \big)$ and $\bar{\phi \big( z \big)}$ have the image,

\begin{align*}
     R \big(  \psi \big( z \big)  \big) =  \bar{\psi \big( r \big( z \big) \big) }  \text{ } \text{ , } \\     R \big(\bar{\psi \big( z \big)}  \big)  =   \psi \big( r \big( z \big) \big)    \text{ } \text{ , } 
\end{align*}

\noindent where in each image under the induced rotation above,

\begin{align*}
      r \big(z \big) = r \big( x+ i y \big) = a + b - \big( x - iy \big)   \text{ } \text{ , } \\    \bar{r \big(z \big) }= \bar{r \big( x+ i y \big) } =    \bar{a+b - \big( x - iy \big) } = a+ b - \big( x + iy \big) 
\text{ } \text{ , } 
\end{align*}

\noindent correspond to two functions of complex arguments $x+iy$, and of $x-iy$.

\bigskip

\noindent With the Ashkin-Teller measure and the correspondence with quantum states described in this section for the fermion operator, the following notion of s-holomorphicity holds. To distinguish between previous notions of massive, and massless, s-holomorphicity which holds for the Ising, Ashkin-Teller and loop models, in the statement below we denote the Ashkin-Teller fermion operator with $\psi \big( z \big) \equiv \psi^{\mathrm{AT-F}}$, $\bar{\psi \big( z \big) } \equiv \bar{\psi^{\mathrm{AT-F}}}$. The complexification of the Ashkin-Teller fermion operator is obtained from entries of the tuple $\big( \psi^{\mathrm{AT-F}}, \bar{\psi^{\mathrm{AT-F}}}\big)$.

\bigskip

\noindent \textbf{Theorem} \textit{3} (\textit{massive s-holmorphicity for Ashkin-Teller fermions, from massive s-holomorphicity for Ising fermions}, \textbf{Theorem} \textit{19}, [13]). Fix $z \in \textbf{I}^{*} \times \textbf{J}$, and the same parameters $\nu$ and $\lambda$ provided in \textbf{Definition} \textit{1}, and \textbf{Definition} \textit{6}. For the Ashkin-Teller fermion operators $\psi^{\mathrm{AT-F}} \big( z \big) $ and $\bar{\psi^{\mathrm{AT-F}}} \big( z \big)$, there exists extensions of $\psi^{\mathrm{AT-F}}$, and of $\bar{\psi^{\mathrm{AT-F}}}$, to $\textbf{I}^{*} \times \textbf{J}$, such that,

\begin{align*}
   \psi^{\mathrm{AT-F}} \big( N \big) + \nu^{-1} \lambda \bar{\psi^{\mathrm{AT-F}} } \big( N \big) = \nu^{-1} \psi^{\mathrm{AT-F}} \big( E \big) + \lambda \bar{\psi^{\mathrm{AT-F}} \big( E \big) } \text{ } \text{ , } \\ \psi^{\mathrm{AT-F}} \big( N \big) + \nu \lambda^{-1} \bar{\psi^{\mathrm{AT-F}} \big( N \big)} = \nu \psi^{\mathrm{AT-F}} \big( W \big) + \lambda^{-1} \bar{\psi^{\mathrm{AT-F}} \big( W \big)} \text{ } \text{ , } \\ \psi^{\mathrm{AT-F}} \big( S \big) + \nu \lambda^3 \bar{\psi^{\mathrm{AT-F}} \big( S \big) } = \nu \psi^{\mathrm{AT-F}} \big( E \big) + \lambda^3 \bar{\psi^{\mathrm{AT-F}} \big( E \big)} \text{ } \text{ , } \\ \psi^{\mathrm{AT-F}} \big( S \big) + \nu^{-1} \lambda^{-3} \bar{\psi^{\mathrm{AT-F}} \big( S \big)} = \nu^{-1} \psi^{\mathrm{AT-F}} \big( W \big) + \lambda^{-3} \bar{\psi^{\mathrm{AT-F}} \big( W \big) } \text{ } \text{ , } 
\end{align*}

\noindent for any face of the complex plane with edges E,N,W,S. At the left, and right, boundary points of the finite volume, the fermion operators also satisfy,

\begin{align*}
  \psi^{\mathrm{AT-F}} \big( a + i y \big) + i \bar{\psi^{\mathrm{AT-F}} } \big( a + iy \big) = 0   \text{ } \text{ , } \\     \psi^{\mathrm{AT-F}} \big( b + i y \big) + i \bar{\psi^{\mathrm{AT-F}} } \big( b + iy \big) = 0            \text{ } \text{ , } 
\end{align*}

\noindent for all $y \in \textbf{J}^{*}$.

\bigskip

\noindent \textit{Proof of Theorem 3}. To exhibit that s-holomorphicity holds for the Ashkin-Teller fermion, we make use of similar notions as those which were used for demonstrating that the extension $h$ exists, and is unique. That is, if we fix some point in the vertical direction, $y \in \textbf{J}^{*}$, and another point in the interior of the horizontal direction, $z \in \big( \textbf{I}_y \backslash \partial \textbf{I}_y \big) \equiv I \big( \textbf{I}_y \big)$, then the equations for the Ashkin-Teller fermion about this point would take the form,

\begin{align*}
   \psi^{\mathrm{AT-F}} \big( z \big) + \nu^{-1} \lambda \bar{\psi^{\mathrm{AT-F}} } \big( z \big) = \nu^{-1} \psi^{\mathrm{AT-F}} \big( z^{\prime} \big) + \lambda \bar{\psi^{\mathrm{AT-F}} \big( z^{\prime} \big) } \text{ } \text{ , } \\ \psi^{\mathrm{AT-F}} \big( z \big) + \nu \lambda^{-1} \bar{\psi^{\mathrm{AT-F}} \big( z \big)} = \nu \psi^{\mathrm{AT-F}} \big( z^{\prime\prime} \big) + \lambda^{-1} \bar{\psi^{\mathrm{AT-F}} \big( z^{\prime\prime} \big)} \text{ } \text{ , } \\ \psi^{\mathrm{AT-F}} \big( z+1  \big) + \nu \lambda^3 \bar{\psi^{\mathrm{AT-F}} \big( z+1  \big) } = \nu \psi^{\mathrm{AT-F}} \big( z^{\prime} \big) + \lambda^3 \bar{\psi^{\mathrm{AT-F}} \big( z^{\prime} \big)} \text{ } \text{ , } \\ \psi^{\mathrm{AT-F}} \big( z+1  \big) + \nu^{-1} \lambda^{-3} \bar{\psi^{\mathrm{AT-F}} \big( z+1 \big)} = \nu^{-1} \psi^{\mathrm{AT-F}} \big( z^{\prime\prime} \big) + \lambda^{-3} \bar{\psi^{\mathrm{AT-F}} \big( z^{\prime\prime}  \big) } \text{ } \text{ , } 
\end{align*}

\noindent where the points in the above system of equations satisfy,

\begin{align*}
  z^{\prime} =  \big\{ z \in \textbf{C} \text{ } \big| \text{ } - 45 \text{ } \mathrm{degree} \text{ } \mathrm{rotation} \text{ } \mathrm{from} \text{ }  \mathrm{N}        \big\}   \text{ } \text{ , }  \\ z^{\prime\prime } =  \big\{ z \in \textbf{C} \text{ } \big| \text{ } 225 \text{ } \mathrm{degree} \text{ } \mathrm{rotation} \text{ } \mathrm{from} \text{ }  \mathrm{N}        \big\} \text{ } \text{ , }  \\        z+1 =    \big\{ z \in \textbf{C} \text{ } \big| \text{ } +1 \text{ } \mathrm{translation} \text{ } \mathrm{of} \text{ } z     \text{ } \mathrm{down}  \big\}           \text{ } \text{ , }
\end{align*}

\noindent Hence, the above relations from the system imply,

\begin{align*}
 \nu^{-1} \psi^{\mathrm{AT-F}} \big( z^{\prime} \big) + \lambda \bar{\psi^{\mathrm{AT-F}} \big( z^{\prime} \big) } - \nu^{-1}  \lambda \bar{\psi^{\mathrm{AT-F}} \big( z \big)} \equiv    \nu \psi^{\mathrm{AT-F}} \big( z^{\prime\prime} \big) + \lambda^{-1} \bar{\psi^{\mathrm{AT-F}} \big( z^{\prime\prime} \big)} -  \nu \lambda^{-1}      \bar{\psi^{\mathrm{AT-F}}\big( z \big) }        \\ \Updownarrow \\        \nu^{-1} \psi^{\mathrm{AT-F}} \big( z^{\prime} \big) + \lambda \bar{\psi^{\mathrm{AT-F}} \big( z^{\prime} \big) }    - \nu \psi^{\mathrm{AT-F} } \big( z^{\prime\prime} \big) - \lambda^{-1} \bar{\psi^{\mathrm{AT-F}} \big( z^{\prime\prime} \big) } =       \big(     \nu^{-1} \lambda - \nu \lambda^{-1}        \big)       \bar{\psi^{\mathrm{AT-F}}\big( z \big) }         \\ \Updownarrow \\        \frac{1}{\big(     \bar{\nu^{-1}} \bar{\lambda } - \bar{\nu} \bar{\lambda^{-1}}        \big)  }    \bigg[    \bar{\nu^{-1}} \bar{\psi^{\mathrm{AT-F}} \big( z^{\prime} \big)} + \bar{\lambda } \bar{\bar{\psi^{\mathrm{AT-F}} \big( z^{\prime} \big) }}    - \bar{\nu} \bar{ \psi^{\mathrm{AT-F} } \big( z^{\prime\prime} \big)} - \bar{\lambda^{-1}} \bar{ \bar{\psi^{\mathrm{AT-F}} \big( z^{\prime\prime} \big) } }   \bigg]       = \psi^{\mathrm{AT-F}} \big( z \big)      \text{ } \text{ , } 
\end{align*}

\noindent from the first and second equations of the system, and,

\begin{align*}
     \nu \psi^{\mathrm{AT-F}} \big( z^{\prime} \big) + \lambda^3 \bar{\psi^{\mathrm{AT-F}} \big( z^{\prime} \big)} -  \nu \lambda^3 \bar{\psi^{\mathrm{AT-F}} \big(  z+1\big) }\equiv \nu^{-1} \psi^{\mathrm{AT-F}} \big( z^{\prime\prime} \big) + \lambda^{-3} \bar{\psi^{\mathrm{AT-F}} \big( z^{\prime\prime}  \big) } - \nu^{-1} \lambda^{-3} \bar{\psi^{\mathrm{AT-F}} \big( z + 1 \big) }
 \\ \Updownarrow \\   \nu \psi^{\mathrm{AT-F}} \big( z^{\prime} \big) + \lambda^3 \bar{\psi^{\mathrm{AT-F}} \big( z^{\prime} \big)}       -   \nu^{-1} \psi^{\mathrm{AT-F}} \big( z^{\prime\prime} \big) - \lambda^{-3} \bar{\psi^{\mathrm{AT-F}} \big( z^{\prime\prime}  \big) }     =              \big( \nu \lambda^3 - \nu^{-1} \lambda^{-3} \big) \bar{\psi^{\mathrm{AT-F}} \big( z+1 \big) }  \\ \Updownarrow \\ \frac{1}{\big( \bar{\nu} \bar{\lambda^3} - \bar{\nu^{-1}} \bar{\lambda^{-3}}\big) }  \bigg[ \bar{\nu} \bar{\psi^{\mathrm{AT-F}} \big( z^{\prime} \big)} + \bar{\lambda^3 } \bar{\bar{\psi^{\mathrm{AT-F}} \big( z^{\prime} \big)}}       -  \bar{ \nu^{-1} } \bar{\psi^{\mathrm{AT-F}} \big( z^{\prime\prime} \big) } - \bar{\lambda^{-3} } \bar{\bar{\psi^{\mathrm{AT-F}} \big( z^{\prime\prime}  \big) }   }     \bigg] =  \psi^{\mathrm{AT-F}} \big( z + 1 \big)  \text{ } \text{ , } 
\end{align*}

\noindent from the third and fourth equations of the system. Further rearranging the second expression above yields an expression for $\bar{\psi^{\mathrm{AT-F}}\big( z \big) }$, which takes the form,

\begin{align*}
     \frac{1}{\big( \nu \lambda^3 - \nu^{-1} \lambda^{-3} \big) }  \bigg[  \nu \psi^{\mathrm{AT-F}} \big( z^{\prime} - 1 \big) + \lambda^3 \bar{\psi^{\mathrm{AT-F}} \big( z^{\prime} -1 \big)}       -   \nu^{-1} \psi^{\mathrm{AT-F}} \big( z^{\prime\prime} -1 \big) - \lambda^{-3} \bar{\psi^{\mathrm{AT-F}} \big( z^{\prime\prime} -1  \big) } \bigg] =    \bar{\psi^{\mathrm{AT-F}}\big( z \big) }\text{ } \text{ . } 
\end{align*}

\noindent Furthermore, each of the the expressions above can be extended to points lying along the horizontal line $\textbf{I}^{*}_{\frac{1}{2}}$, or along the vertical line $\textbf{J}^{*}_{\frac{1}{2}}$, in which,

\begin{align*}
     \frac{1}{\big( \nu \lambda^3 - \nu^{-1} \lambda^{-3} \big) }  \bigg[  \nu \psi^{\mathrm{AT-F}} \big( z^{\prime} - \frac{1}{2} \big) + \lambda^3 \bar{\psi^{\mathrm{AT-F}} \big( z^{\prime} - \frac{1}{2} \big)}       -   \nu^{-1} \psi^{\mathrm{AT-F}} \big( z^{\prime\prime} - \frac{1}{2} \big) - \cdots \\ \lambda^{-3} \bar{\psi^{\mathrm{AT-F}} \big( z^{\prime\prime} - \frac{1}{2} \big) } \bigg]  =    \bar{\psi^{\mathrm{AT-F}}\big( z - \frac{1}{2} \big) }\text{ } \text{ , } 
\end{align*}

\noindent and,

\begin{align*}
     \frac{1}{\big( \nu \lambda^3 - \nu^{-1} \lambda^{-3} \big) }  \bigg[  \nu \psi^{\mathrm{AT-F}} \big( z^{\prime} +  \frac{1}{2} \big) + \lambda^3 \bar{\psi^{\mathrm{AT-F}} \big( z^{\prime} + \frac{1}{2} \big)}       -   \nu^{-1} \psi^{\mathrm{AT-F}} \big( z^{\prime\prime} + \frac{1}{2} \big) - \cdots \\ \lambda^{-3} \bar{\psi^{\mathrm{AT-F}} \big( z^{\prime\prime} + \frac{1}{2} \big) } \bigg] =    \bar{\psi^{\mathrm{AT-F}}\big( z + \frac{1}{2} \big) }\text{ } \text{ , } 
\end{align*}

\noindent corresponding to points which lie on the $\frac{1}{2}$ interval in the horizontal direction,

\begin{align*}
     \frac{1}{\big( \nu \lambda^3 - \nu^{-1} \lambda^{-3} \big) }  \bigg[  \nu \psi^{\mathrm{AT-F}} \big( z^{\prime} - \frac{i}{2} \big) + \lambda^3 \bar{\psi^{\mathrm{AT-F}} \big( z^{\prime} - \frac{i}{2} \big)}       -   \nu^{-1} \psi^{\mathrm{AT-F}} \big( z^{\prime\prime} - \frac{i}{2} \big) - \cdots \\ \lambda^{-3} \bar{\psi^{\mathrm{AT-F}} \big( z^{\prime\prime} - \frac{i}{2} \big) } \bigg]  =    \bar{\psi^{\mathrm{AT-F}}\big( z - \frac{i}{2} \big) }\text{ } \text{ , } 
\end{align*}

\noindent and,

\begin{align*}
     \frac{1}{\big( \nu \lambda^3 - \nu^{-1} \lambda^{-3} \big) }  \bigg[  \nu \psi^{\mathrm{AT-F}} \big( z^{\prime} + \frac{i}{2} \big) + \lambda^3 \bar{\psi^{\mathrm{AT-F}} \big( z^{\prime} + \frac{i}{2} \big)}       -   \nu^{-1} \psi^{\mathrm{AT-F}} \big( z^{\prime\prime} + \frac{i}{2} \big) - \cdots \\ \lambda^{-3} \bar{\psi^{\mathrm{AT-F}} \big( z^{\prime\prime} + \frac{i}{2} \big) } \bigg] =    \bar{\psi^{\mathrm{AT-F}}\big( z + \frac{i}{2} \big) }\text{ } \text{ , } 
\end{align*}

\noindent corresponding to points which lie on the $\frac{1}{2}$ interval in the vertical direction. From each possibility listed above, it is straightforward to verify that an application of $J$ to each equation yields,

\begin{align*}
    \frac{J}{\big(     \bar{\nu^{-1}} \bar{\lambda } - \bar{\nu} \bar{\lambda^{-1}}        \big)  }    \bigg[    \bar{\nu^{-1}} \bar{\psi^{\mathrm{AT-F}} \big( z^{\prime} \big)} + \bar{\lambda } \bar{\bar{\psi^{\mathrm{AT-F}} \big( z^{\prime} \big) }}    - \bar{\nu} \bar{ \psi^{\mathrm{AT-F} } \big( z^{\prime\prime} \big)} - \bar{\lambda^{-1}} \bar{ \bar{\psi^{\mathrm{AT-F}} \big( z^{\prime\prime} \big) } }   \bigg]       = J \bigg[ \psi^{\mathrm{AT-F}} \big( z \big) \bigg] \\ \Updownarrow \\   \frac{1}{\big(     {\nu^{-1}} {\lambda } - {\nu} {\lambda^{-1}}        \big)  }    \bigg[    {\nu^{-1}} {\psi^{\mathrm{AT-F}} \big( z^{\prime} \big)} + {\lambda } \{\bar{\psi^{\mathrm{AT-F}} \big( z^{\prime} \big) }    - {\nu} { \psi^{\mathrm{AT-F} } \big( z^{\prime\prime} \big)} - {\lambda^{-1}} { \bar{\psi^{\mathrm{AT-F}} \big( z^{\prime\prime} \big) } }   \bigg]       =  \bar{\psi^{\mathrm{AT-F}} \big( z \big) } \text{ } \text{ , } 
\end{align*}

\noindent corresponding to the expression obtained from the first and second s-holomorphicity equations for $\psi^{\mathrm{AT-F}} \big( z \big)$,

\begin{align*}
    \frac{J}{\big( \bar{\nu} \bar{\lambda^3} - \bar{\nu^{-1}} \bar{\lambda^{-3}}\big) }  \bigg[ \bar{\nu} \bar{\psi^{\mathrm{AT-F}} \big( z^{\prime} \big)} + \bar{\lambda^3 } \bar{\bar{\psi^{\mathrm{AT-F}} \big( z^{\prime} \big)}}       -  \bar{ \nu^{-1} } \bar{\psi^{\mathrm{AT-F}} \big( z^{\prime\prime} \big) } - \bar{\lambda^{-3} } \bar{\bar{\psi^{\mathrm{AT-F}} \big( z^{\prime\prime}  \big) }   }     \bigg] =  J \bigg[ \psi^{\mathrm{AT-F}} \big( z + 1 \big) \bigg] \\ \Updownarrow \\   \frac{1}{\big( {\nu} {\lambda^3} - {\nu^{-1}} {\lambda^{-3}}\big) }    \bigg[  {\nu} {\psi^{\mathrm{AT-F}} \big( z^{\prime} \big)} + {\lambda^3 } {\bar{\psi^{\mathrm{AT-F}} \big( z^{\prime} \big)}}       -  { \nu^{-1} } {\psi^{\mathrm{AT-F}} \big( z^{\prime\prime} \big) } - {\lambda^{-3} } {\bar{\psi^{\mathrm{AT-F}} \big( z^{\prime\prime}  \big) }   }              \bigg]         =     \bar{\psi^{\mathrm{AT-F}} \big( z + 1 \big) } \text{ } \text{ , } 
\end{align*}

\noindent corresponding to the expression obtained from the first and second s-holomorphicity equations for $\psi^{\mathrm{AT-F}} \big( z +1 \big)$,

\begin{align*}
     \frac{J}{\big( \nu \lambda^3 - \nu^{-1} \lambda^{-3} \big) }  \bigg[  \nu \psi^{\mathrm{AT-F}} \big( z^{\prime} - 1 \big) + \lambda^3 \bar{\psi^{\mathrm{AT-F}} \big( z^{\prime} -1 \big)}       -   \nu^{-1} \psi^{\mathrm{AT-F}} \big( z^{\prime\prime} -1 \big) -\cdots \\  \lambda^{-3} \bar{\psi^{\mathrm{AT-F}} \big( z^{\prime\prime} -1  \big) } \bigg] =   J \bigg[  \bar{\psi^{\mathrm{AT-F}}\big( z \big) } \bigg]  \\ \Updownarrow \\   \frac{1}{\big( \bar{\nu} \bar{\lambda^3} - \bar{\nu^{-1}} \bar{\lambda^{-3}} \big) }  \bigg[  \bar{\nu} \bar{\psi^{\mathrm{AT-F}} \big( z^{\prime} - 1 \big)}  + \bar{\lambda^3}  \bar{\bar{\psi^{\mathrm{AT-F}} \big( z^{\prime} -1 \big)}}       -   \bar{\nu^{-1} } \bar{\psi^{\mathrm{AT-F}} \big( z^{\prime\prime} -1 \big)  }  - \cdots   \\    \bar{\lambda^{-3} } \bar{\bar{\psi^{\mathrm{AT-F}} \big( z^{\prime\prime} -1  \big) }}  \bigg]         = \psi^{\mathrm{AT-F}}\big( z \big)   \text{ } \text{ , } 
\end{align*}

\noindent corresponding to the expression obtained from the first and second s-holomorphicity equations for $\bar{\psi^{\mathrm{AT-F}}\big( z \big) } $,

\begin{align*}
     \frac{J}{\big( \nu \lambda^3 - \nu^{-1} \lambda^{-3} \big) }  \bigg[  \nu \psi^{\mathrm{AT-F}} \big( z^{\prime} - \frac{1}{2} \big) + \lambda^3 \bar{\psi^{\mathrm{AT-F}} \big( z^{\prime} - \frac{1}{2} \big)}       -   \nu^{-1} \psi^{\mathrm{AT-F}} \big( z^{\prime\prime} - \frac{1}{2} \big) - \cdots \\ \lambda^{-3} \bar{\psi^{\mathrm{AT-F}} \big( z^{\prime\prime} - \frac{1}{2} \big) } \bigg] =   J \bigg[  \bar{\psi^{\mathrm{AT-F}}\big( z - \frac{1}{2} \big) } \bigg] \\ \Updownarrow \\                   \frac{1}{\big( \bar{\nu} \bar{\lambda^3} - \bar{\nu^{-1}} \bar{\lambda^{-3}} \big) }  \bigg[ \bar{\nu } \bar{\psi^{\mathrm{AT-F}} \big( z^{\prime} - \frac{1}{2} \big)} + \bar{\lambda^3 } \bar{\bar{\psi^{\mathrm{AT-F}} \big( z^{\prime} - \frac{1}{2} \big)} }      -   \bar{\nu^{-1} } \bar{\psi^{\mathrm{AT-F}} \big( z^{\prime\prime} - \frac{1}{2} \big)}  - \cdots \\ \bar{\lambda^{-3}} \bar{ \bar{\psi^{\mathrm{AT-F}} \big( z^{\prime\prime} - \frac{1}{2} \big) }}  \bigg] =  {\psi^{\mathrm{AT-F}}\big( z - \frac{1}{2} \big) }         \text{ } \text{ , } 
\end{align*}

\noindent corresponding to the expression obtained from the first and second s-holomorphicity equations for $\bar{\psi^{\mathrm{AT-F}}\big( z - \frac{1}{2} \big) }$,

\begin{align*}
     \frac{J}{\big( \nu \lambda^3 - \nu^{-1} \lambda^{-3} \big) }  \bigg[  \nu \psi^{\mathrm{AT-F}} \big( z^{\prime} +  \frac{1}{2} \big) + \lambda^3 \bar{\psi^{\mathrm{AT-F}} \big( z^{\prime} + \frac{1}{2} \big)}       -   \nu^{-1} \psi^{\mathrm{AT-F}} \big( z^{\prime\prime} + \frac{1}{2} \big) - \cdots \\ \lambda^{-3} \bar{\psi^{\mathrm{AT-F}} \big( z^{\prime\prime} + \frac{1}{2} \big) } \bigg] =   J \bigg[  \bar{\psi^{\mathrm{AT-F}}\big( z + \frac{1}{2} \big) } \bigg] \\ \Updownarrow \\                \frac{1}{\big( \bar{\nu} \bar{\lambda^3 }- \bar{\nu^{-1}} \bar{\lambda^{-3}} \big) }  \bigg[ \bar{\nu } \bar{\psi^{\mathrm{AT-F}} \big( z^{\prime} +  \frac{1}{2} \big)} + \bar{\lambda^3}  \bar{ \bar{\psi^{\mathrm{AT-F}} \big( z^{\prime} + \frac{1}{2} \big)} }      -   \bar{\nu^{-1} } \bar{\psi^{\mathrm{AT-F}} \big( z^{\prime\prime} + \frac{1}{2} \big) } - \cdots \\ \bar{\lambda^{-3}} \bar{ \bar{\psi^{\mathrm{AT-F}} \big( z^{\prime\prime} + \frac{1}{2} \big) }} \bigg] =   \psi^{\mathrm{AT-F}}\big( z + \frac{1}{2} \big)      \text{ } \text{ , } 
\end{align*}

\noindent corresponding to the expression obtained from the first and second s-holomorphicity equations for $\bar{\psi^{\mathrm{AT-F}}\big( z + \frac{1}{2} \big) }$,

\begin{align*}
     \frac{J}{\big( \nu \lambda^3 - \nu^{-1} \lambda^{-3} \big) }  \bigg[  \nu \psi^{\mathrm{AT-F}} \big( z^{\prime} - \frac{i}{2} \big) + \lambda^3 \bar{\psi^{\mathrm{AT-F}} \big( z^{\prime} - \frac{i}{2} \big)}       -   \nu^{-1} \psi^{\mathrm{AT-F}} \big( z^{\prime\prime} - \frac{i}{2} \big) - \cdots \\ \lambda^{-3} \bar{\psi^{\mathrm{AT-F}} \big( z^{\prime\prime} - \frac{i}{2} \big) } \bigg] =    J \bigg[ \bar{\psi^{\mathrm{AT-F}}\big( z - \frac{i}{2} \big) } \bigg] \\ \Updownarrow \\     \frac{1}{\big( \bar{\nu } \bar{\lambda^3 } - \bar{\nu^{-1}} \bar{ \lambda^{-3}} \big) }  \bigg[  \bar{\nu} \bar{\psi^{\mathrm{AT-F}} \big( z^{\prime} - \frac{i}{2} \big)}  + \bar{\lambda^3} \bar{ \bar{\psi^{\mathrm{AT-F}} \big( z^{\prime} - \frac{i}{2} \big)} }       -   \bar{\nu^{-1} } \bar{\psi^{\mathrm{AT-F}} \big( z^{\prime\prime} - \frac{i}{2} \big) } - \cdots \\ \bar{\lambda^{-3} } \bar{ \bar{\psi^{\mathrm{AT-F}} \big( z^{\prime\prime} - \frac{i}{2} \big) }} \bigg] =  {\psi^{\mathrm{AT-F}}\big( z - \frac{i}{2} \big) }    \text{ } \text{ , } 
\end{align*}

\noindent corresponding to the expression obtained from the first and second s-holomorphicity equations for $\bar{\psi^{\mathrm{AT-F}}\big( z - \frac{i}{2} \big) }$, and,

\begin{align*}
     \frac{J}{\big( \nu \lambda^3 - \nu^{-1} \lambda^{-3} \big) }  \bigg[  \nu \psi^{\mathrm{AT-F}} \big( z^{\prime} + \frac{i}{2} \big) + \lambda^3 \bar{\psi^{\mathrm{AT-F}} \big( z^{\prime} + \frac{i}{2} \big)}       -   \nu^{-1} \psi^{\mathrm{AT-F}} \big( z^{\prime\prime} + \frac{i}{2} \big) - \cdots \\ \lambda^{-3} \bar{\psi^{\mathrm{AT-F}} \big( z^{\prime\prime} + \frac{i}{2} \big) } \bigg] =   J \bigg[  \bar{\psi^{\mathrm{AT-F}}\big( z + \frac{i}{2} \big) } \bigg] \\ \Updownarrow \\    \frac{1}{\big( \bar{\nu} \bar{ \lambda^3}  - \bar{\nu^{-1}} \bar{ \lambda^{-3} } \big) }  \bigg[  \bar{\nu } \bar{\psi^{\mathrm{AT-F}} \big( z^{\prime} + \frac{i}{2} \big) }+ {\lambda^3} \bar{\bar{ \bar{\psi^{\mathrm{AT-F}} \big( z^{\prime} + \frac{i}{2} \big)}}}       -   \bar{\nu^{-1} } \bar{\psi^{\mathrm{AT-F}} \big( z^{\prime\prime} + \frac{i}{2} \big)}  - \cdots \\ \bar{\lambda^{-3} } \bar{\bar{\psi^{\mathrm{AT-F}} \big( z^{\prime\prime} + \frac{i}{2} \big) } }\bigg] =   \psi^{\mathrm{AT-F}}\big( z + \frac{i}{2} \big)     \text{ } \text{ , } 
\end{align*}

\noindent corresponding to the expression obtained from the first and second s-holomorphicity equations for $\bar{\psi^{\mathrm{AT-F}}\big( z + \frac{i}{2} \big) }$. We conclude the argument. \boxed{}

\subsubsection{Quantum correspondence with the Loop model}

\noindent We define the following. Over a finite volume of the hexagonal lattice which can be expressed with a Cartesian product of rows and columns, $\mathscr{I}_{\textbf{H}} \equiv \mathscr{I}$ and $\mathscr{J}_{\textbf{H}} \equiv \mathscr{J}$, which is embedded into the complex plane, introduce,

\begin{align*}
           \mathscr{P}^{\mathrm{loop},\xi}_{\Lambda_{\textbf{H}}} \big[ \sigma \big] = \frac{\mathrm{exp} \bigg[ \frac{1}{2} \big| \mathrm{log} x \big| \big(          h r \big( \sigma \big) + \mathrm{log} \big( x^{e(\sigma)} \big)      \big)  \bigg]}{\mathscr{Z}^{\mathrm{loop},\xi}_{\Lambda_{\textbf{H}}} \big( \sigma \big) }                     \text{ } \text{ , } 
\end{align*}

\noindent corresponding to the high temperature expansion of the loop $\mathrm{O} \big( 1 \big)$ measure, whose the partition function,

\begin{align*}
\mathscr{Z}^{\mathrm{loop},\xi}_{\Lambda_{\textbf{H}}} \big( \sigma \big)  \equiv \underset{\sigma_{\partial\Lambda_{\textbf{H}}} \equiv + 1}{\underset{e \in \Lambda_{\textbf{H}}}{\underset{    \sigma    \in  \{\pm 1 \}^{\mathscr{I} \times \mathscr{J}} }{\sum}}}     \mathrm{exp} \bigg[ \frac{1}{2} \big| \mathrm{log} x \big| \big(          h r \big( \sigma \big) + \mathrm{log} \big( x^{e(\sigma)} \big)      \big)  \bigg]             \text{ } \text{ , } 
\end{align*}

\noindent can be placed into correspondence with the state,

\begin{align*}
         \mathscr{Z}^{\mathrm{loop},\xi}_{\Lambda_{\textbf{H}}} \big( \sigma \big) \equiv  \mathscr{Z}^{\mathrm{loop},+}_{\Lambda_{\textbf{H}}} \big( \sigma \big) 
 \longleftrightarrow      \bra{f} V^{\mathrm{loop},N} \ket{i}                \text{ } \text{ , } 
\end{align*}

\noindent under $+$ boundary conditions $\xi \equiv +$. Along the lines of the discussion provided in the previous subsection for the spin inversion operation in the Ashkin-Teller model, there exists a spin inversion operator for the high- temperature expansion of the $\mathrm{O} \big( 1 \big)$ model, from which,

\begin{align*}
          \hat{\sigma_j } e_{\sigma} \equiv \sigma_j e_{\sigma}            \text{ } \text{ , } 
\end{align*}

\noindent for the operation from the loop state space into itself,

\begin{align*}
       \hat{\sigma_j}  :   \mathcal{S}^{\mathrm{loop}} \longrightarrow  \mathcal{S}^{\mathrm{loop}}  \text{ } \text{ , } 
\end{align*}

\noindent at site $j$, and,

\begin{align*}
   i \equiv \big( V^{\mathrm{loop}, h} \big)^{\frac{1}{2}} e_{(+)} \text{ } \text{ , } \\ f \equiv  i \equiv \big( V^{\mathrm{loop}, h} \big)^{\frac{1}{2}} e_{(+)} \text{ } \text{ . } 
\end{align*}

\noindent Furthermore, the expectation of a single spin with respect to $\mathscr{P}^{\mathrm{loop},\xi}_{\Lambda_{\textbf{H}}} \big( \cdot \big)$ takes the form,

\begin{align*}
  \mathscr{E}^{\mathrm{loop},\xi}_{\Lambda_{\textbf{H}}} \big(    \text{ }  \sigma_z     \text{ }   \big)  =    \underset{\sigma_{\partial\Lambda_{\textbf{H}}} \equiv + 1}{\underset{e \in \Lambda_{\textbf{H}}}{\underset{    \sigma    \in  \{\pm 1 \}^{\mathscr{I} \times \mathscr{J}} }{\sum}}}  \sigma_z      \text{ }  \mathrm{d}  \mathscr{P}^{\mathrm{loop},\xi}_{\Lambda_{\textbf{H}}} \big( \sigma_z \big)            =  \underset{\sigma_{\partial\Lambda_{\textbf{H}}} \equiv + 1}{\underset{e \in \Lambda_{\textbf{H}}}{\underset{    \sigma    \in  \{\pm 1 \}^{\mathscr{I} \times \mathscr{J}} }{\sum}}}  \frac{\sigma_z }{\mathscr{Z}^{\mathrm{loop},\xi}_{\Lambda_{\textbf{H}}} \big( \sigma \big) }     \text{ }      \mathrm{exp} \bigg[ \frac{1}{2} \big| \mathrm{log} x \big| \big(          h r \big( \sigma \big) + \mathrm{log} \big( x^{e(\sigma)} \big)      \big)  \bigg]                  \text{ } \text{ . } 
\end{align*}

\noindent Similarly, for respect to multiple spins at countably many sites in the complex plane,

\begin{align*}
    \mathscr{E}^{\mathrm{loop},\xi}_{\Lambda_{\textbf{H}}} \bigg[   \underset{z_i \in \textbf{C}}{\prod} \text{ }  \sigma_{z_i}      \bigg] =       \underset{\sigma_{\partial\Lambda_{\textbf{H}}} \equiv + 1}{\underset{e \in \Lambda_{\textbf{H}}}{\underset{    \sigma    \in  \{\pm 1 \}^{\mathscr{I} \times \mathscr{J}} }{\sum}}}  \bigg[    \underset{z_i \in \textbf{C}}{\prod} \text{ }  \sigma_{z_i}      \bigg]    \text{ }  \mathrm{d}  \mathscr{P}^{\mathrm{loop},\xi}_{\Lambda_{\textbf{H}}} \big( \sigma_z \big)            =  \underset{\sigma_{\partial\Lambda_{\textbf{H}}} \equiv + 1}{\underset{e \in \Lambda_{\textbf{H}}}{\underset{    \sigma    \in  \{\pm 1 \}^{\mathscr{I} \times \mathscr{J}} }{\sum}}}  \bigg[ {\mathscr{Z}^{\mathrm{loop},\xi}_{\Lambda_{\textbf{H}}} \big( \sigma \big) } \bigg]^{-1}  \bigg[ \text{ }   \underset{z_i \in \textbf{C}}{\prod} \text{ }  \sigma_{z_i}  \text{ }        \bigg]   \times \cdots \\    \text{ }       \mathrm{exp} \bigg[ \frac{1}{2} \big| \mathrm{log} x \big| \big(          h r \big( \sigma \big) + \mathrm{log} \big( x^{e(\sigma)} \big)      \big)  \bigg]           \text{ } \text{ . } 
\end{align*}

\noindent Under $+$ boundary conditions, the expectation over the product of spins is equivalent to,

\begin{align*}
     \frac{\bra{e_{(+)}} \underset{\mathrm{countably\text{ } many\text{ } } z_i}{\underset{z_i \in \textbf{C}}{\prod}} \sigma_{z_i} \ket{e_{(+)}} }{\bra{e_{(+)}} V^{\mathrm{loop},N} \ket{e_{(+)}}}          \text{ } \text{ . } 
\end{align*}

\noindent From the loop generators, similar relations are satisfied as those from the Ashkin-Teller generators, in which,

\begin{align*}
        p^{\mathrm{loop}}_u              \text{ } \text{ , } 
 \\   q^{\mathrm{loop}}_u                \text{ } \text{ , } 
\end{align*}

\noindent for the basis spanned by the elements,

\begin{align*}
     \psi^{\mathrm{loop}} \big( k + iy \big) = V^{\mathrm{loop},-y} \psi_k V^{\mathrm{loop},y}  \text{ } \text{ , }      \\     \bar{\psi^{\mathrm{loop}} \big( k + iy \big)} = \bar{V^{\mathrm{loop},-y}} \bar{\psi_k} \bar{V^{\mathrm{loop},y}} =   V^{\mathrm{loop},-y} \bar{\psi_k} V^{\mathrm{loop},y}    \text{ } \text{ , } 
\end{align*}

\noindent for $k \in \textbf{I}^{**}$ and $y \in \textbf{J}$.

\bigskip

\noindent From complexification procedure for the loop model, from the tuple $\big( \psi^{\mathrm{loop}} , \bar{\psi^{\mathrm{loop}}} \big)$, under the loop induced rotation,

\begin{align*}
    R \big(\psi^{\mathrm{loop}} \big(z \big) \big) = \bar{\psi^{\mathrm{loop}}} \big( r \big( z \big) \big) \text{ } \text{ , } \\ R \big( \psi^{\mathrm{loop}} \big( z \big)  \big)  = \bar{\psi^{\mathrm{loop}}} \big( r \big( z \big) \big)\text{ } \text{ , } 
\end{align*}

\noindent in which the image of the basis elements under $R$ satisfy the relations above. Similarly, under the other induced rotation $J$,

\begin{align*}
    J \big(\psi^{\mathrm{loop}} \big(z \big) \big) = \bar{\psi^{\mathrm{loop}}} \big( z  \big) \text{ } \text{ . } \\ J \big( \psi^{\mathrm{loop}} \big( z \big)  \big)  = \bar{\psi^{\mathrm{loop}}} \big(  z \big)\text{ } \text{ . } 
\end{align*}

\noindent With properties of the loop measure defined in this subsection, and relations that the loop generators satisfy with respect to the induced rotations $R$ and $J$, below we introduce the massive s-holomorphicity result for the loop fermion.

\bigskip

\noindent \textbf{Theorem} $3^{*}$ (\textit{massive s-holmorphicity for loop fermions}). Fix $z \in \textbf{I}^{**} \times \textbf{J}$, and the same parameters $\nu$ and $\lambda$ provided in \textbf{Definition} \textit{7}, and \textbf{Definition} \textit{8}. For the loop fermion operators $\psi^{\mathrm{loop-F}} \big( z \big) $ and $\bar{\psi^{\mathrm{loop-F}}} \big( z \big)$, there exists extensions of $\psi^{\mathrm{loop-F}}$, and of $\bar{\psi^{\mathrm{loop-F}}}$, to $\textbf{I}^{**} \times \textbf{J}$, such that, 

\begin{align*}
       \psi^{\mathrm{loop-F}} \big( z_1 \big) + \bar{e_1}^{2s} \bar{\psi^{\mathrm{loop-F}} \big( z_1 \big) } = \psi^{\mathrm{loop-F}} \big( z_2 \big) + \bar{e_1}^{2s} \bar{\psi^{\mathrm{loop-F}} \big( z_2 \big) }      \text{ } \text{ , }   \\  
     \psi^{\mathrm{loop-F}} \big( z_2 \big) + \bar{e_2}^{2s} \bar{\psi^{\mathrm{loop-F}} \big( z_2 \big)} = \psi^{\mathrm{loop-F}} \big( z_3 \big) + \bar{e_2}^{2s} \bar{\psi^{\mathrm{loop-F}} \big( z_3 \big)} \text{ } \text{ , }  \\    \psi^{\mathrm{loop-F}} \big( z_3 \big) + \bar{e_3}^{2s} \bar{\psi^{\mathrm{loop-F}} \big( z_3 \big)} = \psi^{\mathrm{loop-F}} \big( z_4 \big) + \bar{e_3}^{2s} \bar{\psi^{\mathrm{loop-F}} \big( z_4 \big) }     \text{ } \text{ , }  \\ \psi^{\mathrm{loop-F}} \big( z_4 \big) + \bar{e_4}^{2s} \bar{\psi^{\mathrm{loop-F}} \big( z_4 \big) } = \psi^{\mathrm{loop-F}} \big( z_1 \big) + \bar{e_4}^{2s} \bar{\psi^{\mathrm{loop-F}} \big( z_1 \big) } \text{ } \text{ , }
\end{align*}

\noindent for any face of the complex plane with edges E,N,W,S. At the left, and right, boundary points of the finite volume, the fermion operators also satisfy,

\begin{align*}
  \psi^{\mathrm{loop-F}} \big( a + i y \big) + i \bar{\psi^{\mathrm{loop-F}} } \big( a + iy \big) = 0   \text{ } \text{ , } \\     \psi^{\mathrm{loop-F}} \big( b + i y \big) + i \bar{\psi^{\mathrm{loop-F}} } \big( b + iy \big) = 0            \text{ } \text{ , } 
\end{align*}

\noindent for all $y \in \textbf{J}^{**}$.

\bigskip

\noindent \textit{Proof of Theorem} $3^{*}$. Appeal to arguments provided in \textbf{Theorem} \textit{3} of the previous subsection, by making use of s-holomorphicity equations provided for the loop fermion observable $\psi^{\mathrm{loop-F}}$ above. \boxed{}

\subsection{Low temperature expansions of parafermionic observables}

\noindent With the massive s-holomorphicity result for fermion operators of the previous section, in this section we define the parafermionic observable in the vanishing temperature limit, for parameters $\alpha^{\mathrm{AT}} \equiv \mathrm{exp} \big(- 2 \beta \big)$, and $\alpha^{\mathrm{loop}} \equiv \mathrm{exp} \big(- 2  \beta^{\mathrm{loop}} \big)$.

\subsubsection{Ashkin-Teller model}

\noindent For the first model, as $\beta \longrightarrow + \infty$ from below, introduce the following two observables related to the Ashkin-Teller parafermionic observable,

\begin{align*}
  f^{\mathrm{AT},\uparrow}_a \big( z \big)    = \frac{1}{\mathscr{C} \big( \mathscr{Z}^{\xi}_{\Lambda} \big) \mathscr{Z}^{\mathrm{low-temp},\xi}_{\Lambda}} \text{ }   \underset{\gamma \in \mathcal{C}^{\mathrm{AT},\uparrow}_a ( z )}{\sum}   \text{ }    \big( \alpha^{\mathrm{AT}} \big)^{L ( \gamma_i )}  \mathrm{exp} \big( - i \sigma_m \theta_{\gamma} \big( r , \vec{r} \big) \big) s_m \big( r \big) \mu_m \big( r \big)                    \text{ } \text{ , } \\          f^{\mathrm{AT},\downarrow}_a \big( z \big)          = \frac{1}{\mathscr{C} \big( \mathscr{Z}^{\xi}_{\Lambda} \big) 
 \mathscr{Z}^{\mathrm{low-temp},\xi}_{\Lambda}}   \text{ }     \underset{\gamma \in \mathcal{C}^{\mathrm{AT},\downarrow}_a ( z )}{\sum}   \text{ }     \big( \alpha^{\mathrm{AT}} \big)^{L ( \gamma_i )}   \mathrm{exp} \big( - i \sigma_m \theta_{\gamma} \big( r , \vec{r} \big) \big) s_m \big( r \big) \mu_m \big( r \big)                 \text{ } \text{ , } 
\end{align*}

\noindent for:

\begin{itemize}
    \item[$\bullet$] (1): A horizontal edge $a \in \mathscr{I}^{*} \times \mathscr{J}$,

     \item[$\bullet$] (2): a rectangle $\mathscr{I} \times \mathscr{J}$,

      \item[$\bullet$] (3): the dual lattice $\big( \textbf{Z}^2 \big)^{*} = \textbf{Z}^2 + \big( \frac{1}{2} , \frac{1}{2} \big)$,

       \item[$\bullet$] (4): a subset $\bigg[ \mathscr{I}^{*} \times \mathscr{J}^{*} \bigg] \subsetneq \bigg[ \textbf{Z}^2 + \big( \frac{1}{2} ,  \frac{1}{2} \big) \bigg] $,

    \item[$\bullet$] (5): a set $\mathcal{C}^{\mathrm{AT},\uparrow}_a \big( z \big)$, consisting of paths such that for any face of the dual lattice excluding points $a + \frac{i}{2}$, the parity, ie the number of edges of the path, adjacent to the face is even,

      \item[$\bullet$] (6): a set $\mathcal{C}^{\mathrm{AT},\downarrow}_a \big( z \big)$, consisting of paths such that for any face of the dual lattice excluding points $a + \frac{i}{2}$, the parity, ie the number of edges of the path, adjacent to the face is odd,

      \item[$\bullet$] (7): a strictly positive constant $\mathscr{C}$ so that $\mathscr{C} \big( \mathscr{Z}^{\xi}_{\Lambda} \big) 
 \mathscr{Z}^{\mathrm{low-temp},\xi}_{\Lambda} \propto \mathscr{Z}^{\xi}_{\Lambda}$.
\end{itemize}

\noindent \textbf{Definition} \textit{11} (\textit{discrete residues and massive s-holomorphicity}, \textbf{Definition} \textit{20}, [13]). For a horizontal edge $a$, there exists a complex valued, massive s-holomorphic function $f$, such that, for $z \neq a$, over the nonempty collection of faces of the lattice containing $a + \frac{i}{2}$, the residue of $f$ is given by $\frac{i}{2\pi} \big( f^{\mathrm{front}} \big( a \big) - f^{\mathrm{back}} \big( a \big) \big)$. The function $f$ can be extended so that it is massive s-holomorphic if $f^{\mathrm{front}} \big( \cdot \big)$ is extended to $a+\frac{i}{2}$, while $f^{\mathrm{back}} \big( \cdot \big)$ is extended to $a - \frac{i}{2}$.

\bigskip

\noindent \textbf{Theorem} \textit{4} (\textit{convergence of Fermion two-point correlations to Ashkin-Teller parafermionic observables defined at the beginning of 2.5.1}). Over $\mathscr{I} \times \mathscr{J}$, one has,

\begin{align*}
        <    \psi^{\mathrm{AT-F}} \big( z \big)     \psi^{\mathrm{AT-F}} \big( a \big)    >_{\mathscr{I} \times \mathscr{J}} =               - f^{\mathrm{AT},\uparrow}_a \big( z \big) + i    f^{\mathrm{AT},\downarrow}_a \big( z \big)     \text{ } \text{ , } \\  < \psi^{\mathrm{AT-F}} \big( z \big)     \bar{\psi^{\mathrm{AT-F}} \big( a \big) }  >_{\mathscr{I} \times \mathscr{J}}  =      f^{\mathrm{AT},\uparrow}_a \big( z \big) + i    f^{\mathrm{AT},\downarrow}_a \big( z \big)   \equiv      - \bar{f^{\mathrm{AT},\uparrow}_a \big( z \big)} - i \bar{f^{\mathrm{AT},\downarrow}_a \big( z \big) }         \text{ } \text{ , } \\  < \bar{\psi^{\mathrm{AT-F}} \big( z \big)}     \bar{\psi^{\mathrm{AT-F}} \big( a \big)} >_{\mathscr{I} \times \mathscr{J}}  =              - \bar{f^{\mathrm{AT},\uparrow}_a \big( z \big)} - i \bar{f^{\mathrm{AT},\downarrow}_a \big( z \big) }          \text{ } \text{ . } 
\end{align*}

\noindent \textit{Proof of Theorem 4}. Fix $z = x+iy$ and $z^{\prime} = x^{\prime} + i y^{\prime}$. By the massive s-holomorphicity of $f^{\mathrm{AT},\uparrow}_a \big( \cdot \big)$, and of $f^{\mathrm{AT},\downarrow}_a \big( \cdot \big)$, compute the expansion of the following state from a product of transfer matrices,

\begin{align*}
  \bra{e_{(+)}}  V^{\mathrm{AT},N-y} \psi^{\mathrm{AT-F}}_x V^{\mathrm{AT}, y- y^{\prime}}       \bar{\psi^{\mathrm{AT-F}_{x^{\prime}}} }    V^{\mathrm{AT}, y^{\prime}}       \ket{e_{(+)}} \text{ } \text{ , } 
\end{align*}

\noindent for $x \neq x^{\prime}$. Before performing the computation, observe,

\begin{align*}
  \bra{e_{(+)}}  V^{\mathrm{AT},N-y} \psi^{\mathrm{AT-F}}_x V^{\mathrm{AT}, y- y^{\prime}}       \bar{\psi^{\mathrm{AT-F}_{x^{\prime}}} }    V^{\mathrm{AT}, y^{\prime}}       \ket{e_{(+)}} \neq 0 \Longleftrightarrow     \exists\text{ }  \hat{\sigma_i} \text{ } \mathrm{st} \text{ }  \tau =  \underset{i < x^{\prime}}{\prod}    \big( -1 \big)^{i} \hat{\sigma_i} \text{ } , \text{ } \mathrm{or} \text{ }    \tau =  \underset{i > x^{\prime}}{\prod}   \big( -1 \big)^{i} \hat{\sigma_i}  \text{ } \text{ , } 
\end{align*}

\noindent from which we write,

\begin{align*}
         \bra{e_{(+)}}  V^{\mathrm{AT},N-y} \psi^{\mathrm{AT-F}}_x V^{\mathrm{AT}, y- y^{\prime}}       \bar{\psi^{\mathrm{AT-F}_{x^{\prime}}} }    V^{\mathrm{AT}, y^{\prime}}       \ket{e_{(+)}}  \propto               \underset{\sigma}{\sum} \bigg[  V^{\mathrm{AT}}_{(+), \sigma^{(N-1)}}\text{ }  \bigg[ \underset{N-1 \leq i \leq y }{\prod}   \text{ } V^{\mathrm{AT}}_{\sigma(i)} \text{ }   \bigg]   \text{ } V^{\mathrm{AT}}_{\sigma(y+1) , \tau(y) } \text{ }  \times \cdots \\  
        \big( \psi_x \big)_{\tau(y) , \sigma(y) } \text{ }  \bigg[    \underset{y \leq i^{\prime} \leq     y^{\prime} }{\prod}     V^{\mathrm{AT}}_{\sigma(i^{\prime})}        \text{ }   \bigg]  \text{ }        \big( \bar{\psi_x} \big)_{\tau(y^{\prime}) , \sigma(y^{\prime})}   \text{ }      V^{\mathrm{AT}}_{\sigma(1),(+)}          \bigg]                    \text{ } \text{ . } 
\end{align*}

\noindent Proceeding with the computation, write the following expression for the entries of the Askin-Teller transfer matrix, in which,

\begin{align*}
   V^{\mathrm{AT}} \equiv     \big( V^{\mathrm{AT},h} \big)^{\frac{1}{2}} V^{\mathrm{AT},V }            \big( V^{\mathrm{AT},h} \big)^{\frac{1}{2}}           \text{ } \text{ , }
\end{align*}

\noindent satisfies the relation,

\begin{align*}
  V^{\mathrm{AT}} \propto            \bigg[   \mathrm{exp} \big[ \text{ }  J^{*} \big(  \underset{k \in \textbf{I}}{ \sum}  p^{\mathrm{AT}}_{k- \frac{1}{2}} q^{\mathrm{AT}}_{k- \frac{1}{2}}   +  \big( p^{\mathrm{AT}}_{k- \frac{1}{2}}\big)^{\prime}    \big( q^{\mathrm{AT}}_{k- \frac{1}{2}} \big)^{\prime}   \big)   \big]  +  \mathrm{exp} \big[  \text{ }  U^{*} \big(  \underset{k \in \textbf{I}}{ \sum}              p^{\mathrm{AT}}_{k-\frac{1}{2}}  \big( p^{\mathrm{AT}}_{k-\frac{1}{2}} \big)^{\prime} q^{\mathrm{AT}}_{k-\frac{1}{2}}  \big( q^{\mathrm{AT}}_{k- \frac{1}{2}} \big)^{\prime}   \big)  \big]  \bigg]                                  \times \cdots \\  \mathrm{exp} \bigg[ \text{ }  J \big(  \underset{k \in \textbf{I}^{*}}{ \sum}  p^{\mathrm{AT}}_kq^{\mathrm{AT}}_k   +  \big( p^{\mathrm{AT}}_k\big)^{\prime}    \big( q^{\mathrm{AT}}_k\big)^{\prime}   \big)   \bigg]  +  \mathrm{exp} \bigg[  \text{ }  U \big(  \underset{k \in \textbf{I}^{*}}{ \sum}              p^{\mathrm{AT}}_k  \big( p^{\mathrm{AT}}_k \big)^{\prime} q^{\mathrm{AT}}_k  \big( q^{\mathrm{AT}}_k \big)^{\prime}   \big)  \bigg] \times \cdots \\   \bigg[   \mathrm{exp} \big[ \text{ }  J^{*} \big(  \underset{k \in \textbf{I}}{ \sum}  p^{\mathrm{AT}}_{k- \frac{1}{2}} q^{\mathrm{AT}}_{k- \frac{1}{2}}   +  \big( p^{\mathrm{AT}}_{k- \frac{1}{2}}\big)^{\prime}    \big( q^{\mathrm{AT}}_{k- \frac{1}{2}} \big)^{\prime}   \big)   \big]  +  \mathrm{exp} \big[  \text{ }  U^{*} \big(  \underset{k \in \textbf{I}}{ \sum}              p^{\mathrm{AT}}_{k-\frac{1}{2}}  \big( p^{\mathrm{AT}}_{k-\frac{1}{2}} \big)^{\prime} q^{\mathrm{AT}}_{k-\frac{1}{2}}  \big( q^{\mathrm{AT}}_{k- \frac{1}{2}} \big)^{\prime}   \big)  \big]  \bigg]                  \text{ } \text{ . } 
\end{align*}

\noindent From the other expectation value,

\begin{align*}
      \mathscr{E}^{\prime} \equiv    \bra{e_{(+)}}  V^{\mathrm{AT},N-y} \psi^{\mathrm{AT-F}}_x V^{\mathrm{AT}, y- y^{\prime}}       \bar{\psi^{\mathrm{AT-F}_{x^{\prime}}} }    V^{\mathrm{AT}, y^{\prime}}       \ket{e_{(+)}} \text{ } \text{ , } \end{align*}

\noindent observe,

\begin{align*}
       \mathscr{E}^{\prime}  \propto     i \underset{\gamma \in \mathcal{C}^{\uparrow}_a (z ) \text{ } \cup \text{ } \mathcal{C}^{\downarrow}_a (z) }{\sum}    \big( \alpha^{\mathrm{AT}} \big)^{L ( \gamma_i )} \big( - 1 \big)^{\# \{  \gamma_i \cap \textbf{I}^{>x}_y \} - \# \{ \gamma_i \cap \textbf{I}^{> x^{\prime}}_y \}}  \mathrm{exp} \big( - i \sigma_m \theta_{\gamma} \big( r , \vec{r} \big) \big) s_m \big( r \big)\text{ } \text{ . } 
\end{align*}

\noindent Furthermore, 

\begin{align*}
     i \underset{\gamma \in \mathcal{C}^{\uparrow}_a (z ) \text{ } \cup \text{ } \mathcal{C}^{\downarrow}_a (z) }{\sum}    \big( \alpha^{\mathrm{AT}} \big)^{L ( \gamma_i )} \big( - 1 \big)^{\# \{  \gamma_i \cap \textbf{I}^{>x}_y \} - \# \{ \gamma_i \cap \textbf{I}^{> x^{\prime}}_y \}}  \mathrm{exp} \big( - i \sigma_m \theta_{\gamma} \big( r , \vec{r} \big) \big) s_m \big( r \big) \propto   i   \underset{\gamma \in \mathcal{C}^{\mathrm{AT},\downarrow}_a ( z )}{\sum}   \text{ }     \big( \alpha^{\mathrm{AT}} \big)^{L ( \gamma_i )}   \times \cdots \\ \mathrm{exp} \big( - i \sigma_m \theta_{\gamma} \big( r , \vec{r} \big) \big) s_m \big( r \big) \mu_m \big( r \big)                   \text{ } \text{ . } 
\end{align*}

\noindent Hence,

\begin{align*}
        <    \psi^{\mathrm{AT-F}} \big( z \big)     \psi^{\mathrm{AT-F}} \big( a \big)    >_{\mathscr{I} \times \mathscr{J}} =               - f^{\mathrm{AT},\uparrow}_a \big( z \big) + i    f^{\mathrm{AT},\downarrow}_a \big( z \big)        \text{ } \text{ , } 
\end{align*}

\noindent for $y > y^{\prime}$, while, 

\begin{align*}
   \mathscr{E}^{\prime} \propto < \bar{\psi^{\mathrm{AT-F}} \big( z \big)}     \bar{\psi^{\mathrm{AT-F}} \big( a \big)} >_{\mathscr{I} \times \mathscr{J}}  =              - \bar{f^{\mathrm{AT},\uparrow}_a \big( z \big)} - i \bar{f^{\mathrm{AT},\downarrow}_a \big( z \big) }        \text{ } \text{ , } 
\end{align*}

\noindent for $y < y^{\prime}$, from which we conclude the argument. \boxed{}

\bigskip

\noindent \textbf{Theorem} \textit{5} (\textit{Pffaffian from multi-point correlations of the Ashkin-Teller fermion operator}, \textbf{Theorem} \textit{23}, [13]). From the Ashkin-Teller fermion operator, $\psi^{\mathrm{AT-F}}$,

\begin{align*}
   <  \underset{\mathrm{countably\text{ } many \text{ }} z_i}{\underset{z_i \in \textbf{C}}{\underset{1 \leq i \leq n}{\prod}}} \psi^{(\mathrm{AT}-F),(i)}_{z_i}      >^{+}_{\mathscr{I} \times \mathscr{J}} = \underset{1 \leq i , j \leq n}{\bigcup} \mathrm{Pf} \bigg[   <   \psi^{(\mathrm{AT-F}),(i)} \big( z_i \big)  \psi^{(\mathrm{AT-F}),(j)} \big( z_j \big)  >_{\mathscr{I} \times \mathscr{J}}           \bigg]  \text{ } \text{ , } 
\end{align*}

\noindent for the expectation with respect to $+$ boundary conditions supported over $_{\mathscr{I} \times \mathscr{J}}$,

\begin{align*}
  < \cdot >^{+}_{\mathscr{I} \times \mathscr{J}}    \text{ } \text{ . } 
\end{align*}

\noindent \textit{Proof of Theorem 5}. By the polarization result in \textit{2.3.1}, \textbf{Lemma} \textit{10}, which works for almost all $\beta$ except for countably many temperatures. Independently of whether $J  \equiv U \equiv \frac{1}{4} \mathrm{log} \big( 3 \big)$, or $J,U < \frac{1}{4} \mathrm{log} \big( 3 \big)$, write,

\begin{align*}
    \big( v^{\prime,(+),N}_{\mathrm{vac}}\big)^{*}  \equiv   \frac{e^{\prime,\mathrm{T}}_{(+)} V^N }{e^{\prime,\mathrm{T}}_{(+)} V^N e^{\prime}_{(+)}} \text{ } \text{ , } 
\end{align*}

\noindent from \textbf{Lemma} $13^{*}$, as the image of some map $\phi$ from the state space into the wedge product of Fock space representations. From the fact that individual terms of quantum states can be expressed in terms of inner and outer products,

\begin{align*}
\bra{e_{(+)}}    V^{\mathrm{AT},N}     \ket{e_{(+)}} \longleftrightarrow              \big( e_{(+)} \big)^{\textbf{T}} V^{\mathrm{AT},N}    e_{(+)}  \text{ } \text{ , } 
\end{align*}

\noindent write,

\begin{align*}
     <  \underset{\mathrm{countably\text{ } many \text{ }} z_i}{\underset{z_i \in \textbf{C}}{\underset{1 \leq i \leq n}{\prod}}} \psi^{(\mathrm{AT}-F),(i)}_{z_i}      >^{+}_{\mathscr{I} \times \mathscr{J}} \equiv \textbf{E}^{+}_{\mathscr{I} \times \mathscr{J}} \bigg[   \underset{\mathrm{countably\text{ } many \text{ }} z_i}{\underset{z_i \in \textbf{C}}{\underset{1 \leq i \leq n}{\prod}}} \psi^{(\mathrm{AT}-F),(i)}_{z_i}        \bigg]   \\ \equiv    \bigg[ \textbf{E}^{+}_{\mathscr{I} \times \mathscr{J}} \bigg]^{\prime} \bigg[  \big( v_{\mathrm{vac}} \big)^{*} \bigg[ \underset{\mathrm{countably\text{ } many \text{ }} z_i}{\underset{z_i \in \textbf{C}}{\underset{1 \leq i \leq n}{\prod}}} \psi^{(\mathrm{AT}-F),(i)}_{z_i}  \text{ } \bigg]   v_{\mathrm{vac}}    \bigg]     \\    \overset{(\textbf{Lemma} \text{ } \textit{11})}{\equiv}                 <           \big(  v_{\mathrm{vac}}  \big)^{*}            \bigg[ \underset{\mathrm{countably\text{ } many \text{ }} z_i}{\underset{z_i \in \textbf{C}}{\underset{1 \leq i \leq n}{\prod}}} \psi^{(\mathrm{AT}-F),(i)}_{z_i}  \text{ } \bigg]          v_{\mathrm{vac}}>             \text{ } \text{ , } 
\end{align*}

\noindent where,

\begin{align*}
   \frac{\textbf{E}^{*}_{\mathscr{I}  \times \mathscr{J}}  \big[ \text{ }  \cdot \text{ }  \big] }{                   \bigg[ \textbf{E}^{*}_{\mathscr{I}  \times \mathscr{J}} \bigg]^{\prime}   \big[ \text{ }  \cdot \text{ }  \big]  } \propto \underset{i \sim j}{\sum} \bigg[  J \bigg[  \tau \big( i \big) \tau \big( j \big)  +  \tau^{\prime} \big( i \big) \tau^{\prime} \big( j \big) \big)  - \big( v_{\mathrm{vac}} \big)^{*} \tau \big( i \big) \tau \big( j \big) v_{\mathrm{vac}} -  \big( v_{\mathrm{vac}} \big)^{*} \tau^{\prime} \big( i \big) \tau^{\prime} \big( j \big) v_{\mathrm{vac}}  \bigg] + \cdots \\            U \bigg[   \tau \big( i \big) \tau \big( j \big) \tau^{\prime} \big( i \big) \tau^{\prime} \big( j \big) - \big( v_{\mathrm{vac}} \big)^{*} \tau \big( i \big) \tau \big( j \big) \tau^{\prime}  \big( j \big) v_{\mathrm{vac}} \bigg]              \bigg]                           \text{ } \text{ , } 
\end{align*}

\noindent from which we conclude the argument. \boxed{}

\bigskip

\noindent We also introduce the multipoint observable, and its connections with the Pfaffian (see \textit{4.5} of [13] for a more extensive overview). For the Ashkin-Teller model, in comparison to the parafermionic observable that is defined over a single point, over multiple points, the observable takes the form,

\begin{align*}
  \psi^{\mathrm{MP-AT}} \bigg[ \epsilon , \big( r_1 , \cdots , r_n \big)  ,  \big( \sigma_1 , \cdots , \sigma_n \big) , \big(  \vec{r}_1 , \cdots , \vec{r}_n \big)   \bigg] \equiv \underset{\gamma_ i \in \mathcal{C}_{\{r_1 , \cdots , r_n  \}}}{\sum}  \text{ } \underset{\gamma_ i \in \mathcal{C}_{\{r_1 , \cdots , r_n  \}}}{\prod} 
 \mathrm{exp} \big( - i \sigma_i \theta_{\gamma_i} \big( r_i , \vec{r}_i \big) \big) s_m \big( r_i \big) \mu_m \big( r_i \big)  \text{ , } 
\end{align*}

\noindent which is a function of the winding number of each path, with respective parameters $\sigma_1 , \cdots , \sigma_n$, from the collection of paths,

\begin{align*}
  \Gamma \equiv \underset{\gamma_i \in \mathcal{C}_{\{ r_1 , \cdots , r_n \}}}{\bigcup} \big\{ \mathrm{paths}\text{ } \gamma_i \text{ } | \text{ }     0 \leq \theta_{\gamma_i} \big( r_i , \vec{r_i} \big)   \leq 2 \pi    \big\}    \text{ } \text{ , }
\end{align*}

\noindent for:

\begin{itemize}
\item [$\bullet$] (1): Countably many points $r_1 , \cdots , r_n$ in the complex plane,

\item[$\bullet$] (2): the position vector, $\vec{r_1} , \cdots , \vec{r_n}$, of each $\gamma_i$,

\item[$\bullet$] (3): the spin of each $\gamma_i$, $s_m \big( r_i \big)$,

\item[$\bullet$] (4): the disorder variable of each $\gamma_i$, $\mu_m \big( r_i \big)$,

\item[$\bullet$] (5): the set of paths $\mathcal{C}_{\{r_1 , \cdots , r_n  \}}$, where each $\gamma_i$ beings at point $r_i$,

\item[$\bullet$] (6): a parameter $\epsilon$, whose form will be provided in \textbf{Theorem} \textit{6} below.

\end{itemize}

\noindent \textbf{Theorem} \textit{6} (\textit{expectation value of multipoint Ashkin-Teller fermion operators is equal to the Ashkin-Teller multipoint parafermionic observable}). Denote, for $\psi^{(\mathrm{MP-AT})} \equiv \psi$,

\begin{align*}
  \psi^{\uparrow} \big( z \big) = \frac{1}{2} \big( \bar{\psi^{\mathrm{AT-F}}} \big( z \big) - \psi^{\mathrm{AT-F}} \big( z \big) \big)  \text{ } \text{ , } 
\end{align*}

\noindent and,

\begin{align*}
    \psi^{\downarrow} \big( z \big)  = \frac{i}{2} \big( {\psi^{\mathrm{AT-F}}} \big( z \big) + \bar{\psi^{\mathrm{AT-F}}} \big( z \big) \big)  \text{ } \text{ , } 
\end{align*}

\noindent for $z \in \textbf{C}$, from which one has,

\begin{align*}
  < \underset{\mathrm{countably\text{ }many}\text{ } z_{(j)} \in \textbf{C}}{\underset{1 \leq j \leq 2m -1}{\underset{1 \leq i \leq 2m }{\prod}}} \psi^{\updownarrow , (2m-i)}_{(j)}   \big( z_{(j)} \big)         >^{+}_{\mathscr{I} \times \mathscr{J}} =    \psi^{\mathrm{MP-AT}}  \bigg[ \big( r_1 , \cdots , r_n \big)  ,  \big( \sigma_1 , \cdots , \sigma_n \big) , \big(  \vec{r}_1 , \cdots , \vec{r}_n \big)   \bigg]  \text{ } \text{ , } 
\end{align*}

\noindent for $\updownarrow \in \big\{ \uparrow , \downarrow \big\}$, and,

\begin{align*}
\epsilon   =  \lambda \text{ } \Longleftrightarrow  \text{ } \updownarrow \equiv \uparrow  \text{ , } \end{align*}

\begin{align*}
\epsilon = \frac{1}{\lambda^2}   \Longleftrightarrow  \text{ } \updownarrow \equiv \downarrow 
 \text{ } \text{ . } 
\end{align*}

\noindent \textit{Proof of Theorem 6}. To obtain the desired expressions for $\epsilon$, observe,

\begin{align*}
        \frac{1}{2} \big( \lambda^{-\eta_j} - i    \lambda^{\eta_j} \big)  = \frac{1}{\lambda} \delta_{\eta_j , +1 }  \text{ } \text{ , }   \\      \frac{1}{2} \big(  i \lambda^{-\eta_j} +   \lambda^{\eta_j} \big)  = \lambda^2 \delta_{\eta_j , +1 }           \text{ } \text{ , } 
\end{align*}

\noindent for real $\eta_j$. To demonstrate that the desired identity holds between the expectation under $+$ boundary conditions and the multipoint Ashkin-Teller parafermionic observable, write,

\begin{align*}
     < \underset{\mathrm{countably\text{ }many}\text{ } z_{(j)} \in \textbf{C}}{\underset{1 \leq j \leq 2m -1}{\underset{1 \leq i \leq 2m }{\prod}}} \psi^{\updownarrow , (2m-i)}_{(j)}   \big( z_{(j)} \big)         >^{+}_{\mathscr{I} \times \mathscr{J}} =   \textbf{E}^{+}_{\mathscr{I} \times \mathscr{J}} \bigg[  \underset{\mathrm{countably\text{ }many}\text{ } z_{(j)} \in \textbf{C}}{\underset{1 \leq j \leq 2m -1}{\underset{1 \leq i \leq 2m }{\prod}}} \psi^{\updownarrow , (2m-i)}_{(j)}   \big( z_{(j)} \big)          \bigg]  \\ =    \textbf{E}^{+}_{\mathscr{I} \times \mathscr{J}} \bigg[   \underset{1 \leq j \leq 2m-1}{\underset{1 \leq i \leq 2m }{\sum} } \bigg[ \underset{\mathrm{countably\text{ }many}\text{ } z_{(j)} \in \textbf{C}}{\underset{i , j }{\prod}} \psi^{\updownarrow , (2m-i)}_{(j)}   \big( z_{(j)} \big)       \bigg]       \bigg]   \\ =      \underset{\omega \in \mathcal{S}}{\underset{\mathscr{I} \times \mathscr{J} \subsetneq \textbf{Z}^2}{\underset{+\text{ }  \mathrm{boundary \text{ } conditions}}{\sum}}} \text{ }     \bigg[            \underset{1 \leq j \leq 2m-1}{\underset{1 \leq i \leq 2m }{\sum} } \bigg[  \underset{\mathrm{countably\text{ }many}\text{ } z_{(j)} \in \textbf{C}}{\underset{i , j }{\prod}} \psi^{\updownarrow , (2m-i)}_{(j)}   \big( z_{(j)} \big)            \bigg]    \bigg]  \text{ } \text{ , }
     \end{align*}
     
     \noindent from which terms can further be rearranged, by separating terms in the product over $i$ and $j$ for which $\psi^{\uparrow}$, or for which $\psi^{\downarrow}$, 
     
     \begin{align*}
     \underset{\omega \in \mathcal{S}}{\underset{\mathscr{I} \times \mathscr{J} \subsetneq \textbf{Z}^2}{\underset{+\text{ }  \mathrm{boundary \text{ } conditions}}{\sum}}} \text{ }     \bigg[           \underset{1 \leq j \leq 2m-1}{\underset{1 \leq i \leq 2m }{\sum} } \bigg[  \underset{\mathrm{countably\text{ }many}\text{ } z_{(j)} \in \textbf{C}}{\underset{i  : \psi^{\uparrow}  , j }{\prod}} \psi^{\uparrow , (2m-i)}_{(j)}   \big( z_{(j)} \big)   \underset{\mathrm{countably\text{ }many}\text{ } z_{(j)} \in \textbf{C}}{\underset{i : \psi^{\downarrow} , j }{\prod}} \psi^{\downarrow , (2m-i)}_{(j)}   \big( z_{(j)} \big)   \bigg]             \bigg]                               \text{ } \text{ , } 
\end{align*}

\noindent which can then be expressed with the Pfaffian,

\begin{align*}
      \underset{\omega \in \mathcal{S}}{\underset{\mathscr{I} \times \mathscr{J} \subsetneq \textbf{Z}^2}{\underset{+\text{ } \mathrm{boundary} \text{ } \mathrm{conditions}}{\sum}}} \text{ }     \bigg[         {\underset{1 \leq i \leq 2m}{\bigcup}}  \text{ } \mathrm{Pf} \bigg[   \underset{\mathrm{countably\text{ }many}\text{ } z_{(j)} \in \textbf{C}}{\underset{i  : \psi^{\uparrow}  , j }{\prod}} \psi^{\uparrow , (2m-i)}_{(j)}   \big( z_{(j)} \big)   \bigg]     \underset{1 \leq j \leq 2m-1}{\bigcup}  \mathrm{Pf} \bigg[ \underset{\mathrm{countably\text{ }many}\text{ } z_{(j)} \in \textbf{C}}{\underset{i : \psi^{\downarrow} , j }{\prod}} \psi^{\downarrow , (2m-i)}_{(j)}   \big( z_{(j)} \big)                     \bigg]           \bigg]  \text{ , } 
\end{align*}

\noindent which is equivalent to,

\begin{align*}
      \underset{\omega \in \mathcal{S}}{\underset{\mathscr{I} \times \mathscr{J} \subsetneq \textbf{Z}^2}{\underset{+\text{ }  \mathrm{boundary \text{ } conditions}}{\sum}}}      \bigg[            \underset{1 \leq j \leq 2m-1}{\underset{1 \leq i \leq 2m}{\bigcup}}  \text{ } \mathrm{Pf} \bigg[    \underset{\mathrm{countably\text{ }many}\text{ } z_{(j)} \in \textbf{C}}{\underset{i  : \psi^{\uparrow}  , j }{\prod}} \psi^{\uparrow , (2m-i)}_{(j)}   \big( z_{(j)} \big)   \underset{\mathrm{countably\text{ }many}\text{ } z_{(j)} \in \textbf{C}}{\underset{i : \psi^{\downarrow} , j }{\prod}} \psi^{\downarrow , (2m-i)}_{(j)}   \big( z_{(j)} \big)                     \bigg]          \bigg]  \\ \equiv    \underset{\omega \in \mathcal{S}}{\underset{\mathscr{I} \times \mathscr{J} \subsetneq \textbf{Z}^2}{\underset{+\text{ }  \mathrm{boundary \text{ } conditions}}{\sum}}} \bigg[            \underset{1 \leq j \leq 2m-1}{\underset{1 \leq i \leq 2m}{\bigcup}} 
      \underset{\mathrm{countably\text{ }many}\text{ } z_{(j)} \in \textbf{C}}{\underset{i  : \psi^{\uparrow}  , i : \psi^{\downarrow} }{\underset{j}{\prod}}}   \mathrm{Pf} \bigg[             \psi^{\updownarrow , (2m-i)}_{(j)} \big( z_{(j)}    \big)    \bigg]         \bigg]     \text{ } \text{ , } 
\end{align*}

\noindent For low temperatues as $\alpha^{\mathrm{AT}} \longrightarrow 0$, one recovers components of the multipoint Ashkin-Teller parafermionic observable from the Pfaffian, in which,

\begin{align*}
     \underset{\gamma_ i \in \mathcal{C}_{\{r_1 , \cdots , r_n  \}}}{\sum}  \text{ } \underset{\gamma_ i \in \mathcal{C}_{\{r_1 , \cdots , r_n  \}}}{\prod} 
 \mathrm{exp} \big( - i \sigma_i \theta_{\gamma_i} \big( r_i , \vec{r}_i \big) \big) s_m \big( r_i \big) \mu_m \big( r_i \big)           \text{ } \text{ , } 
\end{align*}

\noindent from which we conclude the argument. \boxed{}

\subsubsection{Loop model}

\noindent For the second model, as $\beta^{\mathrm{loop}} \longrightarrow + \infty$ from below, introduce the following two observables related to the loop parafermionic observable,

\begin{align*}
     f^{\mathrm{loop},\uparrow}_a \big( z \big)    = \frac{1}{\mathscr{C} \big( \mathscr{Z}^{\mathrm{
 loop},\xi}_{\Lambda} \big) \mathscr{Z}^{\mathrm{low-temp 
 \text{ } loop},\xi}_{\Lambda}}        \text{ }    \underset{\gamma \in \mathcal{C}^{\mathrm{loop},\uparrow}_a ( z )}{\underset{\gamma \subset \Omega}{\underset{\gamma : a \longrightarrow z}{\sum}}} \mathrm{exp} \big( - i \sigma W_{\gamma} \big( a , z \big) \big) x^{l ( \gamma ) }       \text{ } \text{ , } \\                  f^{\mathrm{loop},\downarrow }_a \big( z \big)   = \frac{1}{\mathscr{C} \big( \mathscr{Z}^{\mathrm{ 
 loop},\xi}_{\Lambda} \big) 
 \mathscr{Z}^{\mathrm{low-temp\text{ }  
 loop},\xi}_{\Lambda}}       \text{ }    \underset{\gamma \in \mathcal{C}^{\mathrm{loop},\downarrow}_a ( z )}{\underset{\gamma \subset \Omega}{\underset{\gamma : a \longrightarrow z}{\sum}}} \mathrm{exp} \big( - i \sigma W_{\gamma} \big( a , z \big) \big) x^{l ( \gamma ) }             \text{ } \text{ , }
\end{align*}

\noindent for:

\begin{itemize}
    \item[$\bullet$] (1): A horizontal edge $a \in \mathscr{I}^{**} \times \mathscr{J}$,

     \item[$\bullet$] (2): a rectangle $\mathscr{I}^{**} \times \mathscr{J}$,

      \item[$\bullet$] (3): the dual lattice $\big( \textbf{H} \big)^{*} = \textbf{T}$,

       \item[$\bullet$] (4): a subset $\bigg[ \mathscr{I}^{**} \times \mathscr{J}^{**} \bigg] \subsetneq \textbf{T}$,

        \item[$\bullet$] (5): a set $\mathcal{C}^{\mathrm{loop},\uparrow}_a \big( z \big)$, consisting of paths such that for any face of the dual lattice excluding points $a + \frac{i}{2}$, the parity, ie the number of edges of the path, adjacent to the face is even,

         \item[$\bullet$] (6): a set $\mathcal{C}^{\mathrm{loop},\downarrow}_a \big( z \big)$, consisting of paths such that for any face of the dual lattice excluding points $a + \frac{i}{2}$, the parity, ie the number of edges of the path, adjacent to the face is odd, 

          \item[$\bullet$] (7): a strictly positive constant $\mathscr{C} \big( \mathscr{Z}^{\mathrm{ 
 loop},\xi}_{\Lambda} \big) 
 \mathscr{Z}^{\mathrm{low-temp\text{ }  
 loop},\xi}_{\Lambda} \propto \mathscr{Z}^{\mathrm{loop},. \xi}_{\Lambda}$.
\end{itemize}

\noindent \textbf{Definition} \textit{12} (\textit{discrete residues and massive s-holomorphicity}, \textbf{Definition} \textit{20}, [13])). For a horizontal edge $a$, there exists a complex valued, massive s-holomorphic function $f$, such that, for $z \neq a$, over the nonempty collection of faces of the lattice containing $a + \frac{i}{2}$, the residue of $f$ is given by $\frac{i}{2\pi} \big( f^{\mathrm{front}} \big( a \big) - f^{\mathrm{back}} \big( a \big) \big)$. The function $f$ can be extended so that it is massive s-holomorphic if $f^{\mathrm{front}} \big( \cdot \big)$ is extended to $a+\frac{i}{2}$, while $f^{\mathrm{back}} \big( \cdot \big)$ is extended to $a - \frac{i}{2}$.

\bigskip

\noindent \textbf{Theorem} $4^{*}$ (\textit{convergence of Fermion two-point correlations to loop parafermionic observables defined at the beginning of 2.5.2}). Over $\mathscr{I} \times \mathscr{J}$, one has,

\begin{align*}
       <    \psi^{\mathrm{loop-F}} \big( z \big)     \psi^{\mathrm{loop-F}} \big( a \big)    >_{\mathscr{I} \times \mathscr{J}}  =            - f^{\mathrm{loop},\uparrow}_a \big( z \big) + i    f^{\mathrm{loop},\downarrow}_a \big( z \big)                            \text{ } \text{ , } \\  < \psi^{\mathrm{loop-F}} \big( z \big)     \bar{\psi^{\mathrm{loop-F}} \big( a \big) }  >_{\mathscr{I} \times \mathscr{J}} =                    - \bar{f^{\mathrm{loop},\uparrow}_a \big( z \big)} - i \bar{f^{\mathrm{loop},\downarrow}_a \big( z \big) }                     \text{ } \text{ , } \\  < \bar{\psi^{\mathrm{loop-F}} \big( z \big)}     \bar{\psi^{\mathrm{loop-F}} \big( a \big)} >_{\mathscr{I} \times \mathscr{J}}  =     - \bar{f^{\mathrm{loop},\uparrow}_a \big( z \big)} - i \bar{f^{\mathrm{loop},\downarrow}_a \big( z \big) }           \text{ } \text{ . } 
\end{align*}

\noindent \textit{Proof of Theorem} $4^{*}$. Perform similar computations for each of the four identities of the fermion loop operators as given in the previous subsection for \textbf{Theorem} 
 \textit{4}. \boxed{}

\bigskip

\noindent \textbf{Theorem} $5^{*}$ (\textit{Pffaffian from multi-point correlations of the loop fermion operator}). From the loop fermion operator, $\psi^{\mathrm{loop-F}}$,

\begin{align*}
   <  \underset{\mathrm{countably\text{ } many \text{ }} z_i}{\underset{z_i \in \textbf{C}}{\underset{1 \leq i \leq n}{\prod}}} \psi^{(\mathrm{loop-F}),(i)}_{z_i}      >_{\mathscr{I} \times \mathscr{J}} = \underset{1 \leq i , j \leq n}{\bigcup} \mathrm{Pf} \bigg[   <   \psi^{(\mathrm{loop-F}),(i)} \big( z_i \big)  \psi^{(\mathrm{loop-F}),(j)} \big( z_j \big)  >_{\mathscr{I} \times \mathscr{J}}           \bigg]  \text{ } \text{ . } 
\end{align*}

\bigskip

\noindent \textit{Proof of Theorem} $5^{*}$. Repeat the computations provided for \textbf{Theorem} \textit{5} of the previous subsection, but for the loop fermion operator instead of for the Ashkin-Teller fermion operator. \boxed{}

\bigskip

\noindent We also introduce the multipoint observable, and its connections with the Pfaffian (see \textit{4.5} of [13] for a more extensive overview). For the loop model, in comparison to the parafermionic observable that is defined over a single point, over multiple points, the observable takes the form,

\begin{align*}
    F^{\mathrm{MP-loop}} \bigg[ \epsilon , \big( a_1 , \cdots \big) , \big(  z_1 , \cdots \big) , x ,  \big( x_1 , \cdots \big) , \big( \sigma_1 , \cdots  \big)        \bigg] =        \underset{\gamma_i \in \mathcal{C}_{\{a_1 , \cdots , a_n  \}} }{\underset{\gamma_i \subset \Omega}{\underset{\gamma_i : a_i \longrightarrow z_i}{\sum}}}  \big( \alpha^{\mathrm{loop}} \big)^{L ( \gamma_i )} \times \cdots \\  \underset{\gamma_ i \in \mathcal{C}_{\{a_1 , \cdots , a_n  \}}}{\prod}  \mathrm{exp} \big( - i \sigma_i W_{\gamma_i} \big( a_i , z_i \big) \big) x^{l ( \gamma_i ) }                          \text{ } \text{ , } 
\end{align*}

\noindent which is a function of the winding number of each path, with respective parameters $\sigma_1 , \cdots , \sigma_n$, from the collection of paths,

\begin{align*}
  \Gamma \equiv \underset{\gamma_i \in \mathcal{C}_{\{ a_1 , \cdots , a_n \}}}{\bigcup} \big\{ \mathrm{paths}\text{ } \gamma_i \text{ } | \text{ }     0 \leq W_{\gamma_i} \big( a_i , z_i \big) \leq 2 \pi      \big\}    \text{ } \text{ , }
\end{align*}

\noindent for:

\begin{itemize}
\item [$\bullet$] (1): Countably many points $z_1 , \cdots , z_n$ in the complex plane,

\item[$\bullet$] (2): the beginning points, $a_1 , \cdots , a_n$, of each $\gamma_i$,

\item[$\bullet$] (3): the ending points $x_1 , \cdots , x_n$ of each $\gamma_i$,

\item[$\bullet$] (4): the number of loops, $x^{l(\gamma_i)}$ of each $\gamma_i$,

\item[$\bullet$] (5): an arbitrary point $x$ in the complex plane,

\item[$\bullet$] (6): The set of paths, $\mathcal{C}_{\{a_1 , \cdots , a_n  \}}$, where each $\gamma_i$ begins at $a_i$,

\item[$\bullet$] (7): a parameter $\epsilon$, whose form will be provded in \textbf{Theorem} $6^{*}$ below.

\end{itemize}

\noindent \textbf{Theorem} $6^{*}$ (\textit{expectation value of multipoint loop fermion operators is equal to the loop multipoint parafermionic observable}). Denote, for $\psi^{(\mathrm{loop-F})} \equiv \psi$,

\begin{align*}
  \psi^{\uparrow} \big( z \big) = \frac{1}{2} \bigg[ \bar{\psi^{\mathrm{loop-F}}} \big( z \big) - \psi^{\mathrm{loop-F}} \big( z \big) \bigg]  \text{ } \text{ , } 
\end{align*}

\noindent and,

\begin{align*}
    \psi^{\downarrow} \big( z \big)  = \frac{i}{2} \bigg[ {\psi^{\mathrm{loop-F}}} \big( z \big) + \bar{\psi^{\mathrm{loop-F}}} \big( z \big) \bigg]  \text{ } \text{ , } 
\end{align*}

\noindent for $z \in \textbf{C}$, from which one has,

\begin{align*}
  < \underset{\mathrm{countably\text{ }many}\text{ } z_{(j)} \in \textbf{C}}{\underset{1 \leq j \leq 2m -1}{\underset{1 \leq i \leq 2m }{\prod}}} \psi^{\updownarrow , (2m-i)}_{(j)}   \big( z_{(j)} \big)         >^{+}_{\mathscr{I}^{**} \times \mathscr{J}} =      F^{\mathrm{MP-loop}} \bigg[ \epsilon , \big( a_1 , \cdots \big) , \big(  z_1 , \cdots \big) , x ,  \big( x_1 , \cdots \big) , \big( \sigma_1 , \cdots  \big)        \bigg]  \text{ } \text{ , } 
\end{align*}

\noindent for $\updownarrow \in \big\{ \uparrow , \downarrow \big\}$, and,

\begin{align*}
\epsilon   =  \lambda \Longleftrightarrow  \text{ } \updownarrow \equiv \uparrow  \text{ } \text{ , } \\ \epsilon = \frac{1}{\lambda^2} \Longleftrightarrow  \text{ } \updownarrow \equiv \downarrow  \text{ } \text{ . } 
\end{align*}

\noindent \textit{Proof of Theorem} $6^{*}$. To demonstrate that the desired identity holds above, perform the computations along the same lines of those provided in the previous subsection for \textbf{Theorem} \textit{6}. \boxed{}

\subsection{Riemann-Poincare-Steklov operators}

\noindent In the last subsection before the formula for the convolution operator for each model is provided, we analyze RPS operators.

\subsubsection{Ashkin-Teller model}

\noindent The RPS operator satisfies Riemann boundary conditions which were introduced in \textit{2.1} with \textbf{Definition} \textit{3}.

\bigskip

\noindent \textbf{Lemma} \textit{15} (\textit{Riemann boundary conditions for RPS operators}, \textbf{Lemma} \textit{27}, [13]). For some finite volume $\Omega$ over $\textbf{Z}^2$, with edge set $\mathscr{E}$. If there exists an s-holomorphic function, $h : \mathscr{E} \longrightarrow \textbf{C}$, which satisfies Riemann boundary conditions on the boundary of the edge set, then $h \equiv 0$.

\bigskip

\noindent \textit{Proof of Lemma 15}. Refer to the proof of \textbf{Lemma} \textit{27} in [13]. \boxed{}

\bigskip

\noindent \textbf{Lemma} \textit{16} (\textit{s-holomprhic extension of h}, \textbf{Lemma} \textit{28}, [13]). For a boundary edge in $\Omega$, $u$, there exists another boundary edge in $\mathscr{E}$, $v$, such that $u+v$ has an s-holomorphic extension $h : \mathscr{E} \longrightarrow \textbf{C}$, which satisfies Riemann boundary conditions on $\partial \Omega \backslash v$.

\bigskip

\noindent \textit{Proof of Lemma 16}. Refer to the proof of \textbf{Lemma} \textit{28} in [13]. \boxed{}

\bigskip

\noindent With the each \textbf{Lemma} above, we define the RPS operator below.

\bigskip

\noindent \textbf{Definition} \textit{13} (\textit{Ashkin-Teller RPS operator}). The map,

\begin{align*}
        \big( U^{\textbf{b}}_{\Omega}  \big)^{\mathrm{AT}} : \mathcal{R}^{\textbf{b}}_{\Omega} \longrightarrow \mathcal{I}^{\textbf{b}}_{\Omega} \\ u \mapsto v    \text{ } \text{ , } 
\end{align*}

\noindent constitutes an RPS operator, for ${\textbf{b}} \in \partial \Omega$.

\bigskip

\noindent \textbf{Lemma} \textit{17} (\textit{extending the s-holomorphic function to the convolution kernel}, \textbf{Lemma} \textit{30}, [13]). The convolution kernel is given by the summation,

\begin{align*}
  v \big( x \big) = \underset{y \in \textbf{b}}{\sum}     u \big( y \big) f_{\Omega} \big( y , x \big)        \text{ } \text{ , } 
\end{align*}

\noindent which can be extended to an s-holomorphic function $h$, with,

\begin{align*}
    h \big( x \big) = \underset{y \in \textbf{b}}{\sum}  u  \big( y \big) f_{\Omega}  \big( y , z \big)    \text{ } \text{ , } 
\end{align*}

\noindent for $x \in \mathscr{E}$, and $f_{\Omega}$ belonging to the space of functions,

\begin{align*}
 f^{\mathrm{AT}} =  f : \textbf{b} \longrightarrow \textbf{C}  \text{ } \text{ . } 
\end{align*}

\noindent \textit{Proof of Lemma 17}. ($\Leftarrow$) Suppose that the expansion for $h \big( x \big)$ of the s-holomorphic extension exists. To demonstrate that the expansion for $v \big( x \big)$ holds, for edges which lie on the boundary set $\partial \mathscr{E}$, then the function $u \big( x \big) f_{\Omega} \big(x ,x\big)$, and, for $y \in \partial \mathscr{E} \backslash \big\{ x \big\}$, $u \big( x \big) f_{\Omega} \big(y ,x\big)$ each satisfy Riemann boundary conditions. Hence the existence of the s-holomorphic extension implies that another expansion for $v \big( x \big)$ exists. ($\Rightarrow$) Suppose that the expansion for $v \big( x \big)$ exists. Observe that $h \equiv 0$, from which we conclude the argument. \boxed{}

\bigskip

\noindent In the case of finite volumes that are boxes, either squares or rectangles, the RPS operator introduced in \textbf{Definition} \textit{13} can be expressed in terms of a s-holomorphic propagation operator which satisfies the following properties.

\bigskip

\noindent \textbf{Lemma} \textit{18} (\textit{s-holomorphic propagator of the RPS operator}, \textbf{Lemma} \textit{31}, [13]). For a rectangular domain over the square lattice,

\begin{align*}
  \Omega = \textbf{I} \times \big\{ 0 , \cdots , N \big\}   \text{ } \text{ , } 
\end{align*}

\noindent with the bottom endpoint of 
$\Omega$,

\begin{align*}
         \textbf{b} \times \big\{ 0 \big\}           \text{ } \text{ . } 
\end{align*}

\noindent For the Ashkin-Teller propagator,

\begin{align*}
         P^{\mathrm{AT}} : \big( \textbf{R}^2 \big)^{|\textbf{I}^{**}|}  \longrightarrow   \big( \textbf{R}^2 \big)^{|\textbf{I}^{**}|}      \text{ } \text{ , } 
\end{align*}

\noindent there exists a decomposition of $N$ th power of the propagator matrix,

\[  \big( P^{\mathrm{AT}} \big)^N \equiv 
\begin{bmatrix}
   \big( P^{\mathrm{AT}} \big)^N_{\mathscr{R} \mathscr{R}}    &    \big( P^{\mathrm{AT}} \big)^N_{\mathscr{R} \mathscr{S}}    \\   \big( P^{\mathrm{AT}} \big)^N_{\mathscr{S} \mathscr{R}}      &     \big( P^{\mathrm{AT}} \big)^N_{\mathscr{S} \mathscr{S}}   
\end{bmatrix}
\]

\noindent into four blocks, each of dimension $ \big| \textbf{I}^{**} \big| \times \big| \textbf{I}^{**} \big|$. The decomposition of the $N$ th power of the Ashkin-Teller propagator above holds as a result of the decomposition,

\begin{align*}
  \big( \textbf{R}^2 \big)^{\textbf{I}^{**}} \overset{\sim}{=} \big( \textbf{R} \oplus i \textbf{R} \big)^{\textbf{I}^{**}}  \overset{\sim}{=}    \big( \textbf{R} \big)^{\textbf{I}^{**}} \oplus  \big( i \textbf{R} \big)^{\textbf{I}^{**}}  \text{ } \text{ , } 
\end{align*}

\noindent from which it follows that the convolution operator, from the s-holomorphic extension of \textbf{Lemma} \textit{17}, satisfies,

\begin{align*}
     \big( U^{\textbf{b}}_{\Omega}  \big)^{\mathrm{AT}} =             -\big( \big( P^{\mathrm{AT}} \big)^N_{\mathscr{S}\mathscr{S}} \big)^{-1}     \big( P^{\mathrm{AT}} \big)^N_{\mathscr{R}\mathscr{S}}    \text{ } \text{ . } 
\end{align*}

\bigskip

\noindent \textit{Proof of Lemma 18}. From the definition of the Ashkin-Teller convolution operator, for $u \in \mathcal{R}^{\textbf{b}}_{\Omega}$ st $\mathrm{Im} \big( u \big) \equiv 0$,

\[
    \big( P^{\mathrm{AT}} \big)^N \begin{bmatrix}
    u \\ v 
    \end{bmatrix} = \begin{bmatrix}
    w \\ 0 
    \end{bmatrix} 
\]

\noindent for some $w : \textbf{I} \times \big\{ n \big\} \longrightarrow \textbf{R}$, with $\mathrm{Im} \big( w \big) \equiv 0$. From the definition of $\big( P^{\mathrm{AT}} \big)^N$, it follows that $ \big( P^{\mathrm{AT}} \big)^N_{\mathscr{S} \mathscr{R}}   u +  \big( P^{\mathrm{AT}} \big)^N_{\mathscr{S} \mathscr{S}}   v = 0$, while for the remaining term from the equality above involving $w$, it follows that $v= \big(\big( P^{\mathrm{AT}} \big)^N_{\mathscr{S} \mathscr{S}}\big)^{-1} \big( P^{\mathrm{AT}} \big)^N_{\mathscr{R} \mathscr{S}} u$. We conclude the argument. \boxed{}

\subsubsection{Loop model}

\noindent The RPS operator satisfies Riemann boundary conditions which were introduced in \textit{2.1} with \textbf{Definition} \textit{3}.

\bigskip

\noindent \textbf{Lemma} $15^{*}$ (\textit{Riemann boundary conditions for RPS operators}). For some finite volume $\Omega$ over $\textbf{Z}^2$, with edge set $\mathscr{E}$. If there exists an s-holomorphic function, $h : \mathscr{E} \longrightarrow \textbf{C}$, which satisfies Riemann boundary conditions on the boundary of the edge set, then $h \equiv 0$.

\bigskip

\noindent \textit{Proof of Lemma} $15^{*}$. Refer to the proof of \textbf{Lemma} \textit{27} in [13]. \boxed{}

\bigskip

\noindent \textbf{Lemma} $16^{*}$ (\textit{s-holomprhic extension of h}). For a boundary edge in $\Omega$, $u$, there exists another boundary edge in $\mathscr{E}$, $v$, such that $u+v$ has an s-holomorphic extension $h : \mathscr{E} \longrightarrow \textbf{C}$, which satisfies Riemann boundary conditions on $\partial \Omega \backslash v$.

\bigskip

\noindent \textit{Proof of Lemma 16}. Refer to the proof of \textbf{Lemma} \textit{28} in [13]. \boxed{}

\bigskip

\noindent With the each \textbf{Lemma} above, we define the RPS operator below.

\bigskip

\noindent \textbf{Definition} \textit{13} (\textit{loop RPS operator}). The map,

\begin{align*}
        \big( U^{\textbf{b}}_{\Omega}  \big)^{\mathrm{loop}} : \mathcal{R}^{\textbf{b}}_{\Omega} \longrightarrow \mathcal{I}^{\textbf{b}}_{\Omega} \\ u \mapsto v    \text{ } \text{ , } 
\end{align*}

\noindent constitutes an RPS operator, for ${\textbf{b}} \in \partial \Omega$.

\bigskip

\noindent \textbf{Lemma} $17^{*}$ (\textit{extending the s-holomorphic function to the convolution kernel}). The convolution kernel is given by the summation,

\begin{align*}
  v \big( x \big) = \underset{y \in \textbf{b}}{\sum}     u \big( y \big) f_{\Omega} \big( y , x \big)        \text{ } \text{ , } 
\end{align*}

\noindent which can be extended to an s-holomorphic function $h$, with,

\begin{align*}
    h \big( x \big) = \underset{y \in \textbf{b}}{\sum}  u  \big( y \big) f_{\Omega}  \big( y , z \big)    \text{ } \text{ , } 
\end{align*}

\noindent for $x \in \mathscr{E}$, and $f_{\Omega}$ belonging to the space of functions,

\begin{align*}
 f^{\mathrm{loop}} =  f : \textbf{b} \longrightarrow \textbf{C}  \text{ } \text{ . } 
\end{align*}

\noindent \textit{Proof of Lemma} $17^{*}$. The result above follows from identical observations presented for each direction in \textbf{Lemma} \textit{17} in the previous subsection. \boxed{}

\bigskip

\noindent In the case of finite volumes that are boxes, either squares or rectangles, the RPS operator introduced in \textbf{Definition} \textit{13} can be expressed in terms of a s-holomorphic propagation operator which satisfies the following properties.

\bigskip

\noindent \textbf{Lemma} $18^{*}$ (\textit{s-holomorphic propagator of the RPS operator}). For a rectangular domain over the hexagonal lattice,

\begin{align*}
  \Omega = \textbf{I}^{**} \times \big\{ 0 , \cdots , N \big\}   \text{ } \text{ , } 
\end{align*}

\noindent with the bottom endpoint of 
$\Omega$,

\begin{align*}
         \textbf{b} \times \big\{ 0 \big\}           \text{ } \text{ . } 
\end{align*}

\noindent For the loop propagator,

\begin{align*}
         P^{\mathrm{loop}} : \big( \textbf{R}^2 \big)^{|\textbf{I}^{**}|}  \longrightarrow   \big( \textbf{R}^2 \big)^{|\textbf{I}^{**}|}      \text{ } \text{ , } 
\end{align*}

\noindent there exists a decomposition of $N$ th power of the propagator matrix,

\[  \big( P^{\mathrm{loop}} \big)^N \equiv 
\begin{bmatrix}
   \big( P^{\mathrm{loop}} \big)^N_{\mathscr{R} \mathscr{R}}    &    \big( P^{\mathrm{loop}} \big)^N_{\mathscr{R} \mathscr{S}}    \\   \big( P^{\mathrm{loop}} \big)^N_{\mathscr{S} \mathscr{R}}      &     \big( P^{\mathrm{loop}} \big)^N_{\mathscr{S} \mathscr{S}}   
\end{bmatrix}
\]

\noindent into four blocks, each of dimension $ \big| \textbf{I}^{**} \big| \times \big| \textbf{I}^{**} \big|$. The decomposition of the $N$ th power of the loop propagator above holds as a result of the decomposition,

\begin{align*}
  \big( \textbf{R}^2 \big)^{\textbf{I}^{**}} \overset{\sim}{=} \big( \textbf{R} \oplus i \textbf{R} \big)^{\textbf{I}^{**}}  \overset{\sim}{=}    \big( \textbf{R} \big)^{\textbf{I}^{**}} \oplus  \big( i \textbf{R} \big)^{\textbf{I}^{**}}  \text{ } \text{ , } 
\end{align*}

\noindent from which it follows that the convolution operator, from the s-holomorphic extension of \textbf{Lemma} \textit{17}, satisfies,

\begin{align*}
     \big( U^{\textbf{b}}_{\Omega}  \big)^{\mathrm{loop}} =                -\big( \big( P^{\mathrm{loop}} \big)^N_{\mathscr{S}\mathscr{S}} \big)^{-1}     \big( P^{\mathrm{loop}} \big)^N_{\mathscr{R}\mathscr{S}}       \text{ } \text{ . } 
\end{align*}

\bigskip

\noindent \textit{Proof of Lemma} $18^{*}$. Execute the same computations provided for \textbf{Lemma} \textit{18} in the previous subsection, for $ \big( P^{\mathrm{loop}} \big)^N$ instead of for $ \big( P^{\mathrm{AT}} \big)^N$. \boxed{}

\section{Staggered, and odd, eight-vertex models}

\noindent In the last section, we describe connections between the parafermionic observables of the Ashkin-Teller and staggered eight-vertex models, which are shown to be equivalent in [15]. By virtue of this correspondence between parafermionic observables, we list results which are expected to apply, not only for the parafermionic observable of the staggered eight-vertex model, but also for the parafermionic observable of the odd eight-vertex model.

\bigskip

\noindent There exists a mapping between the eight-vertex model and the Ising model, in which there exists a coupling between the two-colored Ashkin-Teller model, with the following mapping between partition functions, [15],

\begin{align*}
 Z \propto \underset{\sigma \in \{ \pm 1 \}}{\sum} \bigg[    \underset{j , k \in \textbf{N}}{\prod}     \mathrm{exp} \bigg[    K^{+} \big( \sigma_{j,k} \sigma_{j+1 , k+1}\big)  + K^{-} \big(  \sigma_{j,k+1} \sigma_{j+1 , k }\big) 
  + \lambda \sigma_{j,k} \big(  \sigma_{j+1,k+1} \sigma_{j,k+1} \sigma_{j+1,k} \big)  \bigg]                               \bigg]             \text{ } \text{ . }
\end{align*}

\noindent At the same critical point along the line of possible critical points for the Ashkin-Teller model for which massive, and massless s-holomorphicity were given (\textbf{Definition} \textit{5} and \textbf{Definition} \textit{6}), we list the results below which would apply to the staggered eight-vertex parafermionic observable. From the similarity which is shared between the parafermionic observable of the Ashkin-Teller model with the parafermionic observable of the staggered eight-vertex model, propagation mechanisms akin to those provided in \textbf{Lemma} \textit{3}, and for \textbf{Lemma} \textit{4}, are satisfied for eqwual couplings of the two Potts models to $\frac{1}{4} \mathrm{log} \big( 3 \big)$. 

\bigskip

\noindent

\subsection{Staggered eight-vertex parafermionic observable}

\noindent Below, we list results emulating those obtained in the previous section for the Ashkin-Teller model.

\bigskip

\noindent \textbf{Proposition} \textit{S8V 1} (\textit{generator relations for the staggered eight-vertex model, from generator relations for the Ising model}, \textbf{Proposition} \textit{8}, [13]). The components of the Ashkin-Teller transfer matrix satisfy the relations,

\begin{align*}
  V^{\mathrm{S8V},h} \equiv     \mathrm{exp} \bigg[   J \big(  \underset{k \in \textbf{I}^{*}}{ \sum}  p^{\mathrm{S8V}}_kq^{\mathrm{S8V}}_k   +  \big( p^{\mathrm{S8V}}_k\big)^{\prime}    \big( q^{\mathrm{S8V}}_k\big)^{\prime}   \big)   \bigg]  +  \mathrm{exp} \bigg[  U \big(  \underset{k \in \textbf{I}^{*}}{ \sum}              p^{\mathrm{S8V}}_k  \big( p^{\mathrm{S8V}}_k \big)^{\prime} q^{\mathrm{S8V}}_k  \big( q^{\mathrm{S8V}}_k \big)^{\prime}   \big)  \bigg]     \text{ } \text{ , } \\    V^{\mathrm{S8V},V} \equiv          \mathscr{P}       \bigg[   \mathrm{exp} \big[   J^{*} \big(  \underset{k \in \textbf{I}}{ \sum}  p^{\mathrm{S8V}}_{k- \frac{1}{2}} q^{\mathrm{S8V}}_{k- \frac{1}{2}}   +  \big( p^{\mathrm{S8V}}_{k- \frac{1}{2}}\big)^{\prime}    \big( q^{\mathrm{S8V}}_{k- \frac{1}{2}} \big)^{\prime}   \big)   \big]  +  \mathrm{exp} \big[  \text{ }  U^{*} \big(  \underset{k \in \textbf{I}}{ \sum}              p^{\mathrm{S8V}}_{k-\frac{1}{2}}  \big( p^{\mathrm{S8V}}_{k-\frac{1}{2}} \big)^{\prime} q^{\mathrm{S8V}}_{k-\frac{1}{2}}  \big( q^{\mathrm{S8V}}_{k- \frac{1}{2}} \big)^{\prime}   \big)  \big]  \bigg]              \text{ } \text{ , } 
\end{align*}

\noindent for generators $p^{\mathrm{S8V}}_k$, $\big( p^{\mathrm{S8V}}_k\big)^{\prime}$, $q^{\mathrm{S8V}}_k$, and $\big( q^{\mathrm{S8V}}_k \big)^{\prime}$, dual couplings $\big(J^{*} , U^{*} \big)$, obtained from $\big(J, U\big)$ under the relation,

\begin{align*}
       \frac{\mathrm{exp} \big( - 2 J + 2 U \big) -1}{\mathrm{exp} \big( - 2 J^{*} + 2U^{*} \big)-1}  = \mathrm{exp} \big( 2 U \big) \mathrm{sinh} \big( 2 J \big) = \frac{1}{\mathrm{exp} \big( 2 U^{*} \big) \mathrm{sinh} \big( 2 J^{*} \big) }    \text{ } \text{ , } 
\end{align*}

\noindent and the prefactor,

\begin{align*}
  \mathscr{P} \equiv    \mathrm{exp} \big[ \big(  \mathrm{exp} \big[ U^{*} , J , J^{*} \big]  - U  \big)  \big( \text{ }              \underset{k \in \textbf{I}}{\sum}          p^{\mathrm{S8V}}_k \big( p^{\mathrm{S8V}}_k \big)^{\prime} q^{\mathrm{S8V}}_k  \big( q^{\mathrm{S8V}}_k  \big)^{\prime }\big)  \big]                \text{ } \text{ . } 
\end{align*}

\bigskip

\noindent \textbf{Lemma} \textit{S8V 1} (\textit{staggered eight-vertex transfer matrix conjugation}, \textbf{Lemma} in \textit{3.2}, [13]). The action,

\begin{align*}
 \big( V^{\mathrm{S8V},h} \big)^{-\frac{1}{2}} \cdot    p^{\mathrm{S8V}}_k       \cdot  \big(   V^{\mathrm{S8V},h}    \big)^{\frac{1}{2}}  =   c p^{\mathrm{S8V}}_k - i s q^{\mathrm{S8V}}_k  \text{ } \text{ , } \\ \big( V^{\mathrm{S8V},h} \big)^{-\frac{1}{2}} \cdot     q^{\mathrm{S8V}}_k      \cdot  \big(   V^{\mathrm{S8V},h}    \big)^{\frac{1}{2}}  =  is p^{\mathrm{S8V}}_k + c q^{\mathrm{S8V}}_k \text{ } \text{ , } 
\end{align*}

\noindent under the composition operator, $\cdot$, $V^{\mathrm{AT},h}$ takes the form above, while the action,

\begin{align*}
 \big( V^{\mathrm{S8V},V} \big)^{-1} \cdot    p^{\mathrm{S8V}}_k       \cdot  \big(   V^{\mathrm{S8V},V}    \big)^{1}  = \frac{C}{S}  p^{\mathrm{S8V}}_k   + \frac{i}{S} q^{\mathrm{S8V}}_{k+1} \text{ } \text{ , } \tag{*} \\ \big( V^{\mathrm{S8V},V} \big)^{-1} \cdot     q^{\mathrm{S8V}}_k      \cdot  \big(   V^{\mathrm{S8V},V}    \big)^{1}  =     - \frac{i}{S}         p^{\mathrm{S8V}}_{k-1} +    \frac{C}{S}          q^{\mathrm{S8V}}_k  \text{ } \text{ , } \tag{**}
\end{align*}

\noindent under the composition operator $\cdot$, $V^{\mathrm{S8V},V}$ takes the form above, where $\textit{(*)}$ holds for $k \neq k_R$, and where $\textit{(**)}$ holds for $k \neq k_L$. Additionally,

\begin{align*}
   \big( V^{\mathrm{S8V},V} \big)^{-1}    \cdot  p^{\mathrm{S8V}}_k \cdot V^{\mathrm{S8V},V}  = p^{\mathrm{S8V}}_k     \text{ } \text{ , } \tag{***} \\   \big( V^{\mathrm{S8V},V} \big)^{-1}  \cdot q^{\mathrm{S8V}}_k 
 \cdot V^{\mathrm{S8V},V} = q^{\mathrm{S8V}}_k   \text{ } \text{ , } \tag{****}
\end{align*}

\bigskip

\noindent where $\textit{(***)}$ holds for $k \equiv k_R$, and where $\textit{(****)}$ holds for $k \equiv k_L$.

\bigskip

\noindent \textbf{Lemma} \textit{S8V 2} (\textit{commutation rule}, \textbf{Lemma} \textit{9}, [13]). For two other maps , the identities,

\begin{align*}
      T_V \cdot R  = R \cdot    T_V    \text{ } \text{ , } \\   T_V \cdot J = J \cdot T_V           \text{ } \text{ , } 
\end{align*}

\noindent for the rotation hold, where the maps $R$ and $J$ satisfy,

\begin{align*}
             R \big( p_k \big) =     i   q^{\mathrm{S8V}}_{a+b-k}                     \text{ } \text{ , } \\  R \big(  q^{\mathrm{S8V}}_k  \big) =   - i p^{\mathrm{S8V}}_{a+b-k}           \text{ } \text{ , } \\  
             J \big( p^{\mathrm{S8V}}_k  \big) =   i p^{\mathrm{S8V}}_k           
             \text{ } \text{ , } \\ J \big( q^{\mathrm{S8V}}_k   \big) =      - i q^{\mathrm{S8V}}_k       \text{ } \text{ , } 
\end{align*}

\noindent for the bases $\psi_x$ and $\bar{\psi_x}$ spanning,

\begin{align*}
   \mathcal{V}_{\psi_x }  = \big\{   x     \text{ } \big| \text{ }  \psi_x \in \mathcal{V}   \big\} \text{ } \text{ , } \\  \mathcal{V}_{\bar{\psi_x}} = \big\{ x    \text{ } \big| \text{ }  \bar{\psi_x} \in \mathcal{V}  \big\}   \text{ } \text{ , } 
\end{align*}

\noindent which have the following images under $R$ and $J$,

\begin{align*}
      R \big( \psi_x \big) =   \bar{\psi_{a+b-k}}  \text{ } \text{ , } \\   R \big( \bar{\psi_x} \big) =   \psi_{b+a-k}    \text{ } \text{ , } \\   J \big( \psi_x \big) =  \bar{\psi_x}   \text{ } \text{ , } \\  J \big( \bar{\psi_x} \big) =   \psi_k   \text{ } \text{ , } 
\end{align*}

\noindent for the elements,

\begin{align*}
   \psi_x \equiv  \frac{i}{\sqrt{2}} \big( p^{\mathrm{S8V}}_k + q^{\mathrm{S8V}}_k \big) \text{ } \text{ , } \\ \bar{\psi_x} \equiv  
 \frac{1}{\sqrt{2}} \big( p^{\mathrm{S8V}}_k - q^{\mathrm{S8V}}_k \big) \text{ } \text{ . } 
\end{align*}

\bigskip

\noindent For exactly similar properties of $T_V$ for the staggered eight-vertex model from those which are satisfied for the Ashkin-Teller model (\textbf{Theorem} in \textit{2.2.1}), the following version of the Wick's formula for the staggered eight-vertex model, which is formulated from the vacuum vector basis, an identical expression for the Wick's formula can be obtained from the Pfaffian. The desired property stems from the fact that a suitable polarization for the staggered eight-vertex model exists from the generators, which are in correspondence with the Ashkin-Teller generators.

\bigskip

\noindent \textbf{Lemma} \textit{S8V 3} One has, for the staggered eight-vertex polarization, that,

\begin{align*}
  \mathcal{W}^{(+)}_{\mathrm{cr}} \equiv \mathrm{span} \big\{    p^{\mathrm{S8V}}_k - i q^{\mathrm{S8V}}_k \text{ } \big| \text{ } k \in \textbf{I}^{*}      \big\}   \text{ } \text{ , } \\ \mathcal{W}^{(+)}_{\mathrm{ann}} \equiv \mathrm{span} \big\{      p^{\mathrm{S8V}}_k + i q^{\mathrm{S8V}}_k \text{ } \big| \text{ } k \in \textbf{I}^{*}      \big\} 
 \text{ } \text{ . } 
\end{align*}

\bigskip

\noindent Besides the expression for correlations obtained with the Pfaffian, by passing to a suitable polarization, the lattice fermion operator for the staggered eight-vertex model, as does that for the Ashkin-Teller model, satisfy the properties of massive, and massless, s-holomorphicity (coinciding with those provided in \textbf{Definition} \textit{5} and in \textbf{Definition} \textit{6}).

\bigskip

\noindent For the multipoint staggered eight-vertex parafermionic observable, from results akin to \textbf{Theorem} \textit{6}, the multpoint parafermionic observable for the staggered eight-vertex model admits an equivalent identity.

\bigskip

\noindent The kernel of the convoluation operator for the staggered eight-vertex model

\bigskip

\noindent \textbf{Lemma} \textit{S8V 4}. (\textit{extending the s-holomorphic function to the convolution kernel}). The convolution kernel is given by the summation,

\begin{align*}
  v \big( x \big) = \underset{y \in \textbf{b}}{\sum}     u \big( y \big) f_{\Omega} \big( y , x \big)        \text{ } \text{ , } 
\end{align*}

\noindent which can be extended to an s-holomorphic function $h$, with,

\begin{align*}
    h \big( x \big) = \underset{y \in \textbf{b}}{\sum}  u  \big( y \big) f_{\Omega}  \big( y , z \big)    \text{ } \text{ , } 
\end{align*}

\noindent for $x \in \mathscr{E}$, and $f_{\Omega}$ belonging to the space of functions,

\begin{align*}
 f^{\mathrm{S8V}} = f : \textbf{b} \longrightarrow \textbf{C}  \text{ } \text{ . } 
\end{align*}

\bigskip

\noindent The final result provides a block decomposition of the $N$ th power of the staggered eight-vertex model propagator.

\subsection{Odd eight-vertex parafermionic observable, and abelian sandpile parafermionic observable}

\noindent Below, instead of listing results emulating those obtained in the previous section for the Ashkin-Teller, and staggered eight-vertex, models, we refer to the previous section. To avoid being overly repetitive, all results for the staggered eight-vertex model directly carry over for the odd eight-vertex model.

\bigskip

\noindent Furthermore, as described in [15], the fact that the abelian sandpiple model coincides with the $0$-color limit of the self dual Potts model directly implies that the operator formalism can also apply to that model, in which the observable is a boson. Again, for the sake of not being overly repetitive, the arguments for the operator formalism which are expected to carry over are excluded.

\section{References}

\noindent [1] Alam, I.T., Batchelor, M.T. Integrability as a consequence of discrete holomorphicity: the $Z_N$ model. \textit{J. Phys. A: Math. Theor.} 45 494014 (2012).

\bigskip

\noindent [2] Beffara, V., Duminil-Copin, H. Smirnov's Fermionic Observable Away from Criticality. \textit{Annals of Probability} \textbf{40}:6, 2667-2689 (2012). \textit{doi}: 10.1214/11-A0P689

\bigskip

\noindent [3] Blöte, H.W.J., Nienhuis, B. The phase diagram of the $\mathrm{O} \big( n \big)$ model. \textit{Physica A: Statistical Mechanics and its Applications} \textbf{160}(2): 121-134 (1989).

\bigskip

\noindent [4] Cimasoni, D., and Duminil-Copin, H. The Critical Temperature for the Ising Model on Planar Doubly Periodic Graphs. \textit{Electron. J. Probab.} \textbf{18}: 1-18. \textit{doi}: 10.1214/EJP.v18-2352 (2013).

\bigskip

\noindent [5] Delfino, G., Grinza, P. Universal ratios along a line of critical points. The Ashkin-Teller model. \textit{Nuclear Physics B} \textbf{682}(3): 521-550 (2004).

\bigskip

\noindent [6] Duminil-Copin, H. Divergence of the correlation length for critical planar FK percolation with $1 \leq q \leq 4$ via parafermionic observables. \textit{Journal of Physics A Mathematical and Theoretical} \textbf{45}: 49. \textit{doi}: 10.1088/175108113/45/45/494013 (2012).

\bigskip

\noindent [7] Duminil-Copin, H., Van Enter, A.C.D. Sharp Metastability Threshold for an Anisotropic Bootstrap Percolation Model. \textit{Ann. Probab.} \textbf{41}(3A): 1218-1242 (2013). \textit{doi}: 10.1214/11-AOP722

\bigskip

\noindent [8] Duminil-Copin, H., Garban, C., Pete, G. The near-critical FK-Ising model. \textit{Commun. Math. Phys.} 326, 1–35 (2014). https://doi.org/10.1007/s00220-013-1857-0

\bigskip

\noindent [9] Duminil-Copin, H., Kozma, G., Yadin, A. Supercritical self-avoiding walks are space-filling. \textit{Ann. Inst. H. Poincaré Probab. Statist.} \textbf{50}(2): 315-326 (May 2014). \textit{doi}: 10.1214/12-AIHP528.

\bigskip

\noindent [10] Duminil-Copin, H., Hongler, C., Nolin, P. Connection probabilities and RSW-type bounds for the FK Ising model. \textit{Communications on Pure and Applied Mathematics} \textbf{64}: 9, 1165-1198 (2011).

\bigskip

\noindent [11] Duminil-Copin, H., Smirnov, S. Conformal Invariance of Lattice Models. \textit{Clay Mathematics Proceedings} \textbf{15} (2012).

\bigskip

\noindent [12] Duminil-Copin, H., Smirnov, S. The connective constant of the honeycomb lattice equals $\sqrt{2+\sqrt{2}}$. \textit{Annals of Mathematics} \textbf{175}:3, 1653-1665 (2012).

\bigskip

\noindent [13] Hongler, C., Kytola, K., Zahabi, A. Discrete Holomorphicity and Ising Model Operator Formalism. \textit{AMS, Contemporary Mathematics} \textbf{644}: 79-116 (2013).

\bigskip

\noindent [14] Huhtala, P. Loop $\mathrm{O} \big( n \big)$ models: a numerical transfer matrix study of long range order and s-holomorphicity. \textit{Master's Thesis in Mathematics and Operations Research, Aalto University} (2021).

\bigskip

\noindent [15] Tanhayi-Ahari, M., Rouhani, S. Discrete Holomorphic Parafermions in the Eight Vertex Model. \textit{arXiv: 1209.4253} (2012).

\bigskip

\noindent [16] Wu, F.Y., Kunz, H. The odd eight-vertex model. \textit{J. Stat. Phys.} \textbf{116}: 67-78 (2004).

\end{document}